

  \def \ifundef #1{\expandafter \ifx \csname #1\endcsname \relax }


  \input amssym
  \input miniltx
  \input pictex


  \font \bbfive = bbm5
  \font \bbeight = bbm8
  \font \bbten = bbm10
  \font \rs = rsfs10 \font \rssmall = rsfs10 scaled 833 \font \rstiny = rsfs10 scaled 600  
  \font \eightbf = cmbx8
  \font \eighti = cmmi8 \skewchar \eighti = '177
  \font \fouri = cmmi5 scaled 800 
  \font \eightit = cmti8
  \font \eightrm = cmr8
  \font \eightsl = cmsl8
  \font \eightsy = cmsy8 \skewchar \eightsy = '60
  \font \eighttt = cmtt8 \hyphenchar \eighttt = -1

  \font \sixi = cmmi6 \skewchar \sixi = '177
  \font \sixrm = cmr6
  \font \sixsy = cmsy6 \skewchar \sixsy = '60
  \font \tensc = cmcsc10

  \scriptfont \bffam = \bbeight
  \scriptscriptfont \bffam = \bbfive
  \textfont \bffam = \bbten

  \newskip \ttglue

  \def \eightpoint {\def \rm {\fam 0 \eightrm }\relax
  \textfont 0 = \eightrm \scriptfont 0 = \sixrm \scriptscriptfont 0 = \fiverm
  \textfont 1 = \eighti \scriptfont 1 = \sixi \scriptscriptfont 1 = \fouri
  \textfont 2 = \eightsy \scriptfont 2 = \sixsy \scriptscriptfont 2 = \fivesy
  \textfont 3 = \tenex \scriptfont 3 = \tenex \scriptscriptfont 3 = \tenex
  \def \it {\fam \itfam \eightit }\relax
  \textfont \itfam = \eightit
  \def \sl {\fam \slfam \eightsl }\relax
  \textfont \slfam = \eightsl
  \def \bf {\fam \bffam \eightbf }\relax
  \textfont \bffam = \bbeight \scriptfont \bffam = \bbfive \scriptscriptfont \bffam = \bbfive
  \def \tt {\fam \ttfam \eighttt }\relax
  \textfont \ttfam = \eighttt
  \tt \ttglue = .5em plus.25em minus.15em
  \normalbaselineskip = 9pt
  \def \MF {{\manual opqr}\-{\manual stuq}}\relax
  \let \sc = \sixrm
  \let \big = \eightbig
  \let \rs = \rssmall
  \setbox \strutbox = \hbox {\vrule height7pt depth2pt width0pt}\relax
  \normalbaselines \rm }


  \def \text #1{\mathchoice {\hbox {#1}} {\hbox {#1}} {\hbox {\eightrm #1}} {\hbox {\sixrm #1}}}
  \def \varbox #1{\setbox 0\hbox {$#1$}\setbox 1\hbox {$I$}{\ifdim \ht 0< \ht 1 \scriptstyle #1 \else \scriptscriptstyle #1 \fi }}

  \def \rsbox #1{{\mathchoice {\hbox {\rs #1}} {\hbox {\rs #1}} {\hbox {\rssmall #1}} {\hbox {\rstiny #1}}}}
  \def \mathscr #1{\rsbox {#1}}


  \def \myhangindent {1.1truecm}
  \def \VarItem #1#2{\smallskip \par \noindent \kern #1 \hangindent #1 \llap {#2\enspace }\ignorespaces }
  \def \Item #1{\VarItem {\myhangindent }{#1}}

  \newcount \zitemno \zitemno = 0
  \def \iStyle #1{\zitemno =0\relax \def \itemStyle {#1}}  
  \iStyle i 

  \def \roman #1{\ifcase #1 0\or i\or ii\or iii\or iv\or v\or vi\or vii\or viii\or ix\or x\or xi\or xii\or xiii\or xiv\or
xv\or xvi\or xvii\or xviii\or xix\or xx\else zzz\fi }

  \def \latin #1{\ifcase #1 0\or a\or b\or c\or d\or e\or f\or g\or h\or i\or j\or k\or l\or m\or n\or o\or p\or q\or r\or
s\or t\or u\or v\or w\or x\or y\or z\else zzz\fi }

  \def \printItem {\if \itemStyle i\roman \zitemno \else
                   \if \itemStyle a\latin \zitemno \else
                   \if \itemStyle 1\number \zitemno \else
                   \if \itemStyle .$\bullet $\fi \fi \fi \fi }

  \def \item {\global \advance \zitemno by 1 \Item {{\rm (\printItem )}}}


  \newcount \secno \secno = 0
  \newcount \stno \stno = 0

  \def \advseqnumbering {\global \advance \stno by 1}

  \def \section #1 \par {\global \advance \secno by 1 \stno = 0
    \goodbreak \bigbreak
    \noindent {\bf \number \secno .\enspace #1.}
    \edef \rightRunningHead {\ifundef {authorRunninhHead}#1\else \authorRunninhHead \fi }
    \Headlines {#1}{\rightRunningHead }
    \nobreak \medskip \noindent }

  \def \iP {N}

  \def \state #1 #2\par {\begingroup
    \iStyle i
    \global \zitemno = 0\relax
    \def \iP {Y}
    \medbreak \noindent \advseqnumbering {\bf \number \secno .\number \stno .\enspace #1.\enspace \sl #2\par }\medbreak \endgroup }

  \def \Proof {
    \global \def \iP {Y}
    \global \zitemno = 0\relax
    \medbreak \noindent {\it Proof.\enspace }}

  \def \endProof {\if \iP Y \hfill $\endproofmarker $ \looseness = -1 \fi
    \medbreak
    \global \def \iP {N}}

  \def \closeProof {\eqno \endproofmarker
    \global \def \iP {E}}

  \def \undefrule {\kern 2pt \vrule width 2pt height 5pt depth 0pt \kern 2pt}
  \def \UndefLabels {}
  \def \possundef #1{\ifundef {#1}\undefrule {\eighttt #1}\undefrule
    \global \edef \UndefLabels {\UndefLabels #1\par }
  \else \csname #1\endcsname \fi }

  \def \define #1#2{\global \expandafter \edef \csname #1\endcsname {#2}}


  \newcount \bibno \bibno = 0
  \def \newbib #1#2{\define {#1}{#2}}

  \def \Bibitem #1 #2; #3; #4 \par {\global \advance \bibno by 1
    \VarItem {0.5truecm}{[\possundef {#1}]} #2, {``#3''}, #4.\par
    \ifundef {#1PrimarilyDefined}\else
      \fatal {Duplicate definition for bibliography item ``{\tt #1}'',
      already defined in ``{\tt [\csname #1\endcsname ]}''.}
      \fi \ifundef {#1}\else \edef \prevNum {\csname #1\endcsname } \ifnum \bibno =\prevNum \else \error {Mismatch
        bibliography item ``{\tt #1}'', defined earlier (in aux file ?) as ``{\tt \prevNum }'' but should be ``{\tt
        \number \bibno }''.  Running again should fix this.}  \fi \fi
    \define {#1PrimarilyDefined}{#2}\relax
    }

  \def \jrn #1, #2 (#3), #4-#5;{{\sl #1}, {\bf #2} (#3), #4--#5}
  \def \Article #1 #2; #3; #4 \par {\Bibitem #1 #2; #3; \jrn #4; \par }

  \def \references {\begingroup \bigbreak \eightpoint \centerline {\tensc References} \nobreak \medskip \frenchspacing }



  %
  \catcode `@ = 11
  \def \@Nnil {\@Nil }%
  \def \@Fornoop #1\@@ #2#3{}%
  \def \For #1:=#2\do #3{%
     \edef \@Fortmp {#2}%
     \ifx \@Fortmp \empty \else
        \expandafter \@Forloop #2,\@Nil ,\@Nil \@@ #1{#3}%
     \fi
  }%
  \def \@Forloop #1,#2,#3\@@ #4#5{\@Fordef #1\@@ #4\ifx #4\@Nnil \else
         #5\@Fordef #2\@@ #4\ifx #4\@Nnil \else #5\@iForloop #3\@@ #4{#5}\fi \fi
  }%
  \def \@iForloop #1,#2\@@ #3#4{\@Fordef #1\@@ #3\ifx #3\@Nnil
         \let \@Nextwhile =\@Fornoop \else
        #4\relax \let \@Nextwhile =\@iForloop \fi \@Nextwhile #2\@@ #3{#4}%
  }%
  \def \@Forspc { }%
  \def \@Fordef {\futurelet \@Fortmp \@@Fordef }
  \def \@@Fordef {%
    \expandafter \ifx \@Forspc \@Fortmp 
      \expandafter \@Fortrim
    \else
      \expandafter \@@@Fordef
    \fi
  }%
  \expandafter \def \expandafter \@Fortrim \@Forspc #1\@@ {\@Fordef #1\@@ }%
  \def \@@@Fordef #1\@@ #2{\def #2{#1}}%
  %
  %
  \def \citetrk #1{{{\def \a {}\For \name :=#1\do {\a {\bf \possundef {\name }}\def \a {, }}}}} 
  \def \c@ite #1{{\rm [\citetrk {#1}]}}
  \def \sc@ite [#1]#2{{\rm [\citetrk {#2}\hskip 0.7pt:\hskip 2pt #1]}}
  \def \du@lcite {\if \pe@k [\expandafter \sc@ite \else \expandafter \c@ite \fi }
  \def \cite {\futurelet \pe@k \du@lcite }
  \catcode `\@ =12
  %


  \def \Headlines #1#2{\nopagenumbers
    \headline {\ifnum \pageno = 1 \hfil
    \else \ifodd \pageno \tensc \hfil \lcase {#1} \hfil \folio
    \else \tensc \folio \hfil \lcase {#2} \hfil
    \fi \fi }}

  \font \titlefont = cmbx12
  \long \def \title #1{\begingroup
    \titlefont
    \parindent = 0pt
    \baselineskip = 16pt
    \leftskip = 35truept plus 1fill
    \rightskip = \leftskip
    #1\par \endgroup }

  \long \def \Quote #1\endQuote {\begingroup \leftskip 35truept \rightskip 35truept \parindent 17pt \eightpoint #1\par \endgroup }
  \long \def \Abstract #1\endAbstract {\vskip 1cm \Quote \noindent #1\endQuote }
  \def \Address #1#2{\bigskip {\tensc #1 \par \it E-mail address: \tt #2}}
  \def \Authors #1{\bigskip \centerline {\tensc #1}}
  \def \Note #1{\footnote {}{\eightpoint #1}}
  \def \Date #1 {\Note {\it Date: #1.}}


  \def \Case #1:{\medskip \noindent {\tensc Case #1:}}

  \def \fix {\smallskip \noindent $\blacktriangleright $\kern 12pt}

  \def \lcase #1{\edef \auxvar {\lowercase {#1}}\auxvar }
  \def \emph #1{{\it #1}\/}
  \def \reqno #1{\if \iP N \advseqnumbering \fi \eqno {#1}}
  \def \explica #1#2{\mathrel {\buildrel \hbox {\sixrm #1} \over #2}}
  \def \pilar #1{\vrule height #1 width 0pt}

  \newcount \fnctr \fnctr = 0
  \def \fn #1{\global \advance \fnctr by 1
    \edef \footnumb {$^{\number \fnctr }$}\relax
    \footnote {\footnumb }{\eightpoint #1\par \vskip -10pt}}


  \def \mathbb #1{{\bf #1}}
  \def \frac #1#2{{#1\over #2}}
  \def \<{\left \langle \vrule width 0pt depth 0pt height 8pt }
  \def \>{\right \rangle }
  \def \({\left (}
  \def \){\right )}
  \def \ds {\displaystyle }
  \def \and {\mathchoice {\hbox {\quad and \quad }} {\hbox { and }} {\hbox { and }} {\hbox { and }}}

  \def \imply {\mathrel {\Rightarrow }}
  \def \IFF {\kern 7pt\Leftrightarrow \kern 7pt}
  \def \IMPLY {\kern 7pt \Rightarrow \kern 7pt}
  \def \for #1{\mathchoice {\quad \forall \,#1} {\hbox { for all } #1} {\forall #1}{\forall #1}}
  \def \endproofmarker {\square }
  \def \"#1{{\it #1}\/} 
  \def \inv {^{-1}}
  \def \*{\otimes }
  \def \caldef #1{\global \expandafter \edef \csname #1\endcsname {{\cal #1}}}
  \def \mathcal #1{{\cal #1}}
  \def \bfdef #1{\global \expandafter \edef \csname #1\endcsname {{\bf #1}}}
  \bfdef N \bfdef Z \bfdef C \bfdef R
  \def \exists {\mathchar "0239\kern 1pt }
  \def \labelarrow #1{\setbox 0\hbox {\ \ $#1$\ \ }\ {\buildrel \textstyle #1 \over {\hbox to \wd 0 {\rightarrowfill }}}\ }
  \def \subProof #1{\medskip \noindent #1\enspace }
  \def \itmProof (#1) {\subProof {(#1)}}
  \def \itemImply #1#2{\subProof {#1$\Rightarrow $#2}}
  \def \itmImply (#1) > (#2) {\itemImply {(#1)}{(#2)}}
  
  \ifundef {varnothing}  \fi

  \def \close {\end }

  \def \medcup {\mathop {\mathchoice {\raise 1pt \hbox {$\mathchar "1353$}}{\mathchar "1353}{\mathchar "1353}{\mathchar "1353}}}
  \def \medcap {\mathop {\mathchoice {\raise 1pt \hbox {$\mathchar "1354$}}{\mathchar "1354}{\mathchar "1354}{\mathchar "1354}}}
  
  \def \medotimes {\mathop {\mathchoice {\hbox {$\mathchar "134E$}}{\mathchar "134E}{\mathchar "134E}{\mathchar "134E}}}

  \def \clauses #1{\left \{\matrix {#1}\right .}
  \def \cl #1#2{\hfill #1\hfill ,&#2\hfill \vrule width 0pt height 12pt depth 6pt\cr }

\newbib {AraBosa}{1} \newbib {BeckyOne}{2} \newbib {BeckyTwo}{3} \newbib {LisaOne}{4} \newbib {Demeneghi}{5} \newbib {DokuchaExel}{6} \newbib {EffrosHahnOne}{7} \newbib {EffrosHahnTwo}{8} \newbib {Tight}{9} \newbib {ExelPardo}{10} \newbib {ExelPitts}{11} \newbib {FS}{12} \newbib {FeldmanMoore}{13} \newbib {GR}{14} \newbib {IonWill}{15} \newbib {Jacobson}{16} \newbib {Keimel}{17} \newbib {Kumjian}{18} \newbib {Lawson}{19} \newbib {Nystedt}{20} \newbib {Paterson}{21} \newbib {Renault}{22} \newbib {RenaultStruId}{23} \newbib {RenaultCartan}{24} \newbib {Sauvageot}{25} \newbib {SimsWill}{26} \newbib {SteinbergOne}{27} \newbib {SteinbergDisintegr}{28} \newbib {SteinbergEffrosHahn}{29}

 \def \fld {{\bf K}} \def \ex #1{^{\scriptscriptstyle (#1)}} \def \S {{\cal S}}

\def \dom {\text {dom}} \def \ran {\text {ran}} \def \src {\text {src}} \def \tgt {\text {tgt}}

\def \supp {\text {supp}} \def \extra {\star } \def \Ak {A_{\fld }(G,E)} \def \Lc {L_c(G\ex 0)}

\def \st {^\star } 

\def \limitSize #1{\setbox 0\hbox {$#1$}\setbox 1\hbox {$x$} \ifdim \ht 0 > \ht 1 {\scriptscriptstyle #1} \else {\scriptstyle #1} \fi }

\def \Suber #1#2#3{ \def \a {q} \def \b {#1} \edef \spc {\if \a \b \kern 0pt \else \kern -1pt\fi } #1_{\spc \limitSize {#3}}^{\limitSize {#2}}}

\def \SimpleSuber #1#2#3{ \if #2#3 \mathchoice {#1\kern -0.5pt \raise 0.5pt \hbox {$\scriptstyle (#2)$}} {#1\kern -0.5pt \raise 0.5pt \hbox {$\scriptstyle (#2)$}} {#1(#2)} {#1(#2)} \else \Suber {#1}{#2}{#3} \fi }

\def \C #1#2{\Suber C{#1}{#2}} \def \H #1#2{\Suber H{#1}{#2}} \def \B #1#2{\Suber B{#1}{#2}} \def \L #1#2{\Suber L{#1}{#2}} \def \E #1#2{\Suber E{#1}{#2}} \def \N #1#2{\Suber N{#1}{#2}} \def \p #1#2{\Suber p{#1}{#2}} \def \q #1#2{\Suber q{#1}{#2}} \def \G #1#2{\SimpleSuber G{#1}{#2}} \def \M #1#2{\Suber M{#1}{#2}}

\def \ev #1#2{\kern 0.5pt\left \langle #1, #2\right \rangle \kern 0.5pt} \def \evx #1{\ev {#1}x} \def \evy #1{\ev {#1}y}

\def \1#1{1_{\scriptscriptstyle #1}} \def \supp {\text {supp}} \def \Orb #1{\text {Orb}(#1)} \def \span {\text {span}}

\newdimen \htzero

\def \paren #1{\setbox 0\hbox {$#1$} \htzero =\ht 0 \setbox 1\hbox {$($} \ifdim \htzero < \ht 1 (#1) \else \setbox 1\hbox {$\big ($} \ifdim \htzero < \ht 1 \big (#1\big ) \else \left (#1\right ) \fi \fi }

\def \FunctorInd #1{\text {Ind}_{#1}} \def \Ind #1#2{\FunctorInd {#1} \paren {#2}} \def \Indx #1{\Ind x{#1}} \def \Resx #1{\text {Res}_x\paren {#1}} \def \Ann #1#2{\text {Ann}{\raise -3pt\hbox {$\scriptstyle #1$}}\paren {#2}}%

\def \tprod #1{\mathchoice {\mathop \otimes \limits _{\scriptscriptstyle \B {#1}{#1}}} {\otimes _{\scriptscriptstyle \B {#1}{#1}}} {\otimes _{\scriptscriptstyle \B {#1}{#1}}} {\otimes _{\scriptscriptstyle \B {#1}{#1}}} } \def \tprodx {\tprod x} \def \Nbd {\rsbox {N}_{x}} \def \limu {\lim _{\scriptscriptstyle U\downarrow 0}}

\def \V #1{V_{\scriptscriptstyle [#1]}} \def \Vx {\V x}

\title {TWISTED STEINBERG ALGEBRAS,\break REGULAR INCLUSIONS\break AND INDUCTION}

\medskip \centerline {\tensc Dedicated to the memory of Fernando Abadie}

\def \authorRunninhHead {M. Dokuchaev, R. Exel, and H. Pinedo}

\Authors {\authorRunninhHead }

\Note {Partially supported by CNPq and FAPESP.}

\Abstract Given a field $\fld $ and an ample (not necessarily Hausdorff) groupoid $G$, we define the concept of a \emph {line bundle} over $G$ inspired by the well known concept from the theory of C*-algebras. If $E$ is such a line bundle, we construct the associated \emph {twisted Steinberg algebra} in terms of sections of $E$, extending the original construction introduced independently by Steinberg in 2010, and by Clark, Farthing, Sims and Tomforde in a 2014 paper (originally announced in 2011). We also generalize (strictly, in the non-Hausdorff case) the 2023 construction of (cocycle) twisted Steinberg algebras of Armstrong, Clark, Courtney, Lin, Mccormick and Ramagge. We then extend Steinberg's theory of induction of modules, not only to the twisted case, but to the much more general case of \emph {regular inclusions} of algebras. Among our main results, we show that, under appropriate conditions, every irreducible module is induced by an irreducible module over a certain abstractly defined \emph {isotropy algebra}. We also describe a process of \emph {disintegration} of modules and use it to prove a version of the Effros-Hahn conjecture, showing that every primitive ideal coincides with the annihilator of a module induced from isotropy. \endAbstract

\section Introduction

Steinberg algebras were introduced independently by Steinberg in \cite {SteinbergOne}, and by Clark, Farthing, Sims and Tomforde in \cite {LisaOne}, as a purely algebraic counterpart of Renault's \cite {Renault} groupoid C*-algebras.

Steinberg's work in \cite {SteinbergOne} was partly motivated by a desire to describe irreducible representations of inverse semigroups, which in turn led to a wide ranging result classifying \emph {spectral} simple modules over certain algebras associated to ample groupoids, precisely those which later came to be called Steinberg algebras. In \cite [Theorem 7.26]{SteinbergOne} Steinberg shows that all such modules are \emph {induced} from representations of certain \emph {isotropy} groups.

In \cite {BeckyOne} Armstrong, Clark, Courtney, Lin, McCormick and Ramagge introduced a natural generalization of Steinberg algebras, in terms of a \emph {twisted}, ample, Hausdorff groupoid, namely an ample groupoid equipped with a locally constant 2-cocycle which in turn is used to \emph {twist} the standard convolution product on the Steinberg algebra. Many important known results for Steinberg algebras were generalized in \cite {BeckyOne} to the twisted case.

The theory of locally compact (not necessarily ample) groupoids is also awash with twisted groupoids, but here the notion of \emph {twist} is often given in terms of certain groupoid extensions (see e.g.~\cite [Section 4]{RenaultCartan}). As it turns out, 2-cocycles provide examples of groupoid extensions, but these do not exhaust all possible examples due to the fact that the most general groupoid extension may be loosely viewed as a combination of an algebraic twist (a 2-cocycle) with a topological twist (the possible non-splitting of the given extension).

However, over ample, Hausdorff groupoids, the topological part of every twisting is always trivial, as noted in \cite [Theorem 4.10]{BeckyOne}, so this justifies the choice of restricting the definition of twisted Steinberg algebras to algebraically twisted groupoids, i.e.~2-cocycles, provided the applications are targeted at Hausdorff groupoids only.

The result mentioned in the above paragraph is in fact referred to as folklore in \cite {BeckyOne}, perhaps due to the analogy with the fact that, over a totally disconnected, compact, Hausdorff space, every locally trivial bundle (be it a vector bundle or a principal bundle) is in fact trivial since it is possible to patch local sections into a global section.

Speaking of the Hausdorff property of groupoids, it should be mentioned that many ample groupoids naturally occurring in applications fail to be Hausdorff (see e.g. \cite [Corollary 4.11]{ExelPardo}), so the case could be made that studying non Hausdorff groupoids is a worthy endeavor. Indeed Steinberg's work in \cite {SteinbergOne} is based on ample groupoids which are not assumed to be Hausdorff.

Among the motivations of the present work is a desire to extend Steinberg's \cite {SteinbergOne} work on induced representations to the case of non-Hausdorff, twisted, ample groupoids. In the process of pinning down the precise category of groupoids to focus on, we have realized that the failure of the Hausdorff property breaks down the process of patching local sections mentioned above and in fact it surprisingly allows for nontrivial topological twists, as shown by the example discussed in (3.7), below.

We thus develop a theory of (topologically) twisted groupoids from scratch, although we prefer to work with the equivalent notion of line bundles (Definition (2.4)) rather than groupoid extensions. For such twisted groupoids we then introduce a generalized twisted Steinberg algebra.

Incidentally we would like to point out that, despite its greater generality, the line bundle point of view considerably simplifies some arguments, such as e.g.~the somewhat sticky proof of \cite [Proposition 3.2]{BeckyOne}. Indeed the line bundle point of view makes the verification that the convolution product is well defined a trivial task.

Still speaking of our search for the appropriate category to study induced representations, we have further realized that our generalized algebras give rise to regular inclusions of algebras (Definition (5.5)), in the spirit of \cite {AraBosa} and \cite {BeckyTwo}. In fact, regular inclusions of algebras are closely related to groupoid algebras, not only in pure algebra, but also in the theory of von Neumann and C*-algebras \cite {FeldmanMoore, Kumjian, Renault, ExelPitts}, in the sense that the literature contains an ever expanding list of results where conditions are given for a regular inclusion of algebras to be given in terms of a groupoid.

However, without any of those extra conditions, the concept of a regular inclusion is substantially more general. Nevertheless, as we pushed in the direction of studying induced representations in this more general setting, we were startled to realize that, under our main set of standing hypotheses (see (5.6)), the techniques called for by the study of our situation turned out to be incredibly powerfull, allowing us to prove most, if not all, of the results we had hoped to obtain. In particular we have been able to develop a complete theory of induction of representations in the apparently bare context of regular inclusions, involving the identification of isotropy groups, imprimitivity bimodules, as well as the processes of restriction and induction of modules and, to a certain extent, also the process of disintegration of modules.

As a result, this paper is much more about regular inclusions than twisted Steinberg algebras, although, as already indicated, all of our results duly apply to the latter class of algebras.

The class of regular inclusions studied here is precisely described in (5.6), as already mentioned, but it can be equivalently, and very concisely, described as follows: take any associative algebra $B$ containing a linearly spanning, multiplicative subsemigroup $S$, such that $S$ is an inverse semigroup in itself. The linear span of the idempotent semilattice $E(S)$ is then an abelian subalgebra of $B$, and the pair $(A, B)$ forms a regular inclusion in our class. Conversely, any regular inclusion in our class may be described in terms of an associative algebra $B$ spanned by an inverse semigroup, as above. See Proposition (5.9) for more details on this.

The suggestion of working with regular inclusions, rather than Steinberg algebras, came from \cite [Chapter 2]{ExelPitts}, were regular inclusions are studied in the category of C*-algebras, and where a generalized concept of \emph {isotropy algebra} \cite [Definition 2.1.4]{ExelPitts} was first introduced. We in fact borrow extensively from there, although this presents a bit of a dilemma, often faced when one straddles the purely algebraic and the C*-algebraic aspects of the same fundamental ideas: given the differences between these fields, it would not be technically correct to refer to the C*-algebraic results of \cite {ExelPitts} when the object of study here is an associative algebra over an arbitrary field. However, the results, and often the proofs, are very similar to each other and one cannot honestly claim great originality when reproving these C*-algebraic results in a purely algebraic context. Faced with this dilemma we have opted for the perhaps more conservative alternative of offering complete proofs, which sometimes present moderately interesting nuances. The adaptation of the results from \cite {ExelPitts} to our algebraic context are to be found in section (5), which readers acquainted with \cite {ExelPitts}, and well trained in the interplay between pure algebra and C*-algebra, might prefer to skip.

With section (6) we start the development of the theory of induced representations, where we begin by introducing the imprimitivity bimodule in the context of regular inclusions, which is the key ingredient in the process of inducing modules from isotropy algebras via the familiar tensor product construction. The dual notion of restriction is then described and we prove the expected result (8.4) according to which restricting an induced module produces the original module used in the induction process.

A lot more interesting and useful is the reverse procedure of inducing a restricted module, which produces a faithfully represented submodule of the original module, as shown in (10.1), and in turn implies one of our main results, namely (10.2), showing that, under appropriate conditions, an irreducible module is necessarily induced from an irreducible module over an isotropy algebra.

Starting with section (11) we shift our focus from modules to ideals. Observing that every (two sided) ideal in an s-unital algebra is the annihilator of a module, the plan is to apply our theory of induction of modules to obtain results about ideals. We thus focus on a well kown problem in the literature, often referred to as the Effros-Hahn conjecture. Initially posed in \cite {EffrosHahnOne} (see also \cite {EffrosHahnTwo}), this has a long history (for a non-comprehensive list see \cite {Sauvageot, GR, FS, RenaultStruId, IonWill, SimsWill, DokuchaExel, Demeneghi, SteinbergEffrosHahn}) and it can be roughly described as asking for ideals in an algebra to be obtained in terms of ideals induced from isotropy (evidently within a context where the notion of isotropy makes sense). Theorem (12.14), our main result along this line, shows that, if ``$A\subseteq B$'' is a regular inclusion satisfying our main standing hypotheses (5.6), then every ideal $I\trianglelefteq B$ can be written as the intersection of ideals induced from isotropy. Moreover, if $I$ is primitive, that is, if $I$ is the annihilator of an irreducible module, then $I$ coincides with a single ideal (rather than the intersection of ideals) induced from isotropy. It would be highly desirable to have the inducing ideal in this result be primitive as well, but, as already pointed out by Steinberg in \cite {SteinbergEffrosHahn}, this is still an unresolved issue. Incidentally we should say that the techniques adopted in proving Theorem (12.14) are strongly influenced by \cite {SteinbergEffrosHahn}.

As already mentioned, all of this is set in the general context of regular inclusions, so we return to twisted Steinberg algebras in our final section with the purpose of showing that the abstract isotropy algebra is nothing but the twisted group algebra of the concrete isotropy group, as well as showing that the abstract imprimitivity bimodule coincides with the concrete bimodule used in the more popular induction process. Of course, this identification leads to immediate applications of our results to twisted Steinberg algebras.

\section Fell bundles over groupoids

In this section we shall be dealing with \'etale groupoids, so we start by introducing the basic definitions and notations. Given that such an introduction appears in most recent articles dealing with algebras associated to groupoids (see e.g.~\cite {Renault, Tight, SteinbergOne, LisaOne, ExelPitts}) we restrict ourselves to a very brief exposition.

A groupoid is, by definition, a small category in which every morphism is invertible. Therefore any groupoid $G$ comes equipped with \emph {range} and \emph {source} maps: $$ r: \gamma \in G\mapsto \gamma \gamma ^{-1}\in G,\and s: \gamma \in G\mapsto \gamma ^{-1}\gamma \in G. $$ Relevant sets are then defined as follows: $$ \def \crr {\hfill \cr \pilar {12pt}} \matrix { G\ex 0 & = & \big \{\gamma \in G: \gamma = s(\gamma )\big \},\crr G\ex 2 & = & \big \{(\alpha ,\beta )\in G\times G: s(\alpha ) = r(\beta )\big \}. } $$ In addition, if $x$ and $y$ are elements of $G\ex 0$, we define: $$ \def \crr {\hfill \cr \pilar {12pt}} \matrix { \G xx & = & \big \{\gamma \in G: s(\gamma )=x = r(\gamma )\big \}, \crr G_x & = & \big \{\gamma \in G: s(\gamma )=x\big \},\crr G^y & = & \big \{\gamma \in G: r(\gamma )=y\big \},\crr \G yx & = & \big \{\gamma \in G: s(\gamma )=x, \ r(\gamma )=y\big \}. } $$

We call $G\ex 0$ the \emph {unit space}, or the \emph {object space}, $G\ex 2$ is the the set of \emph {composable pairs}, and $\G xx$ is called the \emph {isotropy group} at $x$. The reader should keep in mind that the the placement of the variables $x$ and $y$ above will inspire the notation adopted for many sets introduced in this paper, especially after section (5).

A \emph {topological groupoid} is, by definition, a groupoid equipped with a topology making the operations of composition and inversion continuous. When $G\ex 0$ is locally compact and Hausdorff with its relative topology, and $r$ and $s$ are local homeomorphisms, one says that $G$ is \emph {\'etale}. If, moreover, the topology of $G\ex 0$ admits a basis of compact open subsets, we say that $G$ is \emph {ample}.

We should point out that an ample groupoid $G$ is not assumed to be Hausdorff. However, since $G$ admits local homeomorphisms into a Hausdorff space (e.g.~the range and source maps), it follows that $G$ is locally Hausdorff.

A \emph {bisection} in $G$ is, by definition, any subset $U\subseteq G$, such that the restriction of both the range and source maps to $U$ are injective. It is then a fact that the topology of every \'etale groupoid admits a basis of open bisections, while the topology of every ample groupoid admits a basis of compact open bisections.

A fact that we shall often use is that the product of two open bisections is again an open bisection, as proved in \cite [Proposition 2.2.4]{Paterson}.

\fix Throughout this section we fix an \'etale groupoid $G$ and a field $\fld $.

\state Definition \rm An \emph {algebraic Fell bundle} over $G$ is a set $E$ equipped with: \iStyle a \item a surjective map $\pi :E\to G$, \item a $\fld $-vector space structure on $E_\gamma := \pi ^{-1}(\gamma )$, for each $\gamma $ in $G$, \item a map $m:E\ex 2\to E$, where $$ E\ex 2 = \big \{(u,v)\in E\times E: (\pi (u),\pi (v))\in G\ex 2\big \}, $$ \medskip \noindent satisfying the following conditions: \zitemno = 0 \iStyle i \item for every $(\alpha , \beta )\in G\ex 2$, one has that $m(E_\alpha , E_\beta )\subseteq E_{\alpha \beta }$, \item $m$ is bilinear, \item $m$ is associative in the sense that, if $(u, v)$ and $(v, w)$ lie in $E\ex 2$, then $$ m(m(u,v),w) = m(u,m(v,w)). $$ \medskip \noindent We will say that $E$ is a \emph {topological Fell bundle} if $E$ is moreover equipped with a (not necessarily Hausdorff) topology such that: \smallskip \item $\pi $ is a local homeomorphism, \item the \emph {zero section} of $E$, namely the set $$ Z=\{0_\gamma :\gamma \in G\}, $$ where $0_\gamma $ denotes the zero vector in $E_\gamma $, is closed in $E$, \item addition is continuous as a map $$ +: E\times _\pi E:=\big \{(u,v)\in E\times E: \pi (u)=\pi (v)\big \} \to E, $$ \item for every $t$ in $\fld $, the map $$ u\in E \mapsto tu\in E $$ is continuous, \item $m$ is continuous.

\bigskip We shall mostly be concerned with topological Fell bundles, so, when we refer to Fell bundles without any adjective, we will have the topological version in mind.

\bigskip Given an algebraic Fell bundle $E$ over $G$, and for each $\gamma $ in $G$, we will call the subset $E_\gamma $ mentioned in (2.1.b), the \emph {fiber} of $E$ over $\gamma $. The map $m$ will be called the \emph {multiplication operator} for $E$, and from now on we will adopt the shorthand notation $$ uv := m(u,v), \for (u,v)\in E\ex 2. $$ With this, the associativity axiom (2.1.iii) takes on the more familiar form: $$ (uv)w=u(vw). $$

Recall that every \'etale groupoid is locally Hausdorff (that is, every point admits an open neighborhood that is Hausdorff in the relative topology). Thus, in case $E$ is a topological bundle, the fact that $\pi $ is a local homeomorphism implies that $E$ is also locally Hausdorff.

\state Definition \rm If $U$ is any open subset of $G$, and $\xi :U\to E$ is a function, then $\xi $ is called a \emph {local section} if $\pi \circ \xi $ is the identity map on $U$. The domain and range of $\xi $ will be respectively denoted by $\dom (\xi )$ and $\ran (\xi )$.

Thus a local section may be seen as a selection of one element in $E_\gamma $ for each $\gamma $ in $\dom (\xi )$. When $E$ is a topological bundle, the \emph {continuous} local sections will play a predominant role. In any case, all local sections involved in this work will strictly follow the requirement of having an open domain.

Local sections enjoy many nice properties, such as the following: if $\xi $ and $\eta $ are local sections such that $\xi (\gamma )\in \ran (\eta )$, for some $\gamma $ in $\dom (\xi )$, then $\gamma $ lies in $\dom (\eta )$, and $\xi (\gamma )=\eta (\gamma )$.

The reader used to the concept of Steinberg algebras is certainly aware of the relevance of locally constant functions on groupoids. However, a local section $\xi $ is never locally constant, since each $\xi (\gamma )$ lies in its own fiber $E_\gamma $, and these are pairwise disjoint. Nevertheless, the fact that a continuos local section is a right inverse of $\pi $, and the fact that $\pi $ is a local homeomorphism, should be interpreted as saying that $\xi $ is somehow locally constant. A concrete manifestation of this phenomenon is (2.3.iv), below.

\state Proposition Given a topological Fell bundle $E$ over the \'etale groupoid $G$, the following hold: \item If $\xi $ and $\eta $ are continuous local sections with the same domain, then the pointwise sum $\xi +\eta $ is continuous, and so is the pointwise scalar multiple $t\, \xi $, for every $t$ in $\fld $. \item For every $u$ in $E$, there exists a continuous local section $\xi $ \emph {passing through} $u$, meaning that $u$ lies in the range of $\xi $, and hence that $\xi (\pi (u))=u$. \item The range of any continuous local section is open. \item Given continuous local sections $\xi $ and $\eta $ such that $\xi (\gamma _0)=\eta (\gamma _0)$, for some $\gamma _0$ lying in both of their domains, there exists an open neighborhood $V$ of $\gamma _0$, contained in $ \dom (\xi )\cap \dom (\eta ), $ such that $\xi |_V=\eta |_V$. \item The support of every continuous local section $\xi $, namely the set $$ \supp (\xi ) := \{\gamma \in \dom (\xi ): \xi (\gamma )\neq 0_\gamma \}, $$ is both open and closed relative to $\dom (\xi )$.

\Proof Point (i) follows readily from (2.1.vi) and (2.1.vii). The proof of (ii) is an easy consequence of the fact that $\pi $ is a local homeomorphism. Regarding (iii), let $\xi $ be a continuous local section, let $u_0\in \ran (\xi )$, and put $\gamma _0=\pi (u_0)$, so that necessarily $u_0=\xi (\gamma _0)$. Since $\pi $ is a local homeomorphism, we may choose open sets $U\subseteq E$ and $V\subseteq G$, with $u_0\in U$, $\gamma _0\in V$, and such that $\pi $ is a homeomorphism from $U$ to $V$.

By continuity there exists an open neighborhood $W$ of $\gamma _0$, contained in $\dom (\xi )$, such that $\xi (W)\subseteq U$, and we may clearly assume that $W\subseteq V$.

Since $\pi $ maps $\xi (W)$ onto $W$, which is open, and since $\pi $ is a homeomorphism from $U$ to $V$, it follows that $\xi (W)$ is open. Finally, noting that $$ u_0=\xi (\gamma _0)\in \xi (W)\subseteq \ran (\xi ), $$ we see that $u_0$ belongs to the interior of $\ran (\xi )$, as desired.

In order to prove (iv) we argue as in (iii) and choose open sets $U\subseteq E$ and $V\subseteq G$, with $\xi (\gamma _0)\in U$, $\gamma _0\in V$, and such that $\pi $ is a homeomorphism from $U$ to $V$. As before, we may find an open neighborhood $W$ of $\gamma _0$ contained in $\dom (\xi )\cap \dom (\eta )$, such that $\xi (W)$ and $\eta (W)$ are both contained in $U$. Observing that, for every $\gamma \in W$, $$ \pi (\xi (\gamma )) = \gamma = \pi (\eta (\gamma )), $$ the fact that $\pi $ is injective on $U$ implies that $\xi (\gamma ) = \eta (\gamma )$.

In order to address the last remaining point, recall that the zero section $$ Z:=\{0_\gamma :\gamma \in G\}, $$ defined in (2.1.v), is closed in $E$ by definition. Nevertheless we claim that it is also open, the reason being as follows: given any $\gamma $ in $G$, and using (ii), choose a continuous local section $\eta $ such that $\eta (\gamma )=0_\gamma $. Then the section $\zeta $ defined by $$ \zeta (\alpha )=0{\cdot }\eta (\alpha )=0_\alpha , \for \alpha \in \dom (\eta ), $$ is continuous by (i), and hence its range is an open set by (iii). Observing that $$ 0_\gamma \in \ran (\zeta )\subseteq Z, $$ we deduce that $0_\gamma $ is an interior point of $Z$, proving the claim. Now, since the support of any given continuous local section $\xi $ coincides with $\xi ^{-1}(E\setminus Z)$, the proof follows easily. \endProof

\state Definition \rm A \emph {Fell line bundle} over $G$ (\emph {line bundle for short}) is a topological Fell bundle $E$ such that \item each $E_\gamma $ is a one-dimensional $\fld $-vector space, and \item for every $(\alpha , \beta )$ in $G\ex 2$, the multiplication operator restricted to $E_\alpha \times E_\beta $ is not identically zero.

Given a line bundle $E$, observe that, since each $E_\gamma $ is one dimensional, and since for each $(\alpha ,\beta )$ in $G\ex 2$ the multiplication is a nontrivial bilinear operation on $E_\alpha \times E_\beta $, it is easy to see that the multiplication is also nondegenerate in the sense that if $uv=0_{\alpha \beta }$, for some $(u, v)\in E_\alpha \times E_\beta $, then either $u=0_\alpha $, or $v=0_\beta $.

Important examples of line bundles are obtained from 2-cocycles, as defined below, but before that we need to emphasize that any time the field $\fld $ is viewed as a topological space, it is supposed to have the \emph {discrete topology}. It is then useful to keep in mind that a $\fld $-valued function is continuous if and only if it is locally constant.

\state Definition \rm A \emph {2-cocycle} on $G$ is a continuous map $$ \omega :G\ex 2\to \fld ^*, $$ where $\fld ^* = \fld \setminus \{0\}$, such that, \item for every $\gamma $ in $G$, one has that $ \omega (\gamma ,s(\gamma )) = \omega (r(\gamma ), \gamma ) = 1, $ \item whenever $(\alpha ,\beta )$ and $(\beta ,\gamma )$ are in $G\ex 2$, one has that $$ \omega (\alpha ,\beta )\omega (\alpha \beta ,\gamma ) = \omega (\alpha ,\beta \gamma )\omega (\beta ,\gamma ). $$

\bigskip The requirement that a $\fld $-valued function be continuous, hence locally constant, forces it to be constant in case its domain is connected. So the more disconnected $G$ is, the more room there is for nontrivial 2-cocycles. In fact we will mostly be interested in \emph {ample} groupoids, whose topology admits a basis of compact open sets, and hence are highly disconnected.

Given a $2$-cocycle on the \'etale groupoid $G$, one may build a line bundle as follows: let $E$ be the product topological space $$ E=\fld \times G, $$ and let $\pi $ be the projection onto the second coordinate. Thus $E_\gamma $ becomes $\fld \times \{\gamma \}$, which we view as a $\fld $-vector space in the obvious way.

As for the multiplication operator, we take any $(u,v)$ in $E\ex 2$, write $u=(t,\alpha )$, and $v=(s,\beta )$, so that necessarily $(\alpha ,\beta )\in G\ex 2$, and we put $$ m(u,v) = (t,\alpha )(s,\beta ) = \big (\omega (\alpha ,\beta )ts,\alpha \beta \big ). $$ One may then easily prove that $E$ is a line bundle over $G$.

\state Definition \rm Given a 2-cocycle $\omega $ on $G$, we will denote the line bundle constructed above by $E(\omega )$. In the special case that $\omega $ is identically equal to $1$, we will call $E(\omega )$ the \emph {trivial line bundle}.

The notion of \emph {convolution product} is an important tool in the theory of line bundles over groupoids to be extensively discussed in the later sections of this work. Nevertheless there is an elementary aspect of it that is worth introducing without delay. In order to describe it precisely, suppose that we are given two open bisections $U_1$ and $U_2$ of $G$. Setting $$ U_1 U_2:=\{\gamma _1\gamma _2: (\gamma _1, \gamma _2)\in (U_1\times U_2)\cap G\ex 2\}, $$ observe that any $\gamma $ in $U_1 U_2$ may be uniquely factorized as $\gamma =\gamma _1\gamma _2$, with $\gamma _i\in U_i$. This is because there is at most one $\gamma _1$ in $U_1$ whose range coincides with $r(\gamma )$, and likewise there is at most one $\gamma _2$ in $U_2$ whose source coincides with $s(\gamma )$. Incidentally $\gamma _1$ and $\gamma _2$ may be described in terms of $\gamma $ by $$ \gamma _1=(r|_{U_1})^{-1}(r(\gamma )), \and \gamma _2=(s|_{U_2})^{-1}(s(\gamma )), \reqno{(2.7)} $$ where we are making use of the fact that the restrictions of both $r$ and $s$ to any bisection is injective. We will refer to the pair $(\gamma _1,\gamma _2)$ as the \emph {unique factorization} of $\gamma $ relative to the bisections $U_1$ and $U_2$, the components of which clearly depend continuously on $\gamma $ by (2.7).

\state Definition \rm Given an algebraic Fell bundle $E$ over $G$, and given any two local sections $\xi $ and $\eta $, whose domains are open bisections, we define the \emph {mini convolution product} of $\xi $ and $\eta $ to be the local section defined on the product of the bisections $\dom (\xi )$ and $\dom (\eta )$ by $$ \xi *\eta : \gamma \in \dom (\xi ){\cdot }\dom (\eta ) \mapsto \xi (\gamma _1)\eta (\gamma _2)\in E, $$ where $(\gamma _1,\gamma _2)$ is the unique factorization of $\gamma $ relative to $\dom (\xi )$ and $\dom (\eta )$.

Due to (2.7), it is easy to see that, if $E$ is topological, and $\xi $ and $\eta $ are continuous, then so is $\xi *\eta $.

Another important aspect of line bundles that we must discuss is the presence of certain special elements in the fibers over $G\ex 0$. With this purpose in mind, note that, for every unit $x\in G\ex 0$, the axioms imply that $E_x$ is a $\fld $-algebra. Recalling that any one-dimensional $\fld $-algebra with a non-identically zero multiplication is necessarily isomorphic to the field $\fld $, we see that $E_x$ possesses a distinguished element, namely its unit, which we will denote by $1_x$.

Naturally, we have that $1_xu=u$, for every $u$ in $E_x$, but we claim that this is also true for every $u$ in any $E_\gamma $, provided $r(\gamma )=x$. To see this, notice that the multiplication operator $$ E_{r(\gamma )}\times E_\gamma \to E_\gamma $$ is nonzero by hypothesis, and hence it is necessarily a surjective map. It follows that every $u$ in $E_\gamma $ may be written as $u=vw$, for some $(v, w)$ in $E_{r(\gamma )}\times E_\gamma $, whence $$ 1_xu = 1_x(vw) = (1_xv)w = vw = u. \reqno{(2.9)} $$ By a similar argument we also have that $$ u1_x=u, \reqno{(2.10)} $$ for any $u$ in any $E_\gamma $ such that $s(\gamma )=x$.

Observe that, if $\omega $ is a 2-cocycle on $G$, then the bundle $E(\omega )$ admits a nowhere vanishing continuous global section\fn {By a \emph {global section} we mean a local section whose domain is all of $G$.}, such as the one defined by $$ \xi (\gamma )=(1,\gamma ), \for \gamma \in G. $$ In what follows we will show that this property actually characterizes line bundles arising from $2$-cocycles.

\state Proposition Let $E$ be a line bundle over $G$. \item If $E$ admits a nowhere vanishing continuous global section, then $E$ is isomorphic\fn {In the sense that there exists a bijective function preserving all their structure.} to $E(\omega )$ for some continuous $2$-cocycle on $G$. \item The map $$ \varepsilon :x\in G\ex 0 \mapsto 1_x\in E, $$ is a continuous local section, whence the restriction of $E$ to $G\ex 0$, namely $F:=\pi ^{-1}(G\ex 0)$, admits a nowhere vanishing continuous \emph {global} section, so $F$ is isomorphic to the trivial line bundle over $G\ex 0$.

\Proof Addressing (i), let $\varepsilon $ be a nowhere vanishing continuous global section of $E$ and observe that, for each $(\alpha , \beta )\in G\ex 2$, one has that $\varepsilon (\alpha )\varepsilon (\beta )$ lies in $E_{\alpha \beta }$. Since $\varepsilon (\alpha \beta )$ is nonzero, it is a linear basis for this vector space, so there exists a unique nonzero scalar $\omega (\alpha , \beta )\in \fld $, such that $$ \varepsilon (\alpha )\varepsilon (\beta )=\omega (\alpha , \beta )\varepsilon (\alpha \beta ). \reqno{(2.11.1)} $$

We may then see $\omega $ as a $\fld $-valued function defined on $G\ex 2$, and our next goal is to show that $\omega $ is locally constant, and hence continuous. To see this, fix $(\alpha _0,\beta _0)$ in $G\ex 2$, and choose open bisections $U$ and $V$ containing $\alpha _0$ and $\beta _0$, respectively. Set $\xi =\varepsilon |_U$ and $\eta =\varepsilon |_V$, so that $\xi $ and $\eta $ are continuous local sections defined on bisections and hence the mini convolution product $\xi *\eta $ is a well defined continuous local section defined on $U V$.

Given $\gamma $ in $U V$, let $(\alpha ,\beta )$ be its unique factorization relative to $U$ and $V$, and notice that $$ \omega (\alpha , \beta )\varepsilon (\gamma ) = \varepsilon (\alpha )\varepsilon (\beta )= \xi (\alpha )\eta (\beta ) = (\xi *\eta )(\gamma ). \reqno{(2.11.2)} $$ Defining $\zeta $ to be the pointwise product $\omega (\alpha _0, \beta _0)\varepsilon $ on $U V$, we have that $\zeta $ is a continuous local section by (2.3.i), and clearly $$ \zeta (\gamma _0) = (\xi *\eta )(\gamma _0), $$ where $\gamma _0:= \alpha _0\beta _0$, by (2.11.2). Thanks to (2.3.iv) it then follows that there exists some neighborhood $W$ of $\gamma _0$, contained in $UV$, such that $$ \omega (\alpha _0, \beta _0)\varepsilon (\gamma ) = \zeta (\gamma ) = (\xi *\eta )(\gamma ), \for \gamma \in W. \reqno{(2.11.3)} $$ We next invoke the continuity of the groupoid multiplication operation and pick open sets $U'$ and $V'$ such that $$ \alpha _0\in U'\subseteq U, \and \beta _0\in V'\subseteq V, $$ and such that said operation maps $(U'\times V')\cap G\ex 2$ into $W$. For every $(\alpha ,\beta )$ lying in the former set, we then have that $\gamma :=\alpha \beta $ lies in $W$, whence $$ \omega (\alpha _0, \beta _0)\varepsilon (\gamma ) \explica {(2.11.3)}{=} (\xi *\eta )(\gamma ) \explica {(2.11.2)}{=} \omega (\alpha , \beta )\varepsilon (\gamma ). $$ This in turn implies that $\omega (\alpha _0, \beta _0) = \omega (\alpha , \beta )$, as desired. This concludes the proof of the continuity of $\omega $, and we would now like to prove it to be a $2$-cocycle.

An easy computation based on the associativity property of line bundles gives (2.5.ii), but unfortunately (2.5.i) seems to be out of our reach. Nevertheless we will show that it is possible to replace our global section $\varepsilon $ with a nicer one, so that (2.5.i) is also satisfied.\fn {An easy path to a proof of this claim would be to replace the values of $\varepsilon (x)$ by $1_x$, for every $x$ in $G\ex 0$, in which case the map $\omega $ defined by (2.11.1) clearly satisfies (2.5.i), as a quick computation shows. However, even though $\varepsilon $ would then be continuous on $G\ex 0$, thanks to (ii), one might not be able to show that $\varepsilon $ is continuous on all of $G$ in case $G$ is not Hausdorff. The difficulty arises from the fact that $G\ex 0$ is not closed in a non-Hausdorff groupoid.} For this, fixing any $\gamma $ in $G$, denoting by $x=r(\gamma )$, and plugging $\alpha =x$, and $\beta =\gamma $, in (2.11.1), we get $$ \varepsilon (x)\varepsilon (\gamma )= \omega (x, \gamma ) \varepsilon (\gamma ) \explica {(2.9)}{=} \omega (x, \gamma ) 1_{x}\varepsilon (\gamma ). $$ Since $\varepsilon (\gamma )$ is nonzero, we deduce that $$ \varepsilon (x) = \omega (x, \gamma ) 1_{x}. $$ Repeating the above argument with $\gamma $ replaced by $x$, we obtain $$ \varepsilon (x) = \omega (x, x) 1_{x}, $$ which implies that $\omega (x, \gamma ) = \omega (x, x)$, or equivalently $$ \omega (r(\gamma ), \gamma ) = \omega (r(\gamma ), r(\gamma )), \for \gamma \in G. $$ Likewise, plugging $\alpha =\gamma $, and $\beta =s(\gamma )$ in (2.11.1), a similar argument leads to $$ \omega (\gamma , s(\gamma )) = \omega (s(\gamma ), s(\gamma )), \for \gamma \in G. $$ In order to simplify our notation from now on, we adopt the notation $$ \rho (x) := \omega (x, x), \for x \in G\ex 0, $$ and we observe that, for all $\gamma $ in $G$, we have that $$ \omega (r(\gamma ), \gamma ) = \rho (r(\gamma )), \and \omega (\gamma , s(\gamma )) = \rho (s(\gamma )). \reqno{(2.11.4)} $$ We next consider a new global section $\varepsilon '$, defined by $$ \varepsilon '(\gamma ) = \rho \big (s(\gamma )\big )^{-1}\varepsilon (\gamma ), \for \gamma \in G. $$ Since $\omega $ is locally constant, so is $\rho $, whence $\varepsilon '$ is continuous.

Defining a new $\omega '$ based on $\varepsilon '$, in the same way that $\omega $ was defined based on $\varepsilon $ in (2.11.1), we have that $$ \varepsilon '(\alpha )\varepsilon '(\beta )=\omega '(\alpha , \beta )\varepsilon '(\alpha \beta ). $$ In order to see how does $\omega $ and $\omega '$ relate to each other, we observe that $$ \varepsilon '(\alpha )\varepsilon '(\beta ) = \rho \big (s(\alpha )\big )^{-1}\varepsilon (\alpha ) \rho \big (s(\beta )\big )^{-1}\varepsilon (\beta )= \rho \big (s(\alpha )\big )^{-1} \rho \big (s(\beta )\big )^{-1}\omega (\alpha , \beta )\varepsilon (\alpha \beta ), $$ while $$ \omega '(\alpha , \beta )\varepsilon '(\alpha \beta ) = \omega '(\alpha , \beta )\rho \big (s(\alpha \beta )\big )^{-1}\varepsilon (\alpha \beta ) = \omega '(\alpha , \beta )\rho \big (s(\beta )\big )^{-1}\varepsilon (\alpha \beta ). $$ Comparing the above we deduce that $\omega $ and $\omega '$ relate to each other by means of the formula $$ \omega '(\alpha , \beta ) = \rho \big (s(\alpha )\big )^{-1} \omega (\alpha , \beta ). $$

The same reasons why $\omega $ satisfies (2.5.ii) apply here to show that so does $\omega '$, and we claim that $\omega '$ also satisfies (2.5.i). To see this, pick $\gamma $ in $G$, and observe that $$ \omega '\big (\gamma ,s(\gamma )\big ) = \rho \big (s(\gamma )\big )^{-1} \omega (\gamma , s(\gamma )) \explica {(2.11.4)}{=} 1, $$ while $$ \omega '\big (r(\gamma ), \gamma \big ) = \rho \big (s(r(\gamma ))\big )^{-1} \omega \big (r(\gamma ), \gamma \big ) = \rho \big (r(\gamma )\big )^{-1} \omega \big (r(\gamma ), \gamma \big ) \explica {(2.11.4)}{=} 1. $$ This shows that $\omega '$ is a 2-cocycle on $G$.

Since we have no more use for the original global section $\varepsilon $, nor the failed cocycle candidate $\omega $, we will henceforth denote $\varepsilon '$ by $\varepsilon $, and $\omega '$ by $\omega $, so that, as before, $\varepsilon $ is a nowhere vanishing continuous global section, (2.11.1) still holds, but now $\omega $ is a legitimate $2$-cocycle for $G$.

We will next prove that $E$ is isomorphic to $E(\omega )$, and for this we consider the map $$ \varphi :(t, \gamma )\in \fld \times G \mapsto t\varepsilon (\gamma )\in E. $$ Since $\varepsilon $ is continuous, so is $\varphi $, and we claim that $\varphi $ is also an open map. To see this, let $U$ be an open subset of $\fld \times G$, of the form $U=\{t\}\times V$, where $V$ is open in $G$. Then clearly $\varphi (V)$ coincides with the range of the continuous local section $\xi (\gamma )=t\varepsilon (\gamma )$, defined on $V$, so $\varphi (V)$ is open by (2.3.iii). Since the collection of open sets $U$ considered is a basis for the topology of $\fld \times G$, the claim is proved. Since $\varphi $ is clearly bijective, we see that it is a homeomorphism.

Finally, leaving for the reader to verify that $\varphi $ preserves all of the algebraic operations involved, the proof of (i) is finished.

Addressing (ii), the only outstanding point is now to prove $\varepsilon $ to be continuous, so fix any $x_0$ in $G\ex 0$, and let $\xi $ be a continuous local section of $E$, defined on an open neighborhood $U$ of $x_0$, with $\xi (x_0)=1_{x_0}$. It follows that the map $$ \eta :x\in U\mapsto \xi (x)^2-\xi (x)\in E $$ is another continuous local section, but this time $\eta (x_0)=0$. Letting $\zeta $ be the map defined by $$ \zeta (x)=0{\cdot }\xi (x) = 0_x,\for x\in U, $$ we have that $\zeta $ is continuous as a consequence of (2.1.vii), and clearly $\eta (x_0)=\zeta (x_0)$. By (2.3.iv) there exists a neighborhood $V$ of $x_0$, contained in $U$, such that $$ \eta (x)=\zeta (x)=0_x, \for x\in V. $$ This implies that $\xi (x)$ is idempotent for every $x$ in $V$, whence either $\xi (x)=1_x$, or $\xi (x)=0_x$.

We then claim that there exists a further neighborhood $W$ of $x_0$, contained in $V$, such that $\xi (x)=1_x$, for every $x$ in $W$. Indeed, if this was not true, we could find a net $(x_i)_i$ converging to $x_0$, such that $\xi (x_i)=0_{x_i}$, for all $i$. It would then follow that $$ 0_{x_i} = \xi (x_i)\to \xi (x_0) = 1_{x_0}, $$ contradicting axiom (2.1.v). This shows that $\xi $ coincides with $\varepsilon $ on a neighborhood of $x_0$, so it follows that $\varepsilon $ is continuous at $x_0$, and hence also everywhere. This proves the first part of (ii) which, combined with (i), concludes the proof. \endProof

In \cite {BeckyOne} the authors study Steinberg algebras over a \emph {twisted} ample groupoid, defined to be an ample groupoid equipped with a 2-cocycle. There, the authors also hint at the possibility of generalizing this concept to groupoid extensions (see also \cite {Renault, RenaultCartan}), but in \cite [Theorem 4.10]{BeckyOne} they show that this generalization is pointless over Hausdorff groupoids since all instances of the generalized concept can be realized by 2-cocycles.

Their proof of this Theorem is based on the ability to patch local sections into a global one, but it turns out that this patching procedure breaks down over non Hausdorff groupoids. In fact, we will soon see an example to illustrate this phenomena. It is then our point of view that, after all, the study of more general \emph {topological twists} is worth the effort.

Another obvious consequence of the non-degeneracy of the product in a line bundle $E$, is as follows:

\state Proposition Let $E$ be a line bundle over $G$. Then, for every nonzero $u$ in any $E_\gamma $ there exists a unique $u^{-1}$ in $E_{\gamma ^{-1}}$ such that: \item $u^{-1}u=1_{s(\gamma )}$, \item $uu^{-1}=1_{r(\gamma )}$, \item if $Z$ is the zero section of $E$, then the map $$ u\in E\setminus Z \to u^{-1}\in E\setminus Z $$ is continuous, involutive, and anti-multiplicative.

\Proof Since the multiplication operation restricted to $$ E_{\gamma ^{-1}}\times E_\gamma \to E_{s(\gamma )} $$ is nontrivial, the existence and uniqueness of $u^{-1}$ satisfying (i) is granted. As for (ii), notice that $$ uu^{-1}u=u1_{s(\gamma )} \explica {(2.10)}{=} u \explica {(2.9)}{=} 1_{r(\gamma )}u, $$ so $uu^{-1}=1_{r(\gamma )}$, by nondegeneracy.

In order to prove continuity, pick a nonzero $u$ in $E$, and let $U$ be an open neighborhood of $u$, and $V$ be an open neighborhood of $\pi (u)$, such that $\pi $ defines a homeomorphism from $U$ to $V$. Likewise, let $U'$ be an open neighborhood of $u^{-1}$, and $V'$ be an open neighborhood of $\pi (u^{-1})$, such that $\pi $ defines a homeomorphism from $U'$ to $V'$. Observing that $\pi (u^{-1})=\pi (u)^{-1}$, and using that the inversion map is a homeomorphism on $G$, we may assume that $V'=V^{-1}$.

\begingroup \noindent \hfill \beginpicture \setcoordinatesystem units <0.025truecm, -0.02truecm> \setplotarea x from -50 to 150, y from -30 to 150 \put {\null } at -50 -30 \put {\null } at -50 150 \put {\null } at 150 -30 \put {\null } at 150 150 \put {$U$} at 0 0 \put {$U'$} at 100 0 \arrow <0.11cm> [0.5, 1.8] from 20 0 to 80 0 \put {$\varphi $} at 50 -15 \put {$V$} at 0 100 \put {$V'$} at 100 100 \arrow <0.11cm> [0.5, 1.8] from 20 100 to 80 100 \put {$\text {inv}$} at 50 115 \arrow <0.11cm> [0.5, 1.8] from 0 18 to 0 82 \put {$\pi $} at -10 50 \arrow <0.11cm> [0.5, 1.8] from 100 18 to 100 82 \put {$\pi $} at 90 50 \endpicture \hfill \null \endgroup

One may then fill in the top horizontal map with a unique homeomorphism $\varphi $, making the diagram commute. Observing that $\varphi (u)=u^{-1}$, we conclude that $\varphi (u)u=1_{s(u)}$. Thus the continuous map $$ w\mapsto \varphi (w)w $$ sends $u$ into the range of the continuous section $\varepsilon $ referred to in (2.11.ii). Since this range is open by (2.3.iii), there exists a neighborhood $W$ of $u$, contained in $U$, such that $\varphi (w)w$ lies in said range, for all $w$ in $W$, so that necessarily $$ \varphi (w)w = 1_{s(w)}, \for w\in W. $$ It follows that $w^{-1}=\varphi (w)$, thus proving the continuity of the map in the statement. The remaining assertions in (iii) are left as easy exercises. \endProof

\section Constructing Fell bundles

Let us now describe an efficient way to construct topological Fell bundles over a given \'etale groupoid $G$. Actually, there is not much we have to say about the construction of the algebraic ingredients, namely the projection $\pi $, the vector space structure on each $E_\gamma $, and the multiplication operator $m$, which most of the time may be explicitly constructed without much difficulty. The trickier part is usually the description of the topology, so the purpose of this section is to facilitate the construction of the topology and the verification of the pertinent axioms.

\fix We thus assume that $G$ is an \'etale groupoid, and we are already given an algebraic Fell bundle $E$ over $G$.

\medskip Mimicking (2.3), we now introduce the following concept:

\state Definition \rm By a \emph {fundamental family of local sections} of $E$ we mean a set $\cal S$ of local sections (with open domains) such that: \item every $u$ in $E$\/ lies in the range of some $\xi $ in $\S $, \item for every $\xi $ in $\S $, the restriction of $\xi $ to any open subset of its domain also lies in $\S $, \item given $\xi $ and $\eta $ in $\S $, such that $\xi (\gamma _0)=\eta (\gamma _0)$, for some $\gamma _0$ lying in both of their domains, there exists an open neighborhood $V$ of $\gamma _0$, contained in $ \dom (\xi )\cap \dom (\eta ), $ such that $\xi |_V=\eta |_V$. \item the support of every $\xi $ in $\S $ is open, \item if $\xi $ and $\eta $ are local sections in $\S $ with the same domain, than the pointwise sum $\xi +\eta $ lies in $\S $, and so does the pointwise scalar multiple $t\, \xi $, for every $t$ in $\fld $, \item if $\xi _1$ and $\xi _2$ are local sections in $\S $ whose domains are open bisections of $G$, then the mini convolution product $\xi _1*\xi _2$ lies in \S .

\bigskip Fixing a fundamental family of local sections $\S $, as above, we then consider the topology on $E$ generated by the family of subsets given by $$ \mathcal B = \{\ran (\xi ): \xi \in \S \}, $$ where $\ran (\xi )$ refers to the range of $\xi $.

We observe that $\mathcal B$ is indeed a basis for this topology for the following reason: given $\xi $ and $\eta $ in $\S $, and given $u\in \ran (\xi )\cap \ran (\eta )$, let $\gamma =\pi (u)$, so that $\xi (\gamma )=\eta (\gamma )=u$. Letting $V$ be as in (3.1.iii), and letting $\zeta $ be the restriction of either $\xi $ of $\eta $ to $V$, we have that $$ u\in \ran (\zeta )\subseteq \ran (\xi )\cap \ran (\eta ), $$ thus verifying the well known criterion for a topology basis.

\state Proposition The topology defined above makes $E$ a topological Fell bundle over $G$, such that every $\xi $ in $\S $ is continuous. In addition, for every continuous local section $\eta $ of $E$, and every $\gamma $ in $\dom (\eta )$, there exists some $\xi $ in $\S $ which coincides with $\eta $ on some neighborhood of $\gamma $.

\Proof In order to prove that $E$ is a topological Fell bundle, it clearly suffices to verify the \emph {topological} axioms of (2.1), namely the ones beginning with, and including (2.1.iv).

In order to do so, we first check that $\pi $ is continuous: fix $u$ in $E$, and let $V$ be an open neighborhood of $\gamma := \pi (u)$. Using (3.1.i), choose $\xi $ in $\S $ such that $\xi (\gamma )=u$ and, thanks to (3.1.ii), we may assume that $\dom (\xi )\subseteq V$. It follows that $$ u\in \ran (\xi )\subseteq \pi ^{-1}(V). $$ Observing that $\ran (\xi )$ is an open neighborhood of $u$, this shows that $\pi $ is indeed continuous.

Next observe that $\pi $ maps each basic open set of the form $\ran (\xi )$ onto $\dom (\xi )$, and hence $\pi $ is an open map. Since $\pi $ is one-to-one on $\ran (\xi )$, we then deduce that $\pi $ defines a homeomorphism from $\ran (\xi )$ to $\dom (\xi )$, thus proving that $\pi $ is a local homeomorphism. The inverse of this homeomorphism is clearly $\xi $, whence $\xi $ is continuous.

To check (2.1.v), pick any $u$ in $E$, let $\gamma =\pi (u)$, and assume that $u\neq 0_\gamma $. By (3.1.i), we may fix some $\xi $ in $\S $ such that $\xi (\gamma )=u$, so the set $$ V=\{\alpha \in \dom (\xi ): \xi (\alpha )\neq 0_\alpha \}, $$ contains $\gamma $, and is open by (3.1.iv). By (3.1.ii), we may further assume that $\dom (\xi )=V$, which is the same as saying that $\xi $ vanishes nowhere. The open set $\ran (\xi )$ is therefore a neighborhood of $u$, which does not intercept the set $\{0_\alpha :\alpha \in G\}$ of all zeros, whence the latter is seen to be a closed set.

To prove the continuity of vector addition, choose $u_0$ and $v_0$ in some $E_{\gamma _0}$, and let $\xi $ and $\eta $ be local sections in $\S $, defined on the same open neighborhood $V$ of $\gamma _0$, such that $\xi (\gamma _0)=u_0$, and $\eta (\gamma _0)=v_0$. Considering the basic open neighborhoods $\ran (\xi )$ of $u_0$, and $\ran (\eta )$ of $v_0$, we have that $$ W:= \big (\ran (\xi )\times \ran (\eta )\big )\cap (E\times _\pi E) = \big \{\big (\xi (\gamma ),\eta (\gamma )\big ): \gamma \in V\big \}. $$ Consequently the addition operator ``$+$'' on the above set coincides with the map $$ (u,v) \mapsto (\xi +\eta )(\pi (u)), $$ which is continuous because both $\xi +\eta $ and $\pi $ are continuous. We leave the proof of the continuity of scalar multiplication to the reader.

It now remains to prove (2.1.viii), namely continuity of the multiplication on $E$. We thus fix a pair $(u_0,v_0)\in E\ex 2$, and choose $\xi $ and $\eta $ in $\S $, such that $u_0\in \ran (\xi )$, and $v_0\in \ran (\eta )$. Suitably restricting $\xi $ and $\eta $ we may assume without loss of generality that both their domains are open bisections. Given any $$ (u,v)\in (\ran (\xi )\times \ran (\eta ))\cap E\ex 2, $$ and setting $\alpha =\pi (u)$, and $\beta =\pi (v)$, we have that $(\alpha ,\beta )\in (\dom (\xi )\times \dom (\eta ))\cap G\ex 2$, so the unique factorization of $\alpha \beta $ relative to $\dom (\xi )$ and $\dom (\eta )$ is precisely $(\alpha , \beta )$, and we deduce that $$ uv = \xi (\alpha )\eta (\beta ) = (\xi *\eta )(\alpha \beta ) = (\xi *\eta )(\pi (u)\pi (v)), $$ from where we see that the multiplication operator is continuous at $(u_0,v_0)$. This concludes the proof that $E$ is a topological Fell bundle over $G$.

Addressing the last sentence in the statement, fix any continuous local section $\eta $, and pick $\gamma $ in $\dom (\eta )$. By (2.3.iii) we have that $\ran (\eta )$ is open and hence there exists a basic neighborhood $U$ of $\eta (\gamma )$ contained in $\ran (\eta )$, where by \emph {basic} we of course mean that $U=\ran (\xi )$, for some $\xi $ in $\S $. From $\eta (\gamma )\in \ran (\xi )$ it easily follows that $\eta (\gamma )=\xi (\gamma )$, so the claim is a consequence of (2.3.iv) and (3.1.ii). \endProof

We would now like to present an example to show that not all line bundles are of the form $E(\omega )$, for a cocycle $\omega $, as defined in (2.6). For this let $X$ be the compact, totally disconnected topological space $$ X=\big \{1/n: n\in {\mathbb Z},\ n\neq 0\big \} \cup \{0\} = $$ $$ \def \f #1{\frac 1{\kern 2pt#1\kern 2pt}} = \Big \{-1, -\f 2, -\f 3,-\f 4, \ldots , 0, \ldots ,\f 4, \f 3, \f 2,1 \Big \} $$ seen as a topological subspace of ${\mathbb R}$. Thus $X$ consist of the ranges of two sequences converging to zero, one increasing and another decreasing, together with their common limit zero.

We then consider the groupoid $X\times {\mathbb Z}_2$, seen as a group bundle over $X$, equipped with the product topology and we put $$ G=X\times {\mathbb Z}_2/{\sim }, $$ where ``$\sim $'' is the equivalence relation identifying $(x,0)$ with $(x,1)$, for all $x$ in $X$, except for $x=0$. The unit space of $G$ is therefore given by $$ G\ex 0 = \big \{\overline {(x,0)} : x\in X\big \}, $$ (overline referring to equivalence class), while the only element in $G$ which is not a unit is $\overline {(0,1)}$. In order to simplify our notation, we set $$ \extra = \overline {(0,1)}, $$ and we will tacitly identify $G\ex 0$ with $X$, in the obvious way. So, $$ G=X \cup \{\extra \}. $$

Therefore $G$ is again a group bundle, except that now the isotropy groups are all trivial with the exception of $G(0)$, which turns out to be the cyclic group of order $2$, with generator $\extra $ and neutral element $0$.

Regarding the topology on $G$, namely the quotient topology, observe that a subset $U\subseteq G$ is open if and only if either $U\subseteq G\setminus G(0)$, or there is some $\varepsilon >0$, such that $ I_\varepsilon \subseteq U, $ where $$ I_\varepsilon = X \cap (-\varepsilon ,\varepsilon )\setminus \{0\}. $$ \def \dash {\smallskip \noindent --\ }\relax It is then easy to see that: \dash $X$ is an open subset of $G$, \dash the relative topology of $X$, as a subspace of $G$, coincides with its own topology, \dash except for $0$ and $\extra $, all points of $G$ are isolated points, \dash a neighborhood base for $0$ is given by the sets of the form $\{0\} \cup I_\varepsilon $, for all $\varepsilon >0$, \dash a neighborhood base for $\extra $ is given by the sets of the form $\{\extra \} \cup I_\varepsilon $, for all $\varepsilon >0$.

\medskip It is worth remarking that this topology fails to be Hausdorff since $\extra $ and $0$ cannot be separated by disjoint open subsets. Nevertheless one can easily prove that $G$ is a bona fide ample groupoid with two distinguished bisections $$ U_0 := X = G\ex 0, \and U_1 := X \setminus \{0\}\cup \{\extra \}. \reqno{(3.3)} $$

Fixing any field $\fld $, we will next construct a line bundle over $G$. As a first step, we let $E=\fld \times G$ be the trivial line bundle, as defined in (2.6). However we will ignore its topology, while retaining its algebraic structure, so we will so far view $E$ just as an algebraic Fell bundle.

In order to give $E$ a topology we will provide a fundamental family of local sections, which we will then feed into (3.2). Geometrically, what we intend to do will have the effect of taking the restriction of the (topologoical) trivial bundle over the bisection $U_1$ defined in (3.3), split it in two across the fiber over $\extra $, rotate the right-hand half by half a turn, and then glue it back together. Somewhat surprisingly, all of this will be done while leaving the restriction of $E$ to the bisection $U_0$ untouched, and in fact we are prohibited from messing it up by (2.11.ii).

For every open subset $U\subseteq G$ not containing $\extra $ (so that $U\subseteq X$), and for any locally constant function $\ell :U\to \fld $, we consider the local section $\xi _\ell $ defined on $U$ by $$ \xi _\ell (\gamma ) =\big (\ell (\gamma ), \gamma \big ), \for \gamma \in U. \reqno{(3.4)} $$

On the other hand, for every open subset $U\subseteq G$ containing $\extra $ but not zero, and hence necessarily also containing some $I_\varepsilon $, let $$ V=U\setminus \{\extra \}\cup \{0\}. $$ Observing that $V$ is an open subset of $X$ containng $0$, take any locally constant function $\ell :V\to \fld $, and define a local section $\eta _\ell $ on $U$ by $$ \eta _\ell (\gamma ) = \left \{ \matrix { (\ell (\gamma ), \gamma ), & \text { if $\gamma \in X\cap V$, and $\gamma <0$}, \cr \pilar {12pt} (\ell (0), \extra ), & \text { if } \gamma =\extra ,\hfill \cr \pilar {12pt} (-\ell (\gamma ), \gamma ), & \text { if $\gamma \in X\cap V$, and $\gamma >0$}. } \right . \reqno{(3.5)} $$

The reader should keep an eye on the minus sign above, as it is responsible for the special properties of the line bundle we are about to consider.

It is perhaps worth observing that the parameter $\ell $ used to define both $\xi _\ell $ and $\eta _\ell $, above, is always a locally constant function defined on an open subset of $X$, but, unlike $\xi _\ell $, the domain of $\eta _\ell $ is contained not in $X$, but in the bisection $U_1$ defined in (3.3).

We then let $\S $ be the collection of local sections consisting of all of the $\xi _\ell $ and all of the $\eta _\ell $, as defined above, and we leave it for the reader to show that $\S $ is a fundamental family of local sections for $E$. The only sligthly tricky point is the proof of (3.1.vi), which we believe is worth sketching: Thus, given two local sections $\zeta _1$ and $\zeta _2$ in $\S $, we must prove that the mini convolution product $\zeta _1*\zeta _2$ belongs to $\S $. Incidentally, notice that the domains of all members of $\S $ are open bisections.

If both $\zeta _1$ and $\zeta _2$ are of the form $\xi _\ell $, then their convolution product is just the pointwise product on their common domain, namely $\dom (\zeta _1)\cap \dom (\zeta _2)$, so $\zeta _1*\zeta _2$ clearly lies in $\S $.

If $\zeta _1=\xi _{\ell _1}$ and $\zeta _2=\eta _{\ell _2}$, for locally constant functions $\ell _1$ and $\ell _2$, a quick computation shows that $$ \zeta _1*\zeta _2= \zeta _2*\zeta _1=\eta _\ell , $$ where $\ell $ is the pointwise product of $\ell _1$ and $\ell _2$ on their commond domain.

The most interesting case is when $\zeta _1=\eta _{\ell _1}$, and $\zeta _2=\eta _{\ell _2}$, for locally constant functions $$ \ell _i:V_i\to \fld , $$ for $i=1,2$. In this case, both $V_1$ and $V_2$ are open subsets of $X$ containing zero, while $\dom (\zeta _i)=V_i\setminus \{0\}\cup \{\extra \}$. The domain of $\zeta _1*\zeta _2$, namely $\dom (\zeta _1){\cdot }\dom (\zeta _2)$, then turns out to be $V_1\cap V_2$ (which contains $0$, as the product of $\extra $ times itself), and we have for all $x$ in $V_1\cap V_2$, that $$ (\zeta _1*\zeta _2)(x) = \left \{ \def \quad { } \matrix { \big (\ell _1(x), x\big )\big (\ell _2(x), x\big ) &=& \big (\ell _1(x)\ell _2(x), x\big ), & \text { if $x<0$}, \cr \pilar {12pt} \big (\ell _1(0), \extra \big )\big (\ell _2(0), \extra \big ) &=& \big (\ell _1(0)\ell _2(0), 0\big ), & \text { if } x=0,\hfill \cr \pilar {12pt} \big ({-}\ell _1(x), x\big )\big ({-}\ell _2(x), x\big ) &=& \big (\ell _1(x)\ell _2(x), x\big ), & \text { if $x>0$}. } \right . $$ Therefore $\zeta _1*\zeta _2= \xi _{\ell _1\ell _2}\in \S $, as desired.

The reader might have noticed that the crux of the matter here is that $\extra $ is an involutive element of $G$ (meaning that its square is a unit), and so is the element $-1$ of the field $\fld $.

Employing (3.2), we then have that $E$ becomes a line bundle over $G$.

In line with the geometric motivation given above, observe that, for every $t$ in $\fld $, the sequence $(u_n)_n$ of elements in $E$ given by $u_n=(t, 1/n)$, for $n>0$, converges to $(-t, \extra )$. This is because $$ 1/n\to \extra $$ in $G$, so, regarding the continuous local section $\eta _\ell $ defined in (3.5), with $\ell $ the constant function $-t$, we have that $$ (t, 1/n) = \eta _\ell (1/n) \to \eta _\ell (\extra ) = (-t, \extra ). $$ This is in spite of the fact that the very same sequence $ (t, 1/n) $ converges to $(t, 0)$ as well, which is no surprise as neither $G$ nor $E$ are Hausdorff spaces.

As already mentioned, our motivation for introducing the example above is to show that the concept of line bundles is strictly more general than that of 2-cocycles, meaning that not all line bundles are of the form $E(\omega )$, for a cocycle $\omega $.

\state Proposition Every continuous global section of the line bundle $E$ defined above vanishes on the isotropy group $G(0)$.

\Proof Let $\zeta $ be a continuous global section of $E$. Then, by (3.2), we may pick some $\xi $ in $\S $, such that $\xi $ and $\zeta $ coincide on some neighborhood of $0$. Notice that $\xi $ must necessarily be of the form $\xi _\ell $ (defined in (3.4)), since none of the $\eta _\ell $ (defined in (3.5)) have $0$ in their domain.

Likewise there exists a locally constant function $\ell '$ defined on a neighborhood of $0$, such that $\eta _{\ell '}$ coincides with $\zeta $ on some neighborhood of $\extra $. Since we are allowed to arbitrarily shrink the neighborhoods of $0$ where $\ell $ and $\ell '$ are defined, we may assume that both $\ell $ and $\ell '$ are defined on a single neighborhood $V$ of $0$ where they are in fact constant. Thus, from now on we may view $\ell $ and $\ell '$ as members of $\fld $. Furthermore, we may assume that $\zeta $, $\xi _\ell $ and $\eta _{\ell '}$ are all defined and coincide on $V\setminus \{0\}$.

Choosing $n\in {\mathbb N}$ large enough so that both $1/n$ and $-1/n$ lie in $V$, we then have that $$ (\ell , -1/n) = \xi _\ell (-1/n) = \zeta (-1/n) = \eta _{\ell '}(-1/n) = (\ell ', -1/n), $$ so $\ell =\ell '$. Moreover $$ (\ell , 1/n) = \xi _\ell (1/n) = \zeta (1/n) = \eta _{\ell '}(1/n) = (-\ell ', 1/n). $$ so $\ell =-\ell '$. Consequently $\ell =0$, and thus $\zeta $ vanishes on both $0$ and $\extra $, as required. \endProof

We thus arrive at the reason for introducing this somewhat exotic line bundle.

\state Corollary The Fell bundle $E$ defined above is not isomorphic to any topological Fell bundle of the form $E(\omega )$, where $\omega $ is a 2-cocycle on $G$.

\Proof It is enough to notice that $E(\omega )$ admits nowhere vanishing continuous global sections, such as the one defined by $$ \xi (\gamma )=(1,\gamma ), \for \gamma \in G. \closeProof $$ \endProof

\section Twisted Steinberg algebras

As before we will deal here with \'etale groupoids, but, due to our interest in generalizing Steinberg algebras, we will soon also assume that our groupoids are ample, meaning that their unit space admits a (topological) basis of compact open subsets. This is also equivalent to the fact that the compact open bisections form a basis for the topology of the whole groupoid.

\state Definition \rm By a \emph {twisted groupoid} we shall mean a pair $(G,E)$, where $G$ is an \'etale groupoid, and $E$ is a line bundle over $G$.

In \cite {LisaOne} the same terminology has been used to refer to a Hausdorff, ample groupoid equipped with a continuous 2-cocycle. Although we borrow this terminology, the present concept is significantly more general, not only because we do not assume $G$ to be Hausdorff, but also because our \emph {twisting} is not only algebraic in nature, that is, arising from a 2-cocycle, but also topological, due to the possibility that the line bundle is not topologically trivial.

\fix Throughout this section we fix a field $\fld $, always seen as a discrete topological space, and a twisted groupoid $(G, E)$, such that $G$ is ample.

Our goal is to introduce an associative $\fld $-algebra $\Ak $, generalizing both Steinberg algebras and the twisted Steinberg algebras of \cite {BeckyOne}.

\state Definition \rm Given any open subset $U\subseteq G$, we will denote by $S(U)$ the set of all continuous local sections defined on $U$. Moreover, for every $\xi $ in $S(U)$, we will denote by $\tilde \xi $ the extension of $\xi $ to the whole of $G$, obtained by declaring it to be zero outside of $U$. Precisely $$ \tilde \xi (\gamma ) = \left \{ \matrix { \xi (\gamma ) , & \text { if $\gamma \in U$}, \cr \pilar {12pt} 0_\gamma , & \text { otherwise.} } \right . $$

Although the above definition makes perfect sense for all open sets $U$, we will almost exclusively use it when $U$ is a compact open bisection.

Assuming that $G$ is Hausdorff, and given a compact open bisection $U$, observe that $U$ is necessarily also closed, so it follows that $\tilde \xi $ is continuous on $G$ for every $\xi $ in $S(U)$. However, in the general (non Hausdorff) case, $U$ needs not be closed and hence $\tilde \xi $ might very well be discontinuous. An example of this phenomenon is the local section $\eta _\ell $, introduced in (3.5), with e.g. $\ell \equiv 1$, which admits no extension to a continuous global section by (3.6), because it does not vanish on $\extra $.

One of the technical reasons for assuming $G$ to be ample is the following:

\state Proposition Given a line bundle $E$ over an ample groupoid $G$, and given any $u$ in $E$, there exists a continuous local section $\xi $ passing through $u$, and such that $\dom (\xi )$ is compact.

\Proof Using (3.1.i), choose $\xi $, as above, except that $\dom (\xi )$ might not be compact. Since $G$ is ample, we may find a compact open neighborhood $V$ of $\pi (u)$, such that $V\subseteq \dom (\xi )$. The restriction of $\xi $ to $V$ then satisfies the required conditions. \endProof

\state Definition \rm Given a twisted groupoid $(G, E)$, and given any field $\fld $, we will let $\Ak $ be the linear subspace of the space of all global sections of $E$, linearly spanned by the union of all sets of the form $$ \tilde {S}(U) := \big \{\tilde \xi : \xi \in S(U)\big \}, $$ where $U$ ranges in the family of all compact open bisections of $G$.

We will soon make $\Ak $ into an algebra, equipping it with the multiplication operation defined as follows.

\state Definition \rm Given $f$ and $g$ in $\Ak $, we let $f*g$ be the global section defined by $$ (f*g)(\gamma ) = \sum _{\alpha \beta =\gamma }f(\alpha )g(\beta ). $$ We will say that $f*g$ is the \emph {convolution product} of $f$ and $g$.

The subscript ``$\alpha \beta =\gamma $'' under the summation sign, above, is meant to convey that the sum ranges over the set of all pairs $(\alpha ,\beta )$ in $G\ex 2$ such that $\alpha \beta =\gamma $.

\state Proposition The convolution product defined above is well defined and, with it, $\Ak $ becomes an associative $\fld $-algebra, henceforth called the \emph {twisted Steinberg algebra} associated to the twisted groupoid $(G, E)$.

\Proof We first check that the summation appearing in the definition of the convolution product has at most finitely many nonzero terms and moreover that the function so defined lies in $\Ak $.

For this it clearly suffices to treat the case $f=\tilde \xi $, and $g=\tilde \eta $, where $\xi \in S(U)$ and $\eta \in S(V)$, where $U$ and $V$ are compact open bisections. Given $\gamma \in G$, and given any pair $(\alpha ,\beta )\in G\ex 2$, with $\gamma =\alpha \beta $, leading up to a nonzero summand in said summation, notice that $f(\alpha )$ and $g(\beta )$ are nonzero, whence $\alpha \in U$, and $\beta \in V$. Moreover, since $U$ and $V$ are bisections, $(\alpha ,\beta )$ must be the unique factorization of $\gamma $, as defined near (2.7). In particular, we deduce that the above summation has at most one nonzero term, as required.

Still assuming that $f=\tilde \xi $, and $g=\tilde \eta $, as above, observe that for $\gamma \in U V$, $$ (f*g)(\gamma ) = f(\alpha )g(\beta ) = \xi (\alpha )\eta (\beta ) \explica {(2.7)}{=} \xi \Big ((r|_{U})^{-1}(r(\gamma ))\Big ) \eta \Big ((s|_{V})^{-1}(s(\gamma ))\Big ). $$ This shows that the function $\zeta $ defined by $$ \zeta :=(f*g)|_{U V}, $$ is continuous, and since $$ \zeta (\gamma ) = (f*g)(\gamma ) = f(\alpha )g(\beta ) \in E_\alpha E_\beta \subseteq E_{\alpha \beta } = E_\gamma , $$ we see that $\zeta $ is also a local section, meaning that $\zeta $ lies in $S(U V)$. Noticing that $UV$, is a compact open bisection, we deduce that $$ f*g=\tilde \zeta \in \Ak , $$ as claimed.

This shows that the convolution product is well defined on $\Ak $, and we leave it for the reader to show that it satisfies all of the axioms required to make $\Ak $ an associative $\fld $-algebra. \endProof

We should remark that the convolution product just defined extends the mini convolution product of (2.8), as made clear by the arguments employed in the above proof.

\section Regular inclusions

Much of what we have to say about Steinberg algebras substantially generalizes to regular inclusions, so we will first develop the general theory of regular inclusions before returning to the applications to Steinberg algebras we intend to offer.

We should say that the present section is very closely related to the study of regular inclusions found in \cite {ExelPitts} in the context of C*-algebras, and in fact it may be seen as a routine translation of the results found there to a purely algebraic context. As often happens when ideas arising from C*-algebras are applied in pure algebra, there is a high degree of overlap of results and often also of methods of proof, and such an overlap is definitely very prominent here.

As already mentioned in the introduction, the alternative of simply referring the reader to \cite {ExelPitts} didn't seem to us to be completely honest, so we have opted for the more conservative alternative of offering complete proofs, but the reader acquainted with \cite {ExelPitts}, and well trained in the interplay between pure algebra and C*-algebra, will perhaps prefer to skip to the next section.

\state Definition \rm Let $\fld $ be a field, let $B$ be a $\fld $-algebra, and let $A\subseteq B$ be a subalgebra. \item We shall say that $A$ is a \emph {left-} (resp.~\emph {right-}) \emph {s-unital subalgebra} of $B$ if, for every $b$ in $B$, there exists $u$ in $A$, such that $ub=b$ (resp.~$bu=b$). \item In case $A$ is both a left-s-unital subalgebra and a right-s-unital subalgebra, then $A$ is said to be an \emph {s-unital subalgebra}. \item If $A$ is an s-unital subalgebra of itself, then we say that $A$ is an \emph {s-unital algebra}.

\medskip If $A$ is an s-unital subalgebra of $B$, then the following apparently much stronger property holds: for every finite subset $F\subseteq B$, there exists $u$ in $A$, such that $$ bu=a=ub, \for b\in F. \reqno{(5.2)} $$ In case $A=B$, this result is proved in \cite [Proposition 2.10]{Nystedt}, but the proof given there also works in the general case.

\state Definition \rm Any element $u$ satisfying (5.2) will be called a \emph {dedicated unit} for $F$.

Observe that every s-unital algebra $A$ is idempotent, meaning that $A^2=A$, and if $A$ is an s-unital subalgebra of $B$, then then $AB=B=BA$. Here, the product of two linear subspaces $X$ and $Y$ of an algebra $A$ is defined by $$ XY = \Big \{\sum _{i=1}^n x_iy_i: n\in {\mathbb N}, \ x_i\in X,\ y_i\in Y\Big \}. \reqno{(5.4)} $$

\medskip The following is inspired on \cite {AraBosa}, which in turn is inspired on \cite {FeldmanMoore, Kumjian, RenaultCartan}.

\state Definition \rm Let $B$ be a $\fld $-algebra, and let $A$ be a subalgebra of $B$. \item An element $n\in B$ is said to be a \emph {normalizer} of $A$ in $B$, if there exists an element $n^*\in B$, such that $$ nn^*n=n,\quad n^*nn^*=n^*,\quad nAn^*\subseteq A, \and n^*An\subseteq A. $$ \item Any element $n^*$ satisfying the above condition will be called a \emph {partial inverse of $n$}. \item The set of all normalizers of $A$ in $B$ will be denoted by $N_B(A) $. \item We will say that $A$ is a \emph {regular subalgebra of $B$} if $A$ is an s-unital subalgebra of $B$, and the linear span of $N_B(A) $ coincides with $B$.

In order to be able to quickly refer to the situation we will be dealing with for most of the time, we make the following:

\state {Standing Hypotheses I} We will fix a field\/ $\fld $, a $\fld $-algebra $B$, and a regular subalgebra $A$ of $B$, and we will moreover suppose that: \iStyle a \item $A$ is abelian, \item $A$ is the linear span of the set of idempotent elements in $A$, \item the spectrum of $A$ will be denoted by $X$.

\bigskip Recall from \cite [Theorem IX.6.1]{Jacobson} (see also \cite [Theorem I]{Keimel}) that conditions (a) and (b), above, imply that $X$ is a totally disconnected, locally compact space, and $A$ is isomorphic to the algebra $L_c(X)$ formed by all locally constant, compactly supported, $\fld $-valued functions on $X$. In other words, (a), (b) and (c) of (5.6) may be replaced by

\bigskip \Item {($\bullet $)} $X$ is a totally disconnected, locally compact space, and $A=L_c(X)$. \bigskip

\noindent Accordingly we will adopt the notation $$ \evx a \reqno{(5.7)} $$ to represent the value of a given element $a\in A$, seen as a function on $X$, applied to a given point $x\in X$.

\state Proposition Assuming (5.6), one has that $N_B(A) $ is an inverse semigroup under multiplication, and moreover $A\subseteq N_B(A)$.

\Proof We first claim that, given any normalizer $n$, and any partial inverse $n^*$ of $n$, one has that $nn^*$ and $n^*n$ lie in $A$. This is because we may choose a dedicated unit $e$ for the set $\{n\}$, belonging to $A$, and then $$ nn^* = nen^* \in nAn^* \subseteq A, $$ and similarly $n^*n\in A$.

Given another normalizer $m$, with partial inverse $m^*$, we then deduce that $mm^*$ and $n^*n$ commute, so $$ nm(m^*n^*)nm = n(mm^*)(n^*n)m = n(n^*n)(mm^*)m = nm, $$ and similarly $m^*n^*(nm)m^*n^* = m^*n^*$. This proves that the first two conditions in (5.5.i) hold for $nm$, provided we choose $(nm)^*$ to be $m^*n^*$. The last two conditions are also easy to verify, so we see that $nm$ is again a normalizer. Summarizing, we have proved that $N_B(A) $ is a multiplicative subsemigroup of $B$.

Given any normalizer $n$, and given any partial inverse $n^*$ of $n$, it is clear that $n^*$ is also a normalizer, so $N_B(A) $ is a regular semigroup \cite [Page 6]{Lawson}. As observed above, all idempotents of the form $nn^*$ lie in $A$, which is commutative, so the $nn^*$ commute among themselves and this implies that $N_B(A) $ is an inverse semigroup by \cite [Page 6]{Lawson}, as desired.

Finally, to see that $A\subseteq N_B(A) $, let $a\in A$, and define a $\fld $-valued function $a^*$ on $X$ by $$ a^*(x) = \clauses { \cl {\evx a^{-1}}{\text {if } \evx a \neq 0,} \cl {0} {\text {otherwise.}} } $$ Clearly $a^*$ is a locally constant, compactly supported function, so we may view $a^*$ as an element of $A$, and clearly $$ aa^*a=a, \and a^*aa^*=a^*. $$ Since the remaining two conditions characterizing normalizers are obviously satisfied, we are done. \endProof

A simple consequence of the above is that the partial inverse $n^*$ of any given normalizer $n$ is unique, as this is well known to hold in any inverse semigroup.

Another consequence of (5.8) is that $B$ is linearly spanned by a multiplicative subsemigroup which happens to be an inverse semigroup. In fact we will now see that this provides an equivalent formulation of (5.6).

\state Proposition Let $B$ be an associative $\fld $-algebra, and let $S\subseteq B$ be a multiplicative subsemigroup. Suppose that \item $S$ is an inverse semigroup, and \item the linear span of $S$ coincides with $B$. \medskip \noindent Then, denoting by $E(S)$ the idempotent semilattice of $S$, and letting $$ A = \span \big (E(S)\big ), $$ one has that $A$ is a regular subalgebra of $B$, and the pair of algebras $(A, B)$ satisfies (5.6.a-b). Conversely, if $B$ is an associative algebra, and $A$ is a regular subalgebra of $B$ satisfying (5.6.a-b), then there is a multiplicative subsemigroup of $B$ satisfying (i) and (ii), such that $A = \span \big (E(S)\big )$.

\Proof Assuming the hypotheses of the first part of the statement, observe that $E(S)$ is also a multiplicative subsemigroup, so $A$ is a subalgebra of $B$. Moreover, if $s\in S$, and $e\in E(S)$, it is easy to see that $$ s^*es, \ ses^* \in E(S) \subseteq A, $$ from where one deduces that $s$ lies in $N_B(A)$. Therefore $B$ is spanned by $N_B(A)$.

To see that $A$ is a left s-unital subalgebra of $B$, pick $b\in B$, and write $$ b = \sum _{i=1}^n \lambda _is_i, $$ whith the $\lambda _i$ in $\fld $, and the $s_i$ in $S$. Putting $e_i=s_is_i^*$, we have that $e_i\in E(S)$, so $e_i$ belongs to $A$.

We then claim that there exists an idempotent element $f$ in $A$, such that $$ fe_i=e_i, \for i\in \{1,2,\ldots ,n\}. $$ Using induction, suppose that an idempotent element $f'$ has been found in $A$, satisfying $f'e_i=e_i$, for every $i\leq n-1$. It is then easy to see that $$ f:= e_n+f'-e_nf' $$ satisfies all of the required conditions. Observing that $$ fs_i = fs_is_i^*s_i = fe_is_i = e_is_i = s_i, $$ we see that $fb=b$, thus proving that $A$ is a left s-unital subalgebra of $B$. A similar reasoning proves that $A$ is also a right s-unital subalgebra. The conclusion is then that $A$ is regular in $B$. In order to conclude the proof of the first assertion in the statement we still need to verify (5.6.a-b), but this is evident.

Focusing on the second part of the statement, and hence assuming (5.6), we take $S= N_B(A)$, and all we need to do is prove that $A$ coincides with the linear span of $E(S)$. Given $e$ in $E(S)$, as seen in the proof of (5.8) we have that $$ e = ee^* \in eAe^* \subseteq A, $$ so $E(S)\subseteq A$, and hence also $\span \big (E(S)\big )\subseteq A$. To prove the reverse inclusion, notice that $A$ is spanned by its idempotent elements, by assumption, and it is clear that every such idempotent element lies in the idempotent semilattice of $N_B(A)$. \endProof

Given any element $a$ in $A$, we will denote by $\supp (a)$ the \emph {support} of $a$, namely $$ \supp (a)=\{x\in X:\ev ax\neq 0\}, $$ observing that this is a compact open subset of $X$.

The following is an elementary adaptation of a well known result proved by Kumjian \cite {Kumjian} in the context of C*-algebras.

\state Proposition Assuming (5.6), and given $n\in N_B(A) $, there exists a homeomorphism $$ \beta _n:\src (n) \to \tgt (n), $$ where $$ \src (n) = \supp (n^*n), \and \tgt (n) = \supp (nn^*), $$ such that, for every $a$ in $A$, one has that $$ \ev {n^*an}{x} = \ev {a}{\beta _n(x)}, \for x\in \src (n). $$ Moreover $\beta _{n^*}=\beta _n^{-1}$, and, if $m$ is another normalizer, then $\beta _{nm}=\beta _n\circ \beta _m$, where the composition is defined on the largest possible domain, as is standard for partially defined maps.

\Proof Given $x$ in $\src (n)$, define a linear functional $\varphi $ on $A$ by $$ \varphi (a) = \ev {n^*an}{x}, \for a\in A. $$ It is easy to see that $\varphi $ is multiplicative and $\varphi (nn^*)\neq 0$, so $\varphi $ is a character on $A$, whence there is a unique element in $X$, which we denote by $\beta _n(x)$, such that $$ \varphi (a) = \ev {a}{\beta _n(x)}, \for a\in A. $$

Due to the fact that $\varphi (nn^*)\neq 0$, we have that $\beta _n(x)\in \tgt (n)$. To see that $\beta _n$ is continuous, recall that the topology on $X$ is the \emph {initial topology} determined by the functions $$ x\in X\mapsto \evx a\in \fld , $$ for all $a$ in $A$, where $\fld $ is, as always, viewed as a discrete topological space. This means that, if $Y$ is any topological space, and $$ f:Y \to X $$ is any function, then $f$ is continuous if and only if $$ y\in Y\mapsto \ev a{f(y)}\in \fld , $$ is continuous for all $a$ in $A$. The continuity of $\beta _n$ then easily follows from this criterion. Finally, to see that $\beta _n$ is a homeomorphism, it is enough to observe that $\beta _{n^*}$ is its inverse.

Choosing two normalizers $n$ and $m$, and given $x$ in $\src (nm)$, we have that $$ 0\neq \evx {m^*n^*nm} = \evx {m^*mm^*n^*nm}= \evx {m^*m} \evx {m^*n^*nm}, $$ so $x\in \src (m)$. Moreover, $$ 0\neq \evx {m^*n^*nm} = \ev {n^*n}{\beta _m(x)}, $$ so $\beta _m(x)\in \src (n)$, and we have for all $a$ in $A$ that $$ \ev a{\beta _n(\beta _m(x))} = \ev {n^*an}{\beta _m(x)} = \ev {m^*n^*anm}{x} = \ev {a}{\beta _{nm}(x)}. $$ This implies that $\beta _n(\beta _m(x)) = \beta _{nm}(x)$, so we see that $\beta _n\circ \beta _m$ extends $\beta _{nm}$.

On the other hand, if $x$ lies in the domain of $\beta _n\circ \beta _m$, that is, if $x\in \src (m)$, and $\beta _m(x)\in \src (n)$, then $$ 0\neq \ev {n^*n}{\beta _m(x)} = \ev {m^*n^*nm}{x}, $$ so $x\in \src (nm)$. This concludes the proof. \endProof

The maps introduced above allow for the following very usefull classification of normalizers:

\state Definition \rm Given $x$ and $y$ in $X$, we let \iStyle a \item $N_x =\big \{n\in N_B(A) : x\in \src (n) \}$, \item $\N yx =\big \{n\in N_x : \beta _n(x)=y\big \}$, \item $\Orb x= \{\beta _n(x): n\in N_x\}$.

\bigskip Some properties relating to the notions just defined are as follows:

\state Proposition For every $x$, $y$ and $z$ in $X$, we have that: \item $N_x= \bigcup _{y\in \Orb x} \N yx$, \item $y\in \Orb x$ if and only if $\N yx\neq \emptyset $, \item $\N zy\N yx\subseteq \N zx$, \item $(\N yx)^*=\N xy$, \item if $y\neq z$, then $\big (\N xz\N yx\big ) \cap \N xx = \emptyset $.

\Proof Before we begin, we should perhaps point out that the convention regarding product of sets made in (5.4) does not apply to the product of sets appearing in (iii) and (v), above, due to the fact that neither of these are linear spaces. Thus, when we speak of, e.g.~$\N zy\N yx$ in (iii), we are referring simply to the set of products, rather than the linear span of these.

Returning to the proof, everything is pretty elementary except perhaps for (v), which we prove as follows: suppose by contradiction that $n\in \N xz$, $p\in \N yx$, and $q\in \N xx$ are such that $np=q$. Then $$ n^*np = n^*q \in \N zx \N xx \subseteq \N zx, $$ so $$ z = \beta _{n^*np}(x) = \beta _{n^*n}\big (\beta _p(x)\big ) = \beta _{n^*n}(y), $$ but $\beta _{n^*n}=\beta _n^{-1}\circ \beta _n$ is the identity map on $\src (n)$, so it cannot map $y$ to $z$. \endProof

Let us now temporarily put ourselves in a slightly different context, by assuming instead that:

\state {Standing Hypotheses II} \rm \item $B$ is a fixed $\fld $-algebra and $A$ is an s-unital subalgebra of $B$. \item $I$ and $ J$ are fixed ideals of $A$. \item $I$ and $ J$ will be assumed to be s-unital algebras.

Obviously neither $I$ nor $J$ need to be ideals of $B$, but our goal is to discuss some simple, but crucial properties describing the interactions between these ideals and $B$.

\state Lemma \item For every $b\in IB$, there exists $u$ in $I$, such that $ub=b$. \item For every $b\in BJ$, there exists $v$ in $J$, such that $bv=b$. \item $IBJ= IB\cap BJ$.

\Proof (i) Given $b \in IB$, write $$ b = \sum _{i=1}^n x_ib_i, $$ with $x_i\in I$, and $b_i\in B$. Since $I$ is s-unital, there exists $u$ in $I$ such that $ux_i=x_i$, for every $i$, so $ub=b$. This proves (i), while (ii) is proven in a similar way.

Regarding (iii), let $b\in IB\cap BJ$, and use (i) and (ii) to find $u$ in $I$ and $v$ in $J$ such that $ub = b = bv$. Then $b = ubv\in IBJ$. The reverse inclusion is obvious. \endProof

The following introduces a crucial ingredient leading up to the notion of isotropy algebras.

\state Proposition Under the assumptions in (5.13), and setting $$ \C IJ = \big \{c\in B: cJ\subseteq IB, \ Ic\subseteq BJ\big \}, $$ one has that \item $A\subseteq \C JJ$, \item $ IJ\subseteq I\C IJ = \C IJ J= IBJ\subseteq \C IJ, $ \item $\C IJ \cap IB = IBJ = \C IJ \cap BJ$.

\Proof (i) Obviously $A\subseteq \C JJ $ and $IBJ\subseteq \C IJ $.

\medskip \noindent (ii) Given $x$ in $I$ and $c$ in $\C IJ $, observe that $$ xc\in Ic\subseteq BJ, $$ so we see that $I\C IJ \subseteq BJ$. Since one obviously has $I\C IJ \subseteq IB$, it follows that $$ I\C IJ \subseteq IB\cap BJ\explica {(5.14.iii)}{=} IBJ. $$

Having already proven that $$ I\C IJ \subseteq IBJ\subseteq \C IJ , $$ we may left-multiply everything by $I$ to get $$ I\C IJ \subseteq IBJ\subseteq I\C IJ , $$ whence $I\C IJ = IBJ$. The proof that $\C IJ J= IBJ$ is done along similar lines.

Since $IJ\subseteq IBJ$, the last remaining inclusion in (ii) is proved.

\medskip \noindent (iii) From (ii) it follows that $IBJ$ is contained in both $\C IJ \cap BJ$ and in $\C IJ \cap IB $. In order to prove that $\C IJ \cap BJ \subseteq IBJ$, pick $c\in \C IJ \cap BJ$. Using (5.14.ii), we may find $v$ in $J$ such that $c = cv$, so $$ c = cv \in \C IJ J \explica {(ii)}{=} IBJ. $$ The proof that $\C IJ \cap IB \subseteq IBJ$ goes along similar lines. \endProof

\state Definition \rm The \emph {isotropy module} of the inclusion ``$A\subseteq B$" at the pair of ideals $(I, J) $ is the quotient vector space $$ \B IJ = \frac {\C IJ }{\H IJ }. $$ where $$ \H IJ := I\C IJ = \C IJ J= IBJ. \reqno{(5.16.1)} $$ We will denote the quotient map by $$ \p IJ : \C IJ \to \B IJ . $$ Finally, in case $I=J$, we will refer to $\B JJ$ it as the \emph {isotropy algebra} of the inclusion ``$A\subseteq B$" at the ideal $J$.

\state Proposition If, in addition to the ideals $I$ and $J$ of (5.13), we are given a third ideal $K\trianglelefteq A$, also assumed to be an s-unital algebra, then \item $\C IJ \C JK \subseteq \C IK $, \item there exists a bilinear map (denoted simply by juxtaposition) $$ \B IJ \times \B JK \to \B IK , $$ such that $$ \p IJ (a) \p JK (b) = \p IK (ab), \for a\in \C IJ , \for b\in \C JK . $$

\Proof (i)\enspace Given $c_1\in \C IJ $ and $c_2 \in \C JK $, we have $$ c_1c_2K\subseteq c_1JB\subseteq IBB\subseteq IB, $$ and similarly $$ Ic_1c_2 \subseteq BJc_2 \subseteq BBK\subseteq BK, $$ so $c_1c_2$ lies in $\C IK $.

\itmProof (ii) Multiplying both sides of (i) on the left by $I$, and then again on the right by $K$ gives, $$ I\C IJ \C JK \subseteq I\C IK , \and \C IJ \C JK K\subseteq \C IK K= I\C IK . $$

Denoting by ``$\mu $'' the restriction of the multiplication operation of $B$, we then see that the composition of maps $$ \C IJ \times \C JK \labelarrow \mu \C IK \labelarrow {\p IK} \frac {\C IK } {I\C IK } = \B IK , $$ vanishes on $$ \H IJ \times \C JK , \ \ \hbox {and \ on}\ \ \C IJ \times \H JK, $$ so it factors through the cartesian product of quotient spaces $$ \frac {\C IJ }{\H IJ } \times \frac {\C JK }{\H JK } = \B IJ \times \B JK , $$ producing the required bi-linear map. \endProof

In the special case that $I=J=K$, the above result implies that $\C JJ$ is a $\fld $-algebra, and (5.15) says that that $\H JJ$ is a two-sided ideal in $\C JJ$. Consequently the isotropy algebra $\B JJ$ is indeed a $\fld $-algebra.

\state Lemma Under (5.13), set $$ \L IJ = IB+ BJ. $$ Then: \item given $ b$ in $ B$, we have that $ b\in \L IJ, $ if and only if there are elements $u$ in $I$ and $v$ in $J$, such that $ b= ub+ bv - ubv. $ \item $ \H IJ = \C IJ \cap \L IJ$.

\Proof (i) The ``if'' part being trivial, we focus on the ``only if'' part. Given that $b\in \L IJ$, write $b = c + d$, with $c\in IB$, and $d\in BJ$. Using (5.14.i) and (5.14.ii), choose $u$ in $I$ and $v$ in $J$ such that $uc = c$, and $dv=d$. Then $$ ub+ bv - ubv = u(c+d)+ (c+d)v - u(c+d)v = $$ $$ = c+ud+ cv+d - cv - ud = c + d = b. $$

\medskip \noindent (ii) From (5.15) we have that $\H IJ\subseteq \C IJ$, and it is clear that $\H IJ\subseteq \L IJ$. On the other hand, given $c$ in $\C IJ \cap \L IJ $, we have by (i) that there are elements $u$ in $I$ and $v$ in $J$, such that $ c= uc+ cv - ucv, $ so $$ c= uc+ cv - ucv \in I\C IJ + \C IJ J + I \C IJ J \explica {(5.15.ii)}{\subseteq} IBJ. \closeProof $$ \endProof

\state Proposition Besides the map $\p IJ $ of (5.16), let $$ \q IJ :B\to \frac {B} {\L IJ }, $$ denote the quotient map. Then there is an injective linear map $$ \Psi :\B IJ \to \frac {B}{\L IJ }, $$ such that $\Psi \big (\p IJ (c)\big ) = \q IJ (c)$, for all $c$ in $\C IJ $.

\Proof Since $\H IJ \subseteq \L IJ $, the composition $$ \C IJ \hookrightarrow B\ {\buildrel \q IJ \over \longrightarrow }\ \frac {B}{\L IJ } $$ vanishes on $\H IJ $, and therefore it factors through the quotient providing a linear map $\Psi $, as in the statement, such $$ \Psi \big (\p IJ (c)\big ) = \q IJ (c), \for c\in \C IJ . $$ In order to prove that $\Psi $ is injective suppose that $c\in \C IJ$ and $\Psi \big (\p IJ (c)\big ) =0$. Then $c$ lies in $\L IJ$, whence (5.18.ii) implies that $ c \in \H IJ, $ and it follows that $\p IJ (c)=0$. \endProof

We would now like to discuss a property for pairs of ideals which will be of fundamental importance later on.

\state Definition \rm We shall say that $(I, J)$ is a \emph {regular pair of ideals} if the map $\Psi $ introduced in (5.19) is surjective. In case $(J, J)$ is a regular pair, we shall simply say that $J$ is a \emph {regular ideal}.

A useful characterization of regular pairs is as follows:

\state Proposition $(I, J) $ is regular if and only if $$ B= \C IJ + \L IJ . $$ In this case $\B IJ \simeq B/\L IJ$, as vector spaces, through the map $$ \p IJ(c) = c+\H IJ \in \B IJ \mapsto \q IJ(c) = c+\L IJ \in B/\L IJ. $$

\Proof Given $b$ in $B$ notice that $\q IJ (b)$ lies in the range of $\Psi $ if and only if there exists $c$ in $\C IJ $ such that $\q IJ (b) = \q IJ (c)$, which in turn is equivalent to saying that $b\in \C IJ + \L IJ $. The last sentence in the statement is simply restating that $\Psi $ is surjective. \endProof

\def \C #1#2{\SimpleSuber C{#1}{#2}} \def \H #1#2{\SimpleSuber H{#1}{#2}} \def \B #1#2{\SimpleSuber B{#1}{#2}} \def \L #1#2{\SimpleSuber L{#1}{#2}}

Assuming that $(I,J)$ is regular, we may define a map $$ \E IJ:B\to \B IJ, $$ by $\E IJ = \Psi ^{-1}\circ \q IJ$, which will then satisfy $$ \E IJ(c) = \p IJ(c), \for c\in \C IJ. \reqno{(5.22)} $$

\state Definition \rm The map $\E IJ$ defined above will be called the \emph {isotropy projection} associated to the pair of ideals $(I,J)$.

\state Proposition Under (5.13), suppose that $I$, $J$, and $K$ are ideals of $A$, such that $(I, J)$ and $(J, K)$ are regular pairs. Then, given $$ g \in \C IJ,\quad b \in B, \and h \in \C JK, $$ one has that \item $\p IJ(g)\E JK(b) = \E IK(gb)$, \item $\E IJ(b)\p JK(h) = \E IK(bh)$.

\Proof In order to prove (i), observe that $$ B=\C JK + \L JK = \C JK + JB + BK, $$ due to (5.21), so we may write $$ b=c+a_1b_1+b_2a_2, $$ with $c\in \C JK$, $a_1\in J$, $a_2\in K$, and $b_1, b_2\in B$. Then $$ \p IJ(g)\E JK(b) = \p IJ(g)\E JK(c) = \p IJ(g)\p JK(c) \explica {(5.17.ii)}{=} \p IK(gc). $$ On the other hand, $$ gb = gc+ga_1b_1+gb_2a_2, $$ and we observe that $$ gc\in \C IJ \C JK \subseteq \C IK $$ $$ ga_1b_1\in \C IJ J B \subseteq I B B \subseteq IB, $$ and clearly $gb_2a_2\in BK$. So $$ \E IK(gb) = \E IK(gc) = \p IK(gc). $$ This proves (i), and (ii) may be proved similarly. \endProof

Returning to the context outlined in (5.6), recall that $A$ is identified with $L_c(X)$, where $X$ is the spectrum of $A$. We would now like to specialize our study of ideals to those given by the kernel of characters.

\state Definition \rm Assuming the situation of (5.6), for each $x$ in $X$ we will denote by $J_x$ the ideal of $A$ given by $$ J_x= \{a\in A: \evx a= 0\}. $$ If $y$ is another point in $X$, we will also let $$ \def \crr {\vrule height 10pt depth 1pt width 0pt \cr } \matrix { \C yx & = & \C {J_y}{J_x},\crr \H yx & = & \H {J_y}{J_x},\crr \B yx & = & \B {J_y}{J_x},\crr } \hskip 2cm \matrix { \L yx & = & \L {J_y}{J_x},\crr \p yx & = & \p {J_y}{J_x},\crr \q yx & = & \q {J_y}{J_x}.\crr } \reqno{(5.25.1)} $$ In case $x=y$, we will denote the four sets above simply by $\C xx $, $\H xx $, $\B xx $ and $\L xx $ respectively.

\state Proposition For every $x$ in $X$, one has that \item $J_x$ is an s-unital algebra. \item Given $a\in A$, let $\lambda =\evx a$. Then $$ \lambda b-ab\in J_xb, \and \lambda b-ba\in bJ_x, $$ for all $b$ in $B$.

\Proof (i) Given $ a\in J_x, $ observe that $ U:= \supp (a) $ is a compact open subset of $X$ not containing $x$. Therefore $\1U$ lies in $J_x$, and clearly $\1Ua=a$.

\itmProof (ii) Given $a$ and $b$ as in the statement, choose $u$ in $A$ such that $b=ub$. Then $$ \lambda b-ab = \lambda ub-aub = (\lambda u-au)b \in J_xb, $$ because $\lambda u-au$ lies in $J_x$. Similarly one proves that $\lambda b-ba\in bJ_x$.

\endProof

The next technical Lemma will be of crucial importance. Its purpose is to understand the relationship between the position of a given point $x$ relative to $\src (n)$, on the one hand, and the interplay between the normalizer $n$ and the various subspaces of $B$ determined by $J_x$, on the other.

\state Lemma Under (5.6), given $x$ and $y$ in $X$, and given a normalizer $n$ in $N_B(A) $, one has that: \item If $x\notin \src (n)$, then $n\in BJ_x$. \item If $x\in \src (n)$, then $$ J_{\beta _n(x)} n\subseteq BJ_x, \and nJ_x\subseteq J_{\beta _n(x)}B. $$ \item If $x\in \src (n)$, and $\beta _n(x)\neq y$, then $n\in \L yx $.

\Proof (i) If $\ev {n^*n}{x} = 0$, then $n^*n$ belongs to $J_x$, so $$ n= n(n^*n) \in BJ_x. $$

\itmProof (ii) For every $a$ in $J_{\beta _n(x)}$, notice that $$ \ev {n^*an}{x} = \ev {a}{\beta _n(x)} = 0, $$ so we see that $n^* an\in J_x$, whence $$ an= ann^*n = nn^*an \in BJ_x. $$ On the other hand, if $a$ lies in $J_x$, then $$ 0= \ev {n^*n}{x} \ev {a}{x} \ev {n^*n}{x} = \ev {n^*nan^*n}{x} = \ev {nan^*}{\beta _n(x)}, $$ so $nan^*\in J_{\beta _n(x)}$, whence $$ na = nn^*na = nan^*n \in J_{\beta _n(x)}B. $$

\itmProof (iii) Assuming that $y\neq \beta _n(x)$, we may choose $v_1$ and $v'_2$ in $A$, such that $$ \ev {v_1}y= 1, \quad \ev {\pilar {9pt}v'_2}{\beta _n(x)} = 1, \and v_1v'_2 = 0. $$ Setting $ v_2 = n^*v'_2n, $ we then have $$ \ev {v_2}{x} = \ev {n^*v'_2 n}{x} = \ev {v'_2}{\beta _n(x)} = 1, $$ and $$ v_1nv_2 = v_1nn^*v'_2n= v_1v'_2nn^*n= 0. $$ Observing that $ n-v_1n\in J_yB$, and $n-nv_2 \in BJ_x, $ by (5.26.ii), we have that $$ n = n-v_1n + v_1(n-nv_2) \in J_yB+ BJ_x. \closeProof $$ \endProof

We now have the tools to prove a main result.

\state Theorem With the hypotheses of (5.6), let $x, y\in X$. Then \item $\N yx \subseteq \C yx $, \item for every normalizer $n$ not belonging to $\N yx$, one has that $\q yx (n)=0$. \item if $N$ is a subset of $N_B(A) $ spanning $B$, then $\B yx $ is spanned by $\p yx (N\cap \N yx )$, \item $(J_y, J_x) $ is a regular pair of ideals.

\Proof The first point follows immediately from (5.27.ii). Assuming that $n$ is not in $\N yx$, at least one of the conditions defining $\N yx $ must fail, so we suppose first that $x$ is not in $\src (n)$, which is to say that $\evx {n^*n}=0$. Then by (5.27.i) we have that $n\in BJ_x\subseteq \L yx $, so $\q yx (n)=0$, as desired. If, on the other hand, $n$ lies in $\src (n)$, then necessarily $\beta _n(x)\neq y$, so (5.27.iii) gives that $ n\in \L yx , $ whence $\q yx (n) = 0$.

Focusing now on point (iv) of the statement, we need to show that the map $$ \Psi :\B yx = \frac {\C yx }{\H yx } \ \longrightarrow \ \frac {B}{\L yx }, $$ introduced in (5.19) is surjective.

Let $N\subseteq N_B(A) $ be any subset spanning in $B$, as in (iii). We then have that $\q yx (B)$ is spanned the set of all $\q yx (n)$, as $n$ run in $N$. It is then enough to show that $\q yx (n)$ lies in the range of $\Psi $, for every $n\in N$. But this follows easily since $$ \q yx (n) = \Psi (\p yx (n)), \for n\in N\cap \N yx, $$ by (i), whereas $$ \q yx (n)=0, \for n\in N\setminus \N yx, $$ by (ii).

Observe that the argument above in fact shows that $B/\L yx $ is generated by $\q yx (N\cap \N yx )$. Since we now know that $\Psi $ is a linear isomorphism from $\B yx $ onto $B/\L yx $, it follows that $\B yx $ is generated by $$ \Psi \inv \big (\q yx (N\cap \N yx )\big ) = \p yx (N\cap \N yx ). $$ This proves (iii). \endProof

Now that we know that $(J_y,J_x)$ is regular, the isotropy projection $\E {J_y}{J_x}$ defined in (5.23) becomes available. In line with (5.25) we make the following:

\state Definition \rm Under (5.6), and given $x,y\in X$, we will denote the isotropy projection $\E {J_y}{J_x}$ simply by $\E yx$.

Even though we are not assuming either $A$ or $B$ to be unital algebras, the isotropy algebras are always unital as we will now prove.

\state Lemma Assume (5.6) and let $x, y\in X$. Then, \item for every $h\in \B yx$, one has that $$ \p yy(a)h = \evy a h, \and h\p xx(a) = \evx a h, $$ \item $\B xx $ is a nonzero unital $\fld $-algebra, \item for all $a$ in $A$, one has that $\p xx(a) = \evx a 1_x$, where $1_x$ denotes the unit of $\B xx$.

\Proof Choosing $c$ in $\C yx$, such that $\p yx (c)=h$, we have that $$ \evy ac - ac \explica {(5.26.ii)}{\in} J_yc \subseteq J_y\C yx \explica {(5.16.1)}{=} \H yx, \reqno{(5.30.1)} $$ so $$ \p yy (a) h = \p yy (a) \p yx(c) \explica {(5.17.ii)}{=} \p yx (ac) \explica {(5.30.1)}{=} \evy a \p yx (c) = \evy a h, $$ proving the first equation in (i), while the second one follows by a similar argument.

Choosing any $a_0$ in $A$ such that $\evx {a_0}=1$, we deduce that $$ \p xx (a_0 )h = h = h\p xx (a_0 ), \for h\in \B xx, $$ so it follows that $\p xx (a_0)$ is the unit of $\B xx $, proving (ii). The last point is now evident, so it only remains to prove that $\B xx\neq \{0\}$. Arguing by contradiction, suppose that $\B xx=\{0\}$, meaning that $\C xx=\H xx$, so in particular $$ A\subseteq \H xx=J_xBJ_x\subseteq J_xB. $$ Fixing any $a$ in $A$, we may then write $a = ua$, with $u$ in $J_x$, thanks to (5.14.i). This yields $$ \evx a= \evx {ua} = \evx {u} \evx {a} = 0, $$ because $u\in J_x$, but this is a contradiction due to the fact that one may easily produce $a$ in $A$ with $\evx a=1$. \endProof

\section The imprimitivity bimodule

In this section we will again work in the situation described in (5.6), recalling that $A$ is identified with $L_c(X)$. Our purpose is to define and study the all important bimodules which we will later use as the key ingredient of the induction process.

\state Proposition Letting $x\in X$, and viewing $B$ as a right $\C xx$-module, one has that $BJ_x$ is a sub-module.

\Proof By definition $J_x\C xx\subseteq BJ_x$, so $$ BJ_x\C xx \subseteq BBJ_x \subseteq BJ_x. \closeProof $$ \endProof

As a consequence, the quotient space $$ M_x:= B/BJ_x $$ is a right $\C xx$-module. We next observe that $M_x$ is annihilated by the ideal $$ \H xx=J_xBJ_x\trianglelefteq \C xx, $$ meaning that $M_x\H xx=\{0\}$, so we may view $M_x$ as a right module over $\C xx/\H xx$, namely $\B xx$.

Obviously $M_x$ is also a left $B$-module and it is clear that this module structure is compatible with the right $\B xx$-module structure defined above in the sense that $M_x$ is a $B$-$\B xx$-bimodule.

\state Definition \rm The $B$-$\B xx$-bimodule $M_x$ described above will be called the \emph {imprimitivity bimodule}. If $V$ is any left module over $\B xx$ then the left $B$-module constructed via the familiar tensor product construction $$ M_x\tprodx V $$ will be said to be \emph {induced} from $V$, and it will be denoted by $\Indx V$.

Here is an important remark about our notation: if $U\subseteq M_x$, and $W\subseteq V$ are linear subspaces, then $U\otimes W$ will denote the linear subspace of $M_x\tprodx V$ given by $$ U\otimes W = \span \{u\otimes w: u\in U, w\in W\}. \reqno{(6.3)} $$ Incidentally notice that this is not necessarily equal to $U\tprodx W$, which indeed makes no sense at all unless $U$ is a right $\B xx$-module and $W$ is a left $\B xx$-module, conditions that may or may not be present in what follows.

Recall from (5.15.iii) that $$ \C xx\cap BJ_x = J_xBJ_x = \H xx, $$ so the composition of the inclusion map of $\C xx$ into $B$, with the quotient map from $B$ onto $B/BJ_x$, namely $$ \C xx\hookrightarrow B \to \frac B{BJ_x}, $$ factors through the quotient, providing an injective, $\fld $-linear map $$ \mu _x : \B xx = \frac {\C xx}{\H xx} \to \frac B{BJ_x} = M_x, $$ such that $$ \mu _x(c+\H xx) = c+BJ_x, \for c\in \C xx. \reqno{(6.4)} $$

\state Definition \rm The map $\mu _x$ defined above will be called the \emph {standard inclusion}.

\state Proposition The standard inclusion is a right $\B xx$-module homomorphism.

\Proof Obvious. \endProof

Given the relevance of $M_x$ in what follows, and given that $B$ is spanned by normalizers, it will be important to understand the image of a given normalizer under the quotient map from $B$ to $M_x$. But before that we need to introduce the following useful notation: $$ \Nbd = \{U\subseteq X: U \text { is a compact open neighborhood of }x\}. $$ We should also notice that, for every $U$ in $\Nbd $, the characteristic function of $U$, denoted $\1U$, is a locally constant, compactly supported function on $X$, so $\1U$ may be seen as an element of $A$.

\state Lemma Let $n$ be a normalizer in $N_B(A) $, and let $\xi =n+BJ_x$. \item If $x\notin \src (n)$, then $\xi =0$. \item If $x\in \src (n)$, and $\beta _n(x) \neq x$, then $\1U\xi =0$, for some $U\in \Nbd $. \item If $x\in \src (n)$, and $\beta _n(x) = x$, then $\1U\xi =\xi $, for every $U\in \Nbd $.

\Proof (i) As seen in (5.27.i), in case $x\notin \src (n)$, we have that $n\in BJ_x$, so that $n+BJ_x$ vanishes.

\itmProof (ii) If $x\in \src (n)$, and $\beta _n(x) \neq x$, let $U$ be a compact open neighborhood of $x$ not containing $y:=\beta _n(x)$. Observing that $n\in \N yx\subseteq \C yx$, by (5.28.i), and that $\1U$ belongs to $J_y$, we conclude that $$ \1Un \in J_yn \subseteq BJ_x, $$ so $\1U\xi =0$.

\itmProof (iii) If $x\in \src (n)$, and $\beta _n(x) = x$, then $n\in \C xx$, and, choosing any compact open neighborhood $U$ of $x$, we have that $$ \1Un-n \explica {(5.26.ii)}{\in} J_xn \subseteq BJ_x, $$ so $\1U\xi =\xi $. \endProof

\state Corollary Recalling from (5.11) that $$ N_x = \{n\in N_B(A) : x\in \src (n) \}, $$ one has that $M_x$ is linearly spanned by the set $$ \{n+BJ_x: n\in N_x\}. $$

\Proof Follows immediately from (6.7) and the fact that $B$ is linearly spanned by $N_B(A) $. \endProof

We have already mentioned that $M_x$ is a left $B$-module, so it is also a left $A$-module. In what follows this left $A$-module structure will play an important role, so it is worth discussing it further.

\state Proposition Given a normalizer $n\in N_x$, one has that $$ an+BJ_x=\ev a{\beta _n(x)}n+BJ_x, \for a\in A. $$

\Proof Given $a$ in $A$, we have that $$ an - \ev a{\beta _n(x)}n \explica {(5.26.ii)}{\in} J_{\beta _n(x)}n \subseteq BJ_x, $$ so $$ an+BJ_x = \ev a{\beta _n(x)}n+BJ_x, $$ as desired. \endProof

The following focuses on studying the range of the standard inclusion.

\state Proposition For $x\in X$, the following hold: \item $ \ds \mu _x(\B xx) = \bigcap _{U\in \Nbd } \1UM_x. $ \item There exists a unique idempotent, left $A$-linear map $\pi _x:M_x\to M_x$, whose range coincides with $\mu _x(\B xx)$. \item For all $\xi $ in $M_x$, there exists $U$ in $\Nbd $, such that $$ \pi _x(\xi )=\1V\xi , $$ for every $V$ in $\Nbd $, with $V\subseteq U$. \item $\pi _x$ is right $\B xx$-linear. \item For every $n\in N_B(A) $, one has that $$ \pi _x(n+BJ_x) = \clauses { \cl {0} {\text {if } x\notin \src (n), } \cl {0} {\text {if } x\in \src (n), \text { and } \beta _n(x)\neq x, } \cl {n+BJ_x}{\text {if } x\in \src (n), \text { and } \beta _n(x)=x. } } $$

\Proof Viewing $\Nbd $ as a downward directed set, we claim that, for every $\xi $ in $M_ x$, the net $$ \{\1U\xi \}_{U\in \Nbd } $$ is eventually constant, in the sense that there exists $U$ in $\Nbd $, such that $$ \1U\xi = \1V\xi , $$ for all $V\in \Nbd $, with $V\subseteq U$. In view of (6.8), in order to prove this claim, we may suppose that $\xi =n+BJ_x$, for some $n\in N_x$, but then the claim follows immediately from (6.7.ii-iii). In particular this shows that $$ \limu \1U\xi $$ exists for every $\xi $ in $M_x$, provided $M_x$ is equipped with the discrete topology. Equivalently, this notation may be though of as referring to the eventual value of an eventually constant net.

It is then clear that $$ \pi _x:\xi \in M_x \mapsto \limu \1U\xi \in M_x $$ is a well defined $\fld $-linear map, clearly satisfying (iii), (iv) and (v) (note that the first clause of (v) follows from the fact that $n+BJ_x=0$, when $x\notin \src (n)$). In particular (v) implies that $\pi _x$ is idempotent.

To see that $\pi _x$ is $A$-linear, it suffices to observe that, for every $a$ in $A$, and every $\xi $ in $M_x$, $$ \pi _x(a\xi ) = \limu \1Ua\xi = a\limu \1U\xi = a\pi _x(\xi ). $$

Let us now check that $$ \pi _x(M_x) = \mu _x(\B xx). \reqno{(6.10.1)} $$ For this, observe that (v) implies that the range of $\pi _x$ is spanned by the set $$ \{n+BJ_x: n \in \N xx\}. \reqno{(6.10.2)} $$ On the other hand, we have by (5.28.iii) that $\B xx$ is spanned by $$ \{\p xx(n)=n+\H xx: n\in \N xx\}, $$ so $\mu _x(\B xx)$ is spanned precisely by (6.10.2), and therefore (6.10.1) is proved.

In order to verify (i), and in view of (6.10.1), it is then enough to prove that $$ \pi _x(M_x) = W := \bigcap _{U\in \Nbd } \1 UM_x. $$ With this goal in mind, note that (iii) immediately implies that the range of $\pi _x$ is contained in $W$. Conversely, given $\xi $ in $W$, we have that $\xi =\1U\xi $, for every $U$ in $\Nbd $, so $$ \pi _x(\xi ) = \limu \1U\xi = \xi , $$ so $\xi $ lies in the range of $\pi _x$.

It now remains to prove uniqueness of $\pi _x$, so we suppose that $\pi _x'$ is another $A$-linear projection from $M_x$ onto $\mu _x(\B xx)$, and it suffices to prove that $$ \pi _x(n+BJx) = \pi _x'(n+BJx), \for n\in N_x. \reqno{(6.10.3)} $$ Given $n\in N_x$, suppose first that $\beta _n(x)=x$. Then $n\in \C xx$, so $n+BJ_x$ belongs to $\mu _x(\B xx)$, and (6.10.3) clearly holds.

If, on the other hand, $\beta _n(x)\neq x$, we know from (v) that $n+BJ_x$ lies in the kernel of $\pi _x$, so it suffices to show that $$ \xi := \pi _x'(n+BJ_x)=0. $$ For this we choose $U$ in $\Nbd $ satisfying (iii) relative to $\xi $, and we pick $V\in \Nbd $, with $$ \beta _n(x)\notin V\subseteq U. $$ Since $\xi $ lies in the range of $\pi _x'$, which coincides with the range of $\pi _x$ by hypothesis, we have that $$ \xi = \pi _x(\xi ) = \1V\xi = \1V\pi _x'(n+BJ_x) = $$ $$ = \pi _x'(\1V n+BJ_x) \explica {(6.9)}{=} \pi _x'\big (\ev {\1V}{\beta _n(x)}n+BJ_x\big ) = 0. \closeProof $$ \endProof

The following is another important property of $\pi _x$ to be used later.

\state Lemma Let $y, z\in X$, let $n\in \N yx$, and let $m\in \N zx$. If $y\neq z$, then $$ \pi _x(n^*m + BJ_x) = 0. $$

\Proof Since $n^*m$ is a normalizer, let us discuss whether or not $x\in \src (n^*m)$. If not, then (6.7.i) tells us that actually $n^*m + BJ_x = 0$, and we are done. On the other hand, if $x\in \src (n^*m)$, then $$ \beta _{n^*m}(x) = \beta _{n^*}(\beta _m(x)) = \beta _{n^*}(z) \neq \beta _{n^*}(y) = x, $$ so the result follows from (6.10.v). \endProof

\state Proposition For each $x$ in $X$, and for every $y$ in $\Orb x$, define $$ \M yx = \span \{n + BJ_x : n\in \N yx\}. $$ We then have that \item $a\xi =\ev ay\xi $, for all $a\in A$, and all $\xi \in \M yx$, \item $M_x = \bigoplus _{y\in \Orb x} \M yx$, \item if $n$ is any normalizer in $\N yx$, then $\M yx= (n+ BJ_x)\B xx$.

\Proof The first point follows immediately from (6.9). Regarding (ii), we begin by noticing that $M_x$ is spanned by the set $ \{n + BJ_x : n\in N_x\}, $ thanks to (6.8), and since $$ N_x = \bigcup _{y\in \Orb x}\N yx, $$ it is clear that $$ M_x = \sum _{y\in \Orb x} \M yx. $$ It therefore remains to prove that the $\M yx$ are independent. So we choose pairwise distinct elements $$ y_1,y_1,\ldots ,y_k\in \Orb x, $$ and, for each $i$, we choose $\xi _i$ in $\M {y_i}x$, such that $$ \sum _{i=1}^k \xi _i=0. $$ Our task is then to prove that $\xi _i=0$, for all $i$. For this we fix any $i\leq k$, and pick $a$ in $A$, such that $\ev a{y_i}=1$, while $\ev a{y_j}=0$, for all $j\neq i$. We then have by (i) that $$ 0 = \sum _{j=1}^k a\xi _j = \sum _{j=1}^k \ev a{y_j}\xi _j = \xi _i, $$ as desired. In order to prove (iii), choose any $\xi $ in $\M yx$, so by definition we may write $$ \xi = \sum _{i=1}^k \xi _i, $$ where each $ \xi _i = n_i+BJ_x, $ with $n_i\in \N yx$. Next observe that for every $i$, $$ nn^*\xi _i \explica {(i)}{=} \evy {nn^*}\xi _i = \xi _i, $$ because $nn^*$ is idempotent and $\evy {nn^*}\neq 0$, so $\evy {nn^*}=1$. Noticing that $$ c_i:=n^*n_i \in \N xy \N yx\subseteq \N xx \explica {(5.28.i)}{\subseteq} \C xx, $$ and putting $b_i:= c_i+\H xx \in \B xx$, we then have that $$ \xi = \sum _{i=1}^k \xi _i = \sum _{i=1}^k nn^*\xi _i = \sum _{i=1}^k nn^*(n_i+BJ_x) = \sum _{i=1}^k (nc_i+BJ_x) = $$ $$ = \sum _{i=1}^k (n+BJ_x)(c_i+\H xx)= \sum _{i=1}^k (n+BJ_x)b_i= (n+BJ_x)\sum _{i=1}^k b_i, $$ which belongs to $(n+BJ_x)\B xx$, as desired. \endProof

\state Corollary Fixing $x$ in $X$, choose $n^y$ in $\N yx$, for every $y$ in $\Orb x$. Then the imprimitivity bimodule $M_x$ is free as a right $\B xx$-module, with basis $$ \{n_y+BJ_x\}_{y\in \Orb x}. $$

\Proof In view of (6.12), it suffices to show that, given $y$ in $\Orb x$, and given $n\in \N yx$, the map $$ h\in \B xx \mapsto (n+BJ_x)h\in M_x $$ is injective. For this, suppose that $h$ lies in the kernel of this map, and write $h = c+\H xx$, for some $c$ in $\C xx$. Then $$ 0= (n+BJ_x)h = (n+BJ_x)(c+\H xx) = nc+BJ_x, $$ meaning that $nc\in BJ_x$, so also $n^*nc\in BJ_x$. Observing that $n^*n\in \N xx\subseteq \C xx$, we obtain $$ n^*nc\in \C xx\cap BJ_x \explica {(5.15.iii)}{=} J_xBJ_x = \H xx, $$ from where we conclude that $n^*nc+\H xx = 0$. Recalling that $n^*n + \H xx$ is the unit of $\B xx$, by (5.30.iii), we get $$ h = c+\H xx = (n^*n+\H xx) (c+\H xx) = n^*nc+\H xx = 0, $$ concluding the proof of the injectivity of our map. \endProof

\section Restriction

We again place ourselves in the context outlined in (5.6), and we shall also fix a left $B$-module $V$, which we will assume to be \emph {unital}, in the sense that $BV=V$.

Our goal here is to look for subspaces of $V$ which admit the structure of a $\B xx$-module, and with which we will later attempt to reconstruct $V$ via the induction process. The reader should again keep in mind that $A$ is being identified with $L_c(X)$, as described in the paragraph following (5.6). The following is related to (6.10.i).

\state Proposition Given $x$ in $X$, let $$ V_x = \bigcap _{U\in \Nbd } \1 UV. $$ Then \item $ V_x = \big \{v\in V : J_xv = \{0\} \big \}. $ In particular $J_x V_x=\{0\}$, \item $\H xx V_x=\{0\}$, \item $\C xx V_x\subseteq V_x$.

\Proof Regarding (i), and given $a$ in $J_x$, recall that $a$ is a locally constant function on $X$ vanishing on $x$, so there exists $U\in \Nbd $ such that $a$ vanishes identically on $U$, and hence $a\1U=0$. So, for every $v$ in $V_x$, one has that $$ av = a\1Uv =0. $$ This shows the inclusion ``$\subseteq $'' in the statement. In order to prove the reverse inclusion, suppose that $v$ is an element of $V$ such that $J_xv=\{0\}$. Since $V$ is unital, we may write $v=\sum _{i=1}^k b_iv_i$, with $b_i\in B$, and $v_i\in V$. Choosing a dedicated unit $e$ for the set $\{b_1,b_2,\ldots ,b_k\}$, belonging to $A$, we can easily prove that $ev=v$. Therefore, if $U$ is any member of $\Nbd $, we have that $\1Ue-e$ lies in $J_x$, so $$ 0 = (\1Ue-e)v = \1U v - v, $$ and we see that $v=\1Uv$. Since $U$ is arbitrary, we then conclude that $v$ belongs to $V_x$.

Clearly (i) implies (ii).

In order to verify (iii), given $c$ in $\C xx$ and $v$ in $V_x$, we have $$ J_xcv \subseteq BJ_xv = \{0\}, $$ so $cv\in V_x$, by (i). \endProof

We may then view $V_x$ as a left $\C xx$-module, and since $\H xx$ annihilates $V_x$ by (7.1.ii), we may turn $V_x$ into a left module over $\C xx/\H xx=\B xx$, with the module structure given by $$ (c + \H xx) v = cv, \reqno{(7.2)} $$ for all $c$ in $\C xx$, and all $v$ in $V_x$.

\state Definition \rm The left $\B xx$-module $V_x$ described above will be called the \emph {restriction of\/ $V$ relative to $x$}, and it will be denoted by $\Resx V$.

The reader should notice that $\Resx V$ is a subset of the left $B$-module $V$, but it is not a $B$-submodule. On the other hand, $V_x$ is a left $\B xx$-module, but $V$ is not a $\B xx$-module.

\state Proposition $V_x$ is a unital left $\B xx$-module.

\Proof Given $U$ in $\Nbd $, we have by (5.30.ii) that $\1U+\H xx$ is the unit of $\B xx$, and clearly $\1U v = v$, for all $v$ in $\Resx V$. \endProof

A first and very important example is the case of $M_x$, itself.

\state Proposition Given $x$ in $X$, and regarding $M_x$ as a left $B$-module, we have that $\Resx {M_x}$ is isomorphic to $\B xx$, as left $\B xx$-modules.

\Proof From (6.10.i), it follows that $\Resx {M_x}$ coincides with $\mu _x(\B xx)$, as sets, so is suffices to prove that $\mu _x$ is $\B xx$-linear\fn {Observe that $M_x$ is \emph {not} a left $\B xx$-module in any relevant way, so it makes no sense to ask whether or not $\mu _x$ is $\B xx$-linear, unless we restrict the codomain of $\mu _x$ to $\Resx {M_x}$, which was indeed given the structure of a left $\B xx$-module in (7.2).}, as a map from $\B xx$ onto $\Resx {M_x}$. Technically speaking, our task is to prove that $$ \mu _x\big ((c+\H xx)(d+\H xx)\big ) = c \mu _x(d+\H xx), \for c,d\in \C xx. $$ Developing from the left-hand-side, we have that $$ \mu _x\big ((c+\H xx)(d+\H xx)\big ) = \mu _x\big (cd+\H xx) \explica {(6.4)}{=} cd+BJ_x, $$ while $$ c\mu _x(d+\H xx) \explica {(6.4)}{=} c(d+BJ_x) \explica {(7.2)}{=} cd+BJ_x, $$ as desired. \endProof

It is quite possible that the $V_x$ all vanish, even if $V$ is nonzero. However, this is not so in the finite dimensional case.

\state Proposition Let $V$ be a nonzero, unital, left $B$-module. If\/ $V$ is finite dimensional as a $\fld $-vector space, then there exists $x$ in $X$ such that $V_x\neq \{0\}$.

\Proof Suppose by contradiction that $V_x=\{0\}$, for every $x$ in $X$. We then claim that, given $x$, there exists some $U\in \Nbd $, such that $\1UV=\{0\}$. To prove this claim, pick $U_0$ in $\Nbd $ such that $$ \dim (\1{U_0}V) = \min \big \{\dim (\1{U}V) : U\in \Nbd \big \}. $$ Since, for every $U$ in $\Nbd $, one has that $$ \1{U_0\cap U}V= \1{U_0}\1{U}V\subseteq \1{U_0}V, $$ we deduce that $\1{U_0\cap U}V=\1{U_0}V$, by minimality. Consequently $$ \1{U_0}V = \1{U_0\cap U}V= \1{U}\1{U_0}V\subseteq \1{U}V, $$ and since $U$ is arbitrary, we get $$ \1{U_0}V \subseteq \medcap _{U\in \Nbd }\1{U}V = V_x = \{0\}, $$ proving the claim.

We next claim that every vector $v$ in $V$ has \emph {compact support}, in the sense that there exists a compact open subset $U\subseteq X$, such that $\1Uv=v$. To see this, observe that $V$ is unital, so we may write $$ v=\sum _{i=1}^kb_iv_i, $$ with the $v_i$ in $V$, and the $b_i$ in $B$. Choosing a dedicated unit $u$ for the set $\{b_1,\ldots ,b_k\}$, belonging to $A$, it is easy to see that $uv=v$. Now, letting $U$ be the support of $u$, we have that $U$ is a compact open set, and clearly $\1Uu=u$. Therefore $$ \1Uv = \1Uuv = uv = v, $$ proving the claim.

Let us next fix a nonzero vector $v$ in $V$, and let $U$ be such that $\1Uv = v$, as above. We then use the first part of the proof to obtain a finite cover $$ U = \bigcup _{i=1}^k U_i, $$ such that $\1{U_i}V=\{0\}$, for every $i$. By using the inclusion-exclusion principle we may assume that the $U_i$ are pairwise disjoint, so $$ \1U = \sum _{i=1}^k \1{U_i}, $$ from where we see that $\1UV=\{0\}$, and hence $$ v = \1Uv = 0, $$ a contradiction. This concludes the proof. \endProof

\section Induction

Our goal here will be to study the induction process, already introduced in (6.2), in greater detail. So it is perhaps worth recalling from (6.2) that, if $V$ is any left $\B xx$-module, the left $B$-module induced from $V$ is $$ \Indx V = M_x\tprodx V. $$

\fix Throughout this section, besides assuming the conditions of (5.6), we will fix a point $x$ in $X$, as well as a unital left $\B xx$-module $V$.

\bigskip Recalling that the range of the map $\pi _x$ introduced in (6.10.ii) coincides with the range of the map $\mu _x$ introduced in (6.4), and recalling also that $\mu _x$ is injective, we may define a map $\nu _x : M_x \to \B xx$, by $\nu _x = \mu _x^{-1}\circ \pi _x$,

\begingroup \noindent \hfill \beginpicture \setcoordinatesystem units <0.025truecm, -0.02truecm> \setplotarea x from 70 to 150, y from -30 to 150 \put {\null } at 70 -30 \put {\null } at 70 150 \put {\null } at 150 -30 \put {\null } at 150 150 \put {$\B xx$} at 70 100 \put {$\mu _x(\B xx) = \pi _x(M_x)$} at 200 100 \arrow <0.11cm> [0.5, 1.8] from 93.4 100 to 128.5 100 \put {$\mu _x$} at 110.95 115 \put {$M_x$} at 238 0 \arrow <0.11cm> [0.5, 1.8] from 238 18 to 238 82 \put {$\pi _x$} at 248 50 \arrow <0.11cm> [0.5, 1.8] from 212.36 9.8 to 85.64 80.2 \put {$\nu _x$} at 141.715 31.888 \endpicture \hfill \null \endgroup

\state Proposition The map $\nu _x$ defined above satisfies the following: \item $\mu _x\circ \nu _x = \pi _x$, and $\nu _x\circ \mu _x = id_{\B xx}$, \item $\nu _x$ is right $\B xx$-linear, \item for every normalizer $n$ one has that $$ \nu _x(n+BJ_x) = \clauses { \cl {n+\H xx}{\text {if } n\in \N xx,} \cl {0}{\text {otherwise.}} } $$

\Proof The first point is clear. As for (ii), since both $\pi _x$ and $\mu _x$ are right $\B xx$-linear, so is $\nu _x$. Regarding (iii) in the case that $n\in \N xx$, notice that $\N xx\subseteq \C xx$, so $n+\H xx$ is a legitimate element of $\B xx$, making the expression there syntactically correct.

Notice that, if $x\notin \src (n)$, then (6.7.i) tells us that $n+BJ_x=0$, so (iii) holds. Assuming instead that $x\in \src (n)$, suppose first that $\beta _n(x)\neq x$. Then (6.10.v) implies that $\pi _x(n+BJ_x) = 0$, so again (8.1.iii) checks. Finally, supposing that $\beta _n(x)=x$, that is, supposing that $n\in \N xx$, we have that $$ \nu _x(n+BJ_x) = \mu _x^{-1}\big (\pi _x(n+BJ_x)\big ) \explica {(6.10.v)}{=} \mu _x^{-1}(n+BJ_x) \explica {(6.4)}{=} n+\H xx. $$ This concludes the proof. \endProof

Having already fixed a point $x$ in $X$, let us also fix a normalizer $n_y$ in $\N yx$, for each $y$ in $\Orb x$. Observing that, for every $U$ in $\Nbd $, we have that $\1U\in \N xx$, we shall insist in choosing $n_x=\1U$, where $U$ is fixed beforehand in $\Nbd $. We will then set $$ \zeta _y = n_y+BJ_x \in B/BJ_x = M_x, \for y\in \Orb x. \reqno{(8.2)} $$

\state Proposition Given a unital left $\B xx$-module $V$, consider the map $$ \varphi _x:= \nu _x\otimes id_V : M_x\tprodx V \to \B xx\tprodx V = V, $$ where the identification ``$\B xx\tprodx V = V$'', above, is the usual one, given that $V$ is a unital module. Then: \item For every $\xi \in M_x\tprodx V$, one has that $$ \xi =\sum _{y\in \Orb x}\zeta _y\otimes \varphi _x(n_y^*\xi ), $$ \item For every submodule $W\subseteq V$, one has that $$ M_x \otimes W = \bigoplus _{y\in \Orb x} \zeta _y\otimes W, $$ where $M_x \otimes W$, above, is to be interpreted in light of (6.3). \item For every $y$ in $\Orb x$, the map $$ v\in V\mapsto \zeta _y\otimes v\subseteq M_x\tprodx V $$ is injective.

\Proof We begin by observing that, for all $n\in N_B(A) $, and all $v\in V$, we have by (8.1.iii) that $$ \varphi _x\big ((n+BJ_x)\otimes v\big ) = \clauses { \cl {(n+\H xx)v}{\text {if } n\in \N xx,} \cl {0}{\text {otherwise.}} } \reqno{(8.3.1)} $$

This said, from (6.13) we know that $M_x$ is free with basis $\{\zeta _y\}_{y\in \Orb x}$, so it follows that $$ M_x\tprodx V = \bigoplus _{y\in \Orb x}\zeta _y\otimes V. \reqno{(8.3.2)} $$ Given $\xi \in M_x\tprodx V$, one may then write $$ \xi =\sum _{y\in \Orb x}\zeta _y\otimes v_y, $$ with the $v_y\in V$, and such that $v_y=0$ for all but finitely many $y$ in $\Orb x$. We then claim that $$ \varphi _x(n_z^*\xi ) = v_z, \for z\in \Orb x. \reqno{(8.3.3)} $$ In order to prove it, we fix $z$, and compute $$ \varphi _x(n_z^*\xi ) = \sum _{y\in \Orb x }\varphi _x(n_z^*\zeta _y\otimes v_y) = \sum _{y\in \Orb x }\varphi _x\big ((n_z^*n_y + BJ_x)\otimes v_y\big ). $$ Regarding the summand corresponding to $y=z$, notice that $n_z^*n_z$ lies in $\N xx$, so (8.3.1) implies that $$ \varphi _x\big ((n_z^*n_z + BJ_x)\otimes v_z\big ) = (n_z^*n_z + \H xx)v_z \explica {(5.30.ii)}{=} v_z. $$ With respect to the other summands, namely when $y\neq z$, we have by (5.12.v) that $n_z^*n_y$ is not in $\N xx$, so (8.3.1) says that all such summands vanish so (8.3.3) is proved, and hence so is (i).

Point (ii) follows immediately from (8.3.2), and if $v$ is such that $$ \xi := \zeta _y\otimes v=0, $$ then $$ v \explica {(8.3.3)}{=} \varphi _x(n_y^*\xi )=0, $$ whence (iii). \endProof

The following is one of our main results:

\state Theorem Assuming the hypotheses of (5.6), pick any $x$ in $X$ and let $V$ be a unital left $\B xx$-module. Then $V$ is naturally isomorphic to $\Resx {\Indx V}$. More precisely, letting $\zeta _x=\1{U_0}+BJ_x\in M_x$, for a fixed\fn {Observe that, if $U$ and $V$ lie in $\Nbd $, then $\1U-\1V\in J_x\subseteq BJ_x$, so $\1U+BJ_x=\1V+BJ_x$. In other words, the class of $\1U$ modulo $BJ_x$ does not depend on the choice of $U$.} $U_0$ in $\Nbd $, the map $$ v\in V\mapsto \zeta _x\otimes v\in M_x\tprodx V = \Indx V $$ is a left $\B xx$-module isomorphism onto $\Resx {\Indx V}$.

\Proof For every $y$ in $\Orb x$, let $n_y$ and $\zeta _y$ be as in (8.2), extending the choice already made for $\zeta _x$ in the statement. Given any $\xi $ in $M_x\tprodx V$, and using (8.3.ii), we write $$ \xi =\sum _{y\in \Orb x}\zeta _y\otimes v_y, $$ where the $v_y\in V$ vanish for all but finitely many $y$ in $\Orb x$. Using (6.12.i), for every $U$ in $\Nbd $, and every $y$ in $\Orb x$, we have that $$ \1U \zeta _y\otimes v_y = \evy {\1U}\zeta _y\otimes v_y, $$ so $$ \1U\xi =\sum _{y\in \Orb x \cap U}\zeta _y\otimes v_y, $$ and it is then clear that $\xi \in \Resx {\Indx V}$ if and only if $\xi =\zeta _x\otimes v_x$. This shows that the range of the map in the statement is precisely $\Resx {\Indx V}$. By (8.3.iii) this map is also injective so it remains to show $\B xx$-linearity.

For this, let $c\in \C xx$, and $v\in V$. Recalling that $n_x=1_U$, and hence that $n_x+\H xx$ is the unit of $\B xx$, we have that $$ (c+\H xx)\zeta _x\otimes v = (n_x+\H xx)(c+\H xx)(n_x+BJ_x)\otimes v = (n_xcn_x+BJ_x)\otimes v = $$ $$ = (n_x+BJ_x)(cn_x+\H xx)\otimes v = (n_x+BJ_x)\otimes (cn_x+\H xx)v = $$ $$ = (n_x+BJ_x)\otimes (c+\H xx)(n_x+\H xx)v = \zeta _x\otimes (c+\H xx)v, $$ concluding the proof. \endProof

The following is another important consequence of (8.3).

\state Corollary Under (5.6), one has that $\FunctorInd x $ is an exact functor.

\Proof Since induction consists of tensoring with the free module $M_x$, the conclusion follows. \endProof

\section Irreducibility of induced modules

Here we want to to look at how does the property of being irreducible for a module affects its induced counterpart. As always, we keep working under (5.6).

\def \Subm {\mathcal {S}}

Given any algebra $C$, and given any left $C$-module $V$, let us denote by $\Subm _C(V)$ the family of all of its submodules. In symbols: $$ \Subm _C(V) = \{W\subseteq V: W \text { is a $C$-submodule of } V\}. $$ Regarding the order relation given by inclusion, it is easy to see that $\Subm _C(V)$ is a lattice, where the meet operation is the intersection of modules, and the join operation is the sum of modules. Clearly $\Subm _C(V)$ has a biggest element, namely $V$, and a smallest element, $\{0\}$.

If $V$ is now a left $\B xx$-module, and if $W$ is in $\Subm _{\B xx}(V)$, we may consider the map $$ \Indx W = M_x\tprodx W \labelarrow {id_{M_x}\otimes \iota } M_x\tprodx V = \Indx V, \reqno{(9.1)} $$ where $\iota $ is the inclusion map from $W$ into $V$, which is then an injective $B$-module homomorphism by exactness. Its range, which we have agreed in (6.3) to denote by $M_x\otimes W$, is therefore a $B$-submodule of $\Indx V$, which may be identified with $\Indx W$ through the above map. We then obtain a map $$ W\in \Subm _{\B xx}(V)\ \longmapsto \ \Indx W\in \Subm _{B}\big (\Indx V\big ). \reqno{(9.2)} $$

\state Proposition The above map is a lattice isomorphism. In particular, the submodules of the module induced by $V$ are precisely the modules induced by the submodules of $V$.

\Proof Let $n_y$ and $\zeta _y$ be as in (8.2), where we again insist in choosing $n_x=\1U$, where $U$ is a compact open neighborhood of $x$.

If $W$ is a submodule of $V$, and using (8.3.ii), it is easy to see that $$ (\zeta _x\otimes V)\cap \big (M_x\otimes W\big ) = \zeta _x\otimes W. $$ Since the map $$ v\in V\mapsto \zeta _x\otimes v\in M_x\tprodx V $$ is injective by (8.3.iii), we conclude that $$ W = \{v\in V: \zeta _x\otimes v\in M_x\otimes W\}. $$ The fact that we are able to recover $W$ from $M_x\otimes W$, as above, implies that our correspondence is injective.

In order to show our correspondence to be onto, pick $Z$ in $\Subm _{B}\big (\Indx V\big )$, and let us prove that $$ W:=\{v\in V : \zeta _x\otimes v \in Z\} $$ is a $\B xx$-submodule of $V$. Indeed, given $b$ in $\B xx$, and $w$ in $W$, write $b = c+\H xx$, for some $c$ in $\C xx$. Then $$ \zeta _x\otimes bw = (n_x +BJ_x)\otimes (c+\H xx)w = $$ $$ =(n_x +BJ_x)(c+\H xx)\otimes w = (n_xc +BJ_x)\otimes w = \cdots $$ Recalling that $n_x=1_U$, we have by (5.30.ii) that $n_x+\H xx$ is the unit of $\B xx$, so the above equals $$ \cdots = (n_xc +BJ_x)\otimes (n_x+\H xx)w = (n_xc +BJ_x)(n_x+\H xx)\otimes w = $$ $$ = (n_xcn_x +BJ_x)\otimes w = n_xc(n_x +BJ_x)\otimes w = n_xc\, \zeta _x\otimes w \in Z. $$ This shows that $bw\in W$, as required. We will then show that $Z = M_x\otimes W$. Focusing on proving the inclusion ``$\supseteq $'', let $b\in B$ and $w\in W$. Then, once more using that $n_x+\H xx$ is the unit of $\B xx$, we have $$ (b+BJ_x)\otimes w = (b+BJ_x)\otimes (n_x+\H xx)w = (b+BJ_x)(n_x+\H xx)\otimes w = $$ $$ = (bn_x+BJ_x)\otimes w = b\zeta _x\otimes w \in Z. $$ This proves that $M_x\otimes W\subseteq Z$. To prove the reverse inclusion, pick any $m\in Z$, and write $$ m = \sum _{y\in \Orb x} \zeta _y\otimes v_y, \reqno{(9.3.1)} $$ where the $v_y\in V$ vanish for all but finitely many $y$ in $\Orb x$, according to (8.3.ii). We will first prove that $m\in M_x\otimes W$ under the hypothesis that there is a single nonzero $v_y$, above, that is, $m = \zeta _y\otimes v$, for some $y$ in $\Orb x$, and $v\in V$. For this notice that $$ Z\ni n_y^*m = n_y^*\zeta _y\otimes v = n_y^*(n_y+BJ_x)\otimes v = (n_y^*n_y+BJ_x)\otimes v = \cdots $$ Observing that both $n_y^*n_y$ and $n_x=\1U$ lie in $A$, and that $n_y^*n_y-n_x$ belongs to $J_x$, we see that the above equals $$ \cdots = (n_x+BJ_x)\otimes v = \zeta _x\otimes v. $$ This implies that $v\in W$, whence $$ m = \zeta _y\otimes v \in M_x\otimes W, $$ proving the claim in the special case of a single nonzero $v_y$. In the general case, fix any $y$ in $\Orb x$, and choose $a$ in $A$, such that $\ev ay=1$, and $\ev az=0$, for all $z$ in $F\setminus \{y\}$, where $$ F = \{y\in \Orb x: v_y\neq 0\}, $$ which is a finite set. Then $$ Z \ni am = \sum _{z\in F} a\zeta _z\otimes v_z \explica {(6.12.i)}{=} \sum _{z\in F} \ev az\zeta _z\otimes v_z = \zeta _y\otimes v_y. $$ By the first case treated above, we have that $v_y$ lies in $W$, so $m$ belongs to $M_x \otimes W$, as needed.

Finally, it is now easy to prove that, if $W_1$ and $W_2$ are submodules of $V$, then $M_x\otimes W_1 \subseteq M_x\otimes W_2$ if and only if $W_1 \subseteq W_2$. In other words, our correspondence preserves order relations so it is a lattice isomorphism. \endProof

As a consequence we have:

\state Proposition Given $x$ in $X$, let $V$ be a unital, nontrivial,\fn {That is $V\neq \{0\}$.} left $\B xx$-module. Then: \item $\Indx V$ is irreducible if and only if $V$ is irreducible. \item $\Indx V$ is indecomposable if and only if $V$ is indecomposable.

\Proof If $C$ is any algebra and $V$ is a nontrivial $C$-module, then $V$ is reducible if and only if $\Subm _C(V)$ has at least three elements. On the other hand, $V$ is decomposable if and only if there are nonzero elements $x, y\in \Subm _C(V)$, such that $x\wedge y$ is the smallest element, and $x\vee y$ is the biggest element. In other words, deciding whether or not $V$ is irreducible or indecomposable may be done by looking only at the lattice structure of $\Subm _C(V)$. The conclusion then follows immediately from (9.3). \endProof

\section Inducing restricted modules

Having already understood the restriction of an induced module in (8.4), we now discuss the opposite construction.

\state Theorem Assuming (5.6), let $V$ be a left $B$-module, and let $x\in X$. Then there exists a natural injective $B$-linear map $$ \rho :\Indx {\Resx V} \to V, $$ such that $$ \rho \big ((b+BJ_x)\otimes v\big ) = bv, \for b\in B, \for v\in \Resx V. \reqno{(10.1.1)} $$

\Proof To shorten our notation, we write $V_x$ for $\Resx V$. Given $v\in V_x$, recall that $J_xv=\{0\}$, by (7.1.i), so the map $$ b\in B\mapsto bv\in V $$ vanishes on $BJ_x$ and hence factors through $M_x$, so the map $$ \rho _0: (b+BJ_x,v) \in M_x\times V_x\mapsto bv\in V $$ is well defined. It is also clearly $\fld $-bilinear and $\B xx$-balanced, so it defines a map $$ \rho : \Indx {V_x} = M_x\tprodx V_x \to V, $$ such that $$ \rho (\xi \otimes v) = \rho _0(\xi , v), \for \xi \in M_x, \for v\in V_x, $$ hence satisfying (10.1.1), which in turn implies that $\rho $ is left $B$-linear.

It therefore remains to prove that $\rho $ is one-to-one. Picking any $t$ in $M_x\tprodx V_x$, such that $\rho (t)=0$, write $$ t = \sum _{y\in \Orb x} \zeta _y\otimes v_y, \reqno{(10.1.2)} $$ were the $\zeta _y$ are as in (8.2), and the $v_y\in V_x$ vanish for all but finitely many $y$ in $\Orb x$. Fixing $$ y \in F := \{y\in \Orb x: v_y\neq 0\}, $$ choose $a$ in $A$, such that $\ev ay=1$, and $\ev az=0$, for all $z$ in $F\setminus \{y\}$. We then obtain $$ 0 = a\rho (t) = \rho (at) = \rho \Big (\sum _{z\in F} a\zeta _z\otimes v_z\Big ) \explica {(6.12.i)}{=} $$ $$ = \rho \Big (\sum _{z\in F} \ev az\zeta _z\otimes v_z\Big ) = \rho (\zeta _y\otimes v_y) = \rho \big ((n_y+BJ_x)\otimes v_y\big ) = n_yv_y. $$

Observing that $\evx {n_y^*n_y}=1$, and that $n_y^*n_y$ is locally constant, we may find $U$ in $\Nbd $ such that $\1Un_y^*n_y=\1U$. Therefore, recalling that $v_y$ lies in $V_x$, we have that $$ v_y = \1Uv_y = \1Un_y^*n_yv_y = 0, $$ proving that $t=0$, as required. \endProof

We are now in a position to present another main result.

\state Corollary Assuming the conditions of (5.6), let $V$ be an irreducible left $B$-module such that $\Resx V$ is nonzero for some $x$ in $X$, e.g.~when $V$ is finite dimensional as a $\fld $-vector space (c.f. (7.6)). Then $\Resx V$ is an irreducible left $\B xx$-module, and $V$ is isomorphic to $ \Indx {\Resx V}. $

\Proof The range of the map $\rho $ of (10.1) is a nonzero submodule of $V$, and hence equal to $V$, by irreducibility, so $\rho $ is an isomorphism of $B$-modules. That $\Resx V$ is irreducible follows from (9.4.i), given that it induces an irreducible $B$-module, namely $V$. \endProof

\section The annihilator of induced modules

With this section we start giving more emphasis on the annihilator of a module, rather than the module itself. Since every ideal\fn {In this work the term \emph {ideal} will always mean \emph {2-sided ideal}, unless stated otherwise.} in an s-unital algebra is the annihilator of a module, our results will have consequences for the study of ideals.

Recall that if $V$ is a module over an algebra $C$, then the \emph {annihilator} of $V$ in $C$ is defined by $$ \Ann C{V} = \big \{c\in C: cv=0, \ \forall v\in V\big \}. $$

In this section we would like to discuss the annihilator of induced modules, so we again put ourselves in the context of (5.6).

An important tool will be the isotropy projection $\E xx$, introduced in (5.23) (see also (5.29)). Incidentally we should mention that $\E xx$ is closely related to the map $\nu _x$, introduced in (8.1), as the following makes clear.

\state Proposition Given $x$ in $X$, one has that $$ \E xx (b) = \nu _x(b+BJ_x),\for b\in B. $$

\Proof Since $B$ is spanned by $N_B(A)$, it is enough to prove that $\E xx (n) = \nu _x(n+BJ_x)$, for every normalizer $n$. Given such an $n$, let us first suppose that $n$ is not in $\N xx$. Then $$ \E xx(n)=\Psi ^{-1}(\q xx(n)) \explica {(5.28.ii)}{=} 0, $$ while $\nu _x(n)$ also vanishes by (8.1.iii). In case $n\in \N xx$ we have that $$ \E xx(n)=\Psi ^{-1}(\q xx(n)) \explica {(5.19)}{=} \p xx(n) = n + \H xx \explica {(8.1.iii)}{=} \nu _x(n+BJ_x), $$ and the proof is over. \endProof

The main result of this section is now in order:

\state Proposition Under (5.6), pick any $x$ in $X$, and let $V$ be a unital left $\B xx$-module. Then the annihilator of $\Indx V$ is given by $$ \Ann B{\Indx V} = \big \{b\in B: \E xx(g bh )\in \Ann {\B xx}{V},\ \forall g , h \in B\big \}. $$

\Proof Given $h \in B$, and $v\in V$, let $$ \xi =(h +BJ_x)\otimes v\in M_x\tprodx V. $$

Using (8.3.i) it is easy to see that $\xi =0$ if and only if $\varphi _x(g \xi )=0$, for all $g \in B$. Incidentally, we have $$ \varphi _x(g \xi ) = (\nu _x\otimes id_V)\big ((g h +BJ_x)\otimes v\big ) = \nu _x(g h +BJ_x)\otimes v = $$ $$ = 1\otimes \nu _x(g h +BJ_x)v = 1\otimes \E xx(g h )v. $$ So, $$ (h +BJ_x)\otimes v=0 \iff \E xx(g h )v=0, \ \forall g \in B. $$

For $b\in B$, it then follows that $b$ lies in the annihilator of $\Indx V$ if and only if $$ b(h +BJ_x)\otimes v = 0, \ \forall h \in B,\ \forall v\in V $$ $$ \iff \E xx(g bh )v = 0, \ \forall g , h \in B,\ \forall v\in V $$ $$ \iff \E xx(g bh )\in \Ann {\B xx}{V}, \ \forall g , h \in B. \closeProof $$ \endProof

Observe that the above description of the annihilator of $\Indx V$ depends only on the annihilator of $V$ rather than on $V$ itself. This in turn suggests the following:

\state Definition \rm Given $x$ in $X$, and given an ideal $I\trianglelefteq \B xx$, the \emph {ideal induced} by $I$ is defined to be the ideal of $B$ given by $$ \Ind x I= \big \{b\in B: \E xx(g bh )\in I,\ \forall g , h \in B\big \}. $$

The reader is invited to check that $\Ind x I$ is indeed an ideal in $B$.

Using this terminology, the conclusion of (11.2) can be concisely restated as $$ \Ann B{\Indx V} = \Indx {\Ann {\B xx}{V}}. \reqno{(11.4)} $$

Recall that an ideal is said to be \emph {primitive} if it coincides with the annihilator of an irreducible module.

\state Proposition Under (5.6), pick any $x$ in $X$, and let $I$ be a primitive ideal of $\B xx$. Then $\Indx I$ is a primitive ideal of $B$.

\Proof Write $I=\Ann {\B xx}V$, for some irreducible $\B xx$-module $V$. Then, as seen above, $$ \Indx I = \Indx {\Ann {\B xx}V}= \Ann B{\Indx V}, $$ and the conclusion follows because $\Indx V$ is irreducible by (9.4.i) and hence its annihilator is a primitive ideal. \endProof

\section Germs and the disintegration of modules

We would now like to describe a construction similar to the restriction of a module defined in Section (7), but which has a better chance of producing nonzero modules.

\fix Again working under (5.6), we shall now fix a unital left $B$-module $V$.

\state Definition \rm Given $x$ in $X$, and given $v$ and $w$ in $V$, let us say that $v$ and $w$ are \emph {$x$-equivalent}, in symbols $v\sim _xw$, if there exists a compact open neighborhood $U$ of $x$, such that $\1Uv=\1Uw$. The corresponding equivalence class of any given $v$ in $V$ will be called the \emph {germ} of $v$ at $x$, and it will be denoted by $[v,x]$.

\state Proposition Given $v$ and $w$ in $V$, one has that $$ v\sim _xw \iff v-w\in J_xV. $$

\Proof Assuming that $v\sim _xw$, choose $U$ in $\Nbd $ such that $\1Uv=\1Uw$. Using that $V$ is unital, and that $A$ is an s-unital subalgebra of $B$, it is easy to see that there exists $a$ in $A$, such that $v=av$, and $w= aw$. Observing furthermore that $a-\1Ua$ lies in $J_x$, we then have that $$ J_xV \ni (a-\1Ua)(v-w) = v-w - \1Uv + \1Uw = v-w. $$

Conversely, if $v-w\in J_xV$, then we may write $$ v-w=\sum _{i=1}^k a_iv_i, $$ where the $a_i\in J_x$, and the $v_i\in V$. Since the $a_i$ are locally constant functions vanishing at $x$, we may find a compact open neighborhood $U$ of $x$, such that $\1Ua_i=0$, for all $i$, and this in turn implies that $\1U(v-w)=0$, as desired. \endProof

In conclusion we see that the space of germs is nothing but the quotient space of $V$ by $J_xV$, while $$ [v,x] = v + J_xV, \for v\in V. $$

\state Definition \rm Given a unital left $B$-module $V$, and given $x$ in $X$, the $\fld $-vector space $$ \Vx := V/J_xV $$ will be called the \emph {space of germs} for $V$ at the point $x$.

The brackets around $x$ in the above notation have no reason other than distinguishing this from the notation $V_x$ introduced in (7.1), while at the same time hinting that germs are involved.

Since $J_xV$ is not a $B$-submodule of $V$, one cannot expect $\Vx $ to have the structure of a $B$-module. However, viewing the disjoint union of the $\Vx $ as a bundle over $X$ (it is actually a sheaf, but this will not be relevant here), we will see that it is acted upon by the isotropy modules $\B xy$ in a way to be made precise in a moment.

For this, let us fix $x,y\in X$, and let $c\in \C yx$. Then, $$ cJ_xV\subseteq J_yBV \subseteq J_yV, $$ so the operation of multiplication by $c$ on $V$ factors through the quotient spaces providing a $\fld $-linear map $$ \lambda _c: v+J_xV \in \Vx \mapsto cv+J_yV \in \V y. $$ In case $c$ lies in $\H yx = J_yBJ_x$, observe that $cv\in J_yV$, for every $v$, meaning that $\lambda _c$ vanishes. So it follows that the map $\lambda :c\mapsto \lambda _c$ factors through $\C yx/\H yx = \B yx$, providing a map $$ \bar \lambda : \B yx \to \rsbox {L}\, (\Vx ,\V y), $$ such that $$ \bar \lambda (c+\H yx)(v+J_x V) =cv+J_y V, $$ for all $c\in \C yx$, and all $v\in V$. Rather than emphasizing the \emph {operator} nature of $\bar \lambda $, we will simply use juxtaposition, thus obtaining a clearly bilinear \emph {multiplication operation} $$ \B yx \times \Vx \to \V y, \reqno{(12.4)} $$ such that $$ (c+\H yx)(v+J_x V) =cv+J_y V, \for c\in \C yx, \for v\in V. $$

\state Proposition The above collection of multiplication operators is associative, in the sense that, given $x,y,z\in X$, and given $g\in \B zy$, $h\in \B yx$, and $u\in \Vx $, one has that $$ (gh)u = g(hu). $$

\Proof Left for the reader. \endProof

In view of \cite {SteinbergDisintegr}, the collection of multiplication operators given by (12.4) could be called the \emph {disintegration} of the $B$-module $V$. However we believe it to be highly unlikely that a complete theory of disintegration be developped in the present level of generality. In particular, the reverse procedure of \emph {integrating} representations, such as the one described in \cite [Section 3]{SteinbergDisintegr}, is probably out of reach for us here. Nevertheless we will be able to use this procedure to our advantage by proving a version of Effros-Hahn conjecture below.

Observe that, in case $x=y$, the operation $$ \B xx \times \Vx \to \Vx , $$ provided by (12.4), turns $\Vx $ into a left $\B xx$-module.

\state Proposition If\/ $V$ is a unital left $B$-module, and if $x\in X$, then $\Vx $ is a unital left $\B xx$-module.

\Proof Let $a\in A$ be such that $\evx a=1$. Fixing $b\in B$, and $v\in V$, it is easy to see that $ab-b\in J_xB$, so $$ abv -bv \in J_x V, $$ whence $$ (a+\H xx)(bv + J_xV) = abv + J_xV = bv + J_xV, $$ proving that $a+\H xx$ acts as a unit on vectors of the form $bv + J_xV$. Since these span $\Vx $, we have proved the statement. \endProof

In the preamble of this section, when we claimed that the present construction has a better chance of producing nonzero modules, this is what we meant:

\state Proposition Under (5.6), let $V$ be a unital left $B$-module. Then, for each nonzero element $v$ in $V$, there exists $x$ in $X$ such that the germ of $v$ at $x$ is nonzero.

\Proof Suppose by way of contradiction that $$ v+J_xV=0, \for x\in X. $$ Then, for every $x$, we may find a compact open neighborhood $U_x$ of $x$, such that $\1{U_x}v=0$. Writing $v=av$, for some $a$ in $A$, we then have that the $U_x$ form an open cover of the compact set $\supp (a)$, so there is a finite subcover, say $$ \supp (a) \subseteq \medcup _{i=1}^k U_{x_i}. $$ Employing the inclusion-exclusion principle one obtains another cover for $\supp (a)$, namely $$ \supp (a) \subseteq \medcup _{i=1}^l W_i, $$ such that the $W_i$ are now pairwise disjoint, and such that $\1{W_i}v$ still vanishes for all $i$. It then follows that $a = \sum _{i=1}^l a\1{W_i}$, whence $$ v = av = \sum _{i=1}^l a\1{W_i}v = 0, $$ a contradiction. \endProof

Fixing $x$ in $X$, and in possession of the module of germs $\Vx $, one can consider the induced module $\Indx {\Vx }$.

\state Proposition Working under (5.6), let $V$ be a unital left $B$-module. Given $x$ in $X$, and given $y$ in $\Orb x$, one has that $$ \Indx {\Vx } \simeq \Ind y {\V y } $$ as left $B$-modules.

\Proof Since $y$ lies in $\Orb x$, we may choose a normalizer $n$ in $\N yx$. We then have that $n^*\in \N xy\subseteq \C xy$, so $$ BJ_xn^* \subseteq BBJ_y \subseteq BJ_y, $$ and hence the map between imprimitivity bimodules $ \varphi _n : M_x \to M_y, $ given by $$ \varphi _n(b+BJ_x) = bn^* + BJ_y, \for b\in B, $$ is well defined. Likewise, $$ nJ_xV\subseteq J_yBV\subseteq J_yV, $$ so the map between germ spaces $ \psi _n : \Vx \to \V y, $ given by $$ \psi _n(v+J_xV) = nv + J_yV, \for v\in B, $$ is well defined. This allows us to define the map $$ \varphi \times \psi :(\xi , u)\in M_x\times \Vx \mapsto \varphi (\xi )\otimes \psi (u)\in M_y\tprod y\V y, $$ which we claim to be $\B xx$ balanced. Indeed, given $c\in \C xx$, $b\in B$, and $v\in V$, we have that $$ (\varphi \times \psi )\big ((b+BJ_x) (c+\H xx), v + J_xV\big ) = (\varphi \times \psi )(bc+BJ_x, v + J_xV) = $$ $$ = (bcn^*+BJ_y)\otimes (nv + J_yV) = \cdots $$ Since $n$ is in $\N yx$, we have that $\evx {n^*n}=1$, whence $b-bn^*n$ lies in $BJ_x$, thanks to (5.26.ii). Consequently $$ (b-bn^*n)cn^* \in BJ_x\C xx\N xy \subseteq BJ_x\C xy \subseteq BBJ_y \subseteq BJ_y, $$ so $bcn^* + BJ_y = bn^*ncn^* + BJ_y$, and hence the above equals $$ \cdots = (bn^*ncn^*+BJ_y)\otimes (nv + J_yV) = $$ $$ = (bn^*+BJ_y)(ncn^*+\H yy)\otimes (nv + J_yV) = $$ $$ = (bn^*+BJ_y)\otimes (ncn^*+\H yy)(nv + J_yV) = $$ $$ = (bn^*+BJ_y)\otimes (ncn^*nv + J_yV) = (bn^*+BJ_y)\otimes (ncv + J_yV), $$ where the last step is proved similarly to the what was done above, based on the fact that $n^*nv-v$ lies in $J_xV$, and hence $nc(n^*nv-v)$ lies in $J_yV$. Being balanced, $\varphi \times \psi $ factors through a well defined map $$ \varphi \otimes \psi :M_x\otimes \Vx \to M_y\tprod y\V y, $$ which is easily seen to be left $B$-linear. In addition, it is easy to see that $\varphi \otimes \psi $ is bijective since its inverse may be obtained by running this proof after interchanging $x$ and $y$, and substituting $n^*$ for $n$. This concludes the proof. \endProof

We have already characterized the annihilator of induced modules in (11.2), but, in case the inducing module is a germ space, there is some more we can say. However we first need a simple technical result.

\state Lemma Given $b$ in $B$, there exists $U$ in $\Nbd $ such that $$ \E xx(b)(v + J_xV) = b\1U v + J_xV, \for v\in V. $$

\Proof Using (5.21) we may write $$ b = c + a_1b_1 + b_2a_2, $$ with $$ c\in \C xx,\quad a_1,a_2\in J_x, \and b_1,b_2\in B. $$ By the definition of $\E xx$, we then have that $$ \E xx(b) = \Psi ^{-1}\big (\q xx(b)\big ) = \Psi ^{-1}\big (\q xx(c)\big ) \explica {(5.19)}{=} \p xx(c) = c + \H xx. $$ Since $a_2$ is a locally constant function vanishing at $x$, we may find $U$ in $\Nbd $ such that $a_2\1U=0$. Given $v$ in $V$, we then observe that $ v-\1Uv\in J_xV, $ so $$ v+J_xV = \1Uv+J_xV, $$ whence $$ \E xx(b)(v + J_xV) = \E xx(b)(\1Uv + J_xV) = $$ $$ = (c + \H xx)(\1Uv + J_xV) = c\1Uv + J_xV = $$ $$ = b\1Uv - a_1b_1\1Uv - b_2a_2\1Uv+ J_x V = b\1Uv + J_x V. \closeProof $$ \endProof

Still considering our fixed $B$-module $V$, pick $x$ in $X$ and, viewing $V$ as a left $\C xx$-module, suppose that we are given a $\C xx$-submodule $W$ of $V$, containing $J_xV$, so that $$ J_xV\subseteq W\subseteq V. $$ The quotient $V/W$ is then clearly a $\C xx$-module annihilated by $\H xx$, and hence also a $\B xx$-module.

The following characterization of the annihilator of $\Indx {V/W}$ is a fundamental technical result. When $W=J_xV$, we have that $V/W=\Vx $, so this result also characterizes the annihilator of $\Indx {\Vx }$.

\state Theorem Working under (5.6), let $V$ be a unital left $B$-module. Pick $x$ in $X$, let $W$ be a $\C xx$-submodule of\/ $V$ containing $J_xV$, and consider $V/W$ as a $\B xx$-module. Then, for every $b$ in $B$, the following are equivalent: \item $b\in \Ann B{\Indx {V/W }}$, \item for every $d$ in $B$, one has that $dbV\subseteq W$.

\Proof Assuming (i), we will initially prove that $$ bV\subseteq W. \reqno{(12.10.1)} $$ Using that $A$ is a regular subalgebra of $B$, write $$ b = \sum _{i=1}^k n_i, $$ where the $n_i$ are normalizers. We then consider the partition of the set $I = \{1,2,\ldots ,k\}$, given by $$ I= I_0 \sqcup \bigsqcup _{y\in \Orb x} I_y, $$ where $$ I_0 = \big \{i\in I: \evx {n_in_i^*} =0\big \}, $$ and $$ I_y = \big \{i\in I: n_i \in \N xy\big \}, \for y\in \Orb x. $$

The reader might have noticed that, so far, the elements represented by the letter $x$ in this paper have mainly been in the \emph {source} of normalizers, but it so happens that, here, it is the membership of $x$ in the \emph {target} of normalizers that is at stake.

We then have that $$ b = \sum _{i\in I_0} n_i + \sum _{y\in \Orb x} \sum _{i\in I_y} n_i = c_0 + \sum _{y\in \Orb x} c_y, $$ where $c_0$ and $c_y$ are implicitly defined above. Incidentally notice that $n_in_i^*\in J_x$, for $i$ in $I_0$, so $$ c_0 = \sum _{i\in I_0} n_i = \sum _{i\in I_0} n_in_i^*n_i \in J_xB, $$ while $$ c_y \in \C xy, \for y\in \Orb x, $$ by (5.28.i). Observing that $c_0V\subseteq J_xV\subseteq W$, in order to verify (12.10.1), all we need to do is prove that $c_yV\subseteq W$, for $y$ in $\Orb x$. Given any $z$ in $\Orb x$, and fixing some $n$ in $\N xz$, we claim that $$ \E xx(bn^*) = \E xx(c_zn^*). \reqno{(12.10.2)} $$ To prove it we compute $$ \E xx(bn^*) = \E xx(c_0n^*) + \sum _{y\in \Orb x} \E xx(c_yn^*), $$ and we want to show that the only surviving term is $\E xx(c_zn^*)$. With respect to the first summand, we have that $$ c_0n^* \in J_xB, $$ so $\E xx(c_0n^*)=0$. Fixing any $y$ in $\Orb x$, with $y\neq z$, observe that for all $i$ in $I_y$, one has that $ n_in^* \in \N xy\N zx, $ and hence $ n_in^* \notin \N xx, $ thanks to (5.12.v). Therefore $\E xx(n_in^*) = 0$, by (5.28.ii) and we deduce that $$ \E xx(c_yn^*) = \sum _{i\in I_y}\E xx(n_in^*) = 0. $$ This concludes the proof of (12.10.2). Still referring to the elements $z$ and $n$ fixed above, observe that $$ c_zn^*\in \C xz \N zx \subseteq \C xz \C zx \subseteq \C xx, $$ so we have that $$ c_zn^* + \H xx = \p xx (c_zn^*) = \E xx(c_zn^*) = \E xx(bn^*) \explica {(11.2)}{\in} \Ann {\B xx}{V/W}. $$ For every $v$ in $V$, we then have that $$ 0 = (c_zn^* + \H xx)(v + W) = c_zn^*v+ W, $$ meaning that $c_zn^*v \in W$. Since this holds for every $v$, it ought to hold also for $nv$, so $$ c_zn^*nv \in W. $$ This gets us very close to our goal, but we still need to deal with the unwanted term $n^*n$ in this expression. For this observe that, since $z\in \src (n)$, we have that $\ev {n^*n}z=1$, so $$ c_z-c_zn^*n \explica {(5.26.ii)}{\in} c_zJ_z \subseteq J_xB, $$ because $c_z\in \C xz$. With this we finally get $$ c_zv = (c_z - c_zn^*n)v + c_zn^*nv \in W, \for v\in V, $$ as desired.

This proves (12.10.1), and, in order to obtain (ii), we just need to realize that $\Ann B{\Indx {V/W }}$ is an ideal in $B$, so $db$ is also in $\Ann B{\Indx {V/W }}$, for every $d$ in $B$. Therefore the first part of the proof, with $db$ in place of $b$, gives (ii).

\itmProof (ii)$\imply $(i) Assuming (ii), we will prove that $$ \E xx(gbh)\in \Ann {\B xx}{V/W }, \for g,h\in B, \reqno{(12.10.3)} $$ and (i) will follow from (11.2). We thus fix $g,h\in B$, and, use (12.9) to find $U$ in $\Nbd $ such that, for every $v$ in $V$, $$ \E xx(gbh)v + J_xV = gbh\1U v + J_xV, $$ and consequently $$ \E xx(gbh)v + W= gbh\1U v + W. $$ The conclusion then follows because $gbh\1U v \in W$, by hypotheses. \endProof

As already mentioned, the following is an immediate consequence, obtained by choosing $W=J_xV$ in (12.10):

\state Corollary Working under (5.6), let $V$ be a unital left $B$-module, and pick $x$ in $X$. Then, for every $b$ in $B$, the following are equivalent: \item $b\in \Ann B{\Indx {\Vx }}$, \item for every $d$ in $B$, one has that $dbV\subseteq J_xV$.

Another easy consequence is as follows:

\state Proposition Under (5.6), let $V$ be a unital left $B$-module. Then $$ \Ann B{V} = \medcap _{x\in X}\Ann B{\Ind x{\Vx }}. $$

\Proof Given $b$ in $\Ann B{V}$, it is obvious that (12.11.ii) holds for every $x$ in $X$, so $b$ is in the annihilator of every $\Ind x{\Vx }$, proving the inclusion ``$\subseteq $'' between the sets in the statement. Conversely, if $b$ lies in the above intersection of annihilators, then, given $v$ in $V$, we have by (12.11) that $$ bv\in J_xV, \for x\in X. $$ In other words, this says that the germ of $bv$ at every $x$ vanishes, so $bv=0$ by (12.7). Since $v$ is arbitrary, we conclude that $b$ lies in $\Ann B{V}$. \endProof

For the case of anihilators of irreducible modules we have the following main result. Its proof is an adaptation to our abstract setting of regular subalgebras of the key idea used in proving \cite [Theorem 7]{SteinbergEffrosHahn}.

\state Lemma Under (5.6), let $V$ be an irreducible left $B$-module. Suppose we are given $x$ in $X$, as well as a $\B xx$-submodule $T \varsubsetneq \Vx $. Then $$ \Ann B {V} = \Ann B {\Indx {\Vx }}= \Ann B {\Indx {\Vx /T }}. $$

\Proof Before we address the statement proper, let us consider the subset of $V$ given by $$ W = \{v\in V : v+J_xV \in T\}. $$ It is elementary to check that $W$ is a $\C xx$-submodule of $V$, containing $J_xV$, incidentally exactly as called for in the statement of (12.10). Moreover, $W$ is clearly the kernel of the composition of quotient maps $$ V \to \frac V{J_xV} = \Vx \to \frac {\Vx }T, $$ so it follows that $$ \frac VW \simeq \frac {\Vx }T, \reqno{(12.13.1)} $$ as $\B xx$-modules. Back to our main goal, let us write $$ I= \Ann B {V}, \quad L= \Ann B {\Indx {\Vx }}, \and M= \Ann B {\Indx {\Vx /T }}, $$ for short, and let us first prove that $M\subseteq I$. So we assume by contradiction that there exists some $b\in M\setminus I$, whence there is some $v$ in $V$, such that $bv\neq 0$. Since $M$ is an ideal in $B$, we see that $$ V_1 := Mv $$ is a submodule of $V$, and since $bv\in V_1$, we have that $V_1$ is nonzero. As $V$ is irreducible, we then deduce that $V_1=V$. Since $T\varsubsetneq \Vx $, by hypothesis, we may choose some $u$ in $V$ whose germ at $x$ is not in $T$, which is the same as saying that $u\in V\setminus W$. Given that $V_1=V$, we may furthermore write $$ u = bv, $$ for some $b$ in $M$. As $b$ annihilates $\Ind x{\Vx /T }$, by assumption, we deduce from (12.13.1) and (12.10) that $$ u= bv \in bV \subseteq W, $$ a contradiction, thus proving that $M\subseteq I$.

In order to prove that $I\subseteq L$, pick any $b$ in $I$. Recalling that $I=\Ann B{V}$, it follows from (12.12) that $b$ lies in $\Ann B{\Ind y{\V y}}$, for every $y$ in $X$, including of course the case $y=x$.

It then remains to prove that $L\subseteq M$. For this pick $b$ in $L$, meaning that $b$ annihilates $\Indx {\Vx }$. So, for every $d$ in $B$, one has that $$ dbV\explica {(12.11)}{\subseteq} J_xV\subseteq W, $$ therefore $b$ annihilates $\Ind x{V/W}$, by (12.10), and the conclusion follows from (12.13.1). \endProof

With this we may generalize several important results in the literature describing ideals in term of isotropy, such as \cite [Theorem 5.23]{DokuchaExel}, \cite [Proposition 5.2.9]{Demeneghi} and \cite [Theorem 7]{SteinbergEffrosHahn}.

\state Theorem Under the conditions of (5.6), we have that: \item Every ideal of $B$ coincides with the intersection of a family of induced ideals\fn {Recall from (11.3) that an induced ideal is the annihilator of an induced module.}. \item Every primitive ideal of $B$ coincides with an induced ideal.

\Proof Given an ideal $I\trianglelefteq B$, and focusing on (i), consider the quotient algebra $$ V:=B/I $$ as a left $B$-module, with the obvious module structure. It is then easy\fn {Using that $A$ is an s-unital subalgebra of $B$.} to see that $\Ann B{V}=I$, so the conclusion follows from (12.12).

Assuming now that $I$ is a primitive ideal, write $I=\Ann B{V}$, for some irreducible left $B$-module $V$. Employing (12.7), there exists some $x$ in $X$ such that $\Vx \neq \{0\}$. The result then follows from (12.13), upon choosing $T= \{0\}$. \endProof

In view of (12.14.ii) and (11.5), it makes sense to ask:

\state Question \rm Under (5.6), is every primitive ideal of $B$ induced by a \underbar {primitive} ideal of some $\B xx$?

The version of this question for Steinberg algebras is also discussed in \cite {SteinbergEffrosHahn}, and a solution is given in \cite [Theorem 8]{SteinbergEffrosHahn} under strong hypotheses. In our setting of regular subalgebras we can give positive answers in two special cases, the first one of which is inspired by \cite [Theorem 8]{SteinbergEffrosHahn}.

\state Proposition Under (5.6), let $I$ be a primitive ideal of $B$, written as $I=\Ann BV$, where $V$ is an irreducible left $B$-module. Suppose that there exists $x$ in $X$, such that $\Vx $ is nonzero and contains a maximal submodule. Then $I$ coincides with the ideal induced by a primitive ideal of $\B xx$.

\Proof If $T$ is a maximal submodule of $\Vx $, then then $$ I = \Ann B {V} \explica {(12.13)}{=} \Ann B {\Indx {\Vx /T }} \explica {(11.4)}{=} \Indx {\Ann {\B xx} {\Vx /T }}, $$ so the result follows because $\Vx /T$ is an irreducible module, and hence its annihilator is a primitive ideal. \endProof

The reader might want to read about the concept of \emph {left max rings}, which forms the basis of the proof of \cite [Theorem 8]{SteinbergEffrosHahn}, precisely by providing maximal submodules, as in the above proof.

Our second answer to a special case of question (12.15) also has strong hypotheses.

\state Proposition Under (5.6), let $I$ be a primitive ideal of $B$, written as $I=\Ann BV$, where $V$ is a nontrivial irreducible left $B$-module. Suppose that there exists an \underbar {isolated} point $x$ in $X$, such that $\Vx $ is nonzero (notice that, by (12.7), this condition is automatic in case $X$ is discrete). Then $I$ coincides with the ideal induced by a primitive ideal of $\B xx$.

\Proof As a first step we claim that $\Vx $ is irreducible as a $\B xx$-module. Arguing by contradiction, let $T\varsubsetneq \Vx $ be a nonzero submodule, and consider $$ W = \big \{v\in V: bv+J_xV\in T,\ \forall b\in B\big \}. $$ It is then obvious that $W$ is a $B$-submodule of $V$, and that $W\neq V$. We will then reach a contradiction by proving that $W$ is nonzero. For this, pick any nonzero element in $T$ and write it as $v+J_xV$, for some $v$ in $V$. As $x$ is isolated, we have that the characteristic function $1_{\{x\}}$ lies in $A$, and it is easy to see that $$ 1_{\{x\}}v + J_xV = v + J_xV, $$ so necessarily $1_{\{x\}}v\neq 0$. We will next prove that $$ n1_{\{x\}}v+J_xV\in T,\for n\in N_B(A). \reqno{(12.17.1)} $$ Assuming first that $\evx {nn^*}=0$, we have that $n = nn^*n \in J_xB$, so (12.17.1) follows. Otherwise we have that $x\in \tgt (n)$, so $n\in \N xy$, where $y=\beta _{n^*}(x)$. Assuming next that $y\neq x$, we have that $$ \evx {n1_{\{x\}}n^*} = \ev {1_{\{x\}}}{\beta _{n^*}(x)} = \ev {1_{\{x\}}}y = 0, $$ so $n1_{\{x\}}n^*\in J_x$, and then $$ n1_{\{x\}}v = nn^*n1_{\{x\}}v = n1_{\{x\}}n^*nv \in J_xV, $$ so again (12.17.1) holds. Finally, assuming that $y=x$, we have that $n\in \C xx$, so $$ n1_{\{x\}}v + J_xV = (n+ \H xx)(1_{\{x\}}v + J_xV) \in T, $$ because $T$ is a $\B xx$-submodule. This proves (12.17.1), and since $B$ is spanned by normalizers, this implies that $$ b1_{\{x\}}v+J_xV\in T, $$ for every $b$ in $B$, meaning that $1_{\{x\}}v$ lies in $W$, and hence that $W$ is nonzero, as claimed. This violates the hypothesis that $V$ is irreducible, so we have proved that $\Vx $ is irreducible. Therefore $$ I = \Ann B {V} \explica {(12.13)}{=} \Ann B {\Indx {\Vx }} \explica {(11.4)}{=} \Indx {\Ann {\B xx} {\Vx }}. $$ so the result follows because $\Vx $ is an irreducible module, and hence its annihilator is a primitive ideal. \endProof

\section Applications to twisted Steinberg Algebras

In this section we intend to show that the theory developed up to now applies to twisted Steinberg algebras. However, before we start moving in the direction we have in mind, let us briefly describe the theory of induced modules studied by Steinberg in \cite {SteinbergOne}, which is the main inspiration for our work, and which we plan to generalize.

As already mentioned, Steinberg only deals with untwisted groupoids, but, once we have put in place the basic ingredients pertaining to twisted groupoids in sections (2--4), above, it is not difficult to develop the first few steps of a theory of induction in the twisted case.

\fix We therefore fix a twisted, ample groupoid $(G,E)$, with bundle projection $\pi :E\to G$, and let $$ A=\Lc , \and B = \Ak . \reqno{(13.1)} $$

Fixing a point $x$ in $G\ex 0$, observe that the restriction of $E$ to the isotropy group $\G xx$, namely $\pi ^{-1}(\G xx)$, is a bona fide line bundle over $\G xx$, so we can speak of the associated twisted Steinberg algebra $$ A_\fld \big (\G xx,\pi ^{-1}(E)\big ), $$ which we shall simply denote by $A_\fld (\G xx,E)$.

Since $\G xx$ is discrete, there is certainly no trouble finding a nowhere vanishing continuous global section, so $\pi ^{-1}(E)$ is necessarily isomorphic to the line bundle $E(\omega )$, where $\omega $ is a 2-cocycle over $\G xx$, thanks to (2.11.i). This said, one clearly has that $$ A_\fld \big (\G xx,E\big ) \simeq \fld (\G xx,\omega ), $$ where $\fld (\G xx,\omega )$ denotes the twisted group algebra. Nevertheless we will continue to emphasize the Steinberg algebra point of view which is perhaps more in line with our current standpoint.

The inducing bimodule $KL_x$ described in \cite [Proposition 7.8]{SteinbergOne} in itself is not relevant here, but one doesn't need much imagination in order to reformulate it in the present context, and the reader might have already guessed that the correct replacement is the set of all finitely supported sections of the bundle restricted to $G_x$, which we will denote by $L_c(G_x, E)$. We may then make $L_c(G_x, E)$ into an $\Ak $-$A_\fld \big (\G xx,E\big )$-bimodule as follows: the left module structure is defined by $$ (f \xi )(\gamma ) = \sum _{\alpha \beta =\gamma } f(\alpha )\xi (\beta ), \for \gamma \in G_x, $$ for $f\in \Ak $, and $\xi \in L_c(G_x, E)$, while the right module structure is given by $$ (\xi g)(\gamma ) = \sum _{\buildrel \scriptstyle \alpha \beta =\gamma \over {\beta \in \G xx}} \xi (\alpha )g(\beta ), \for \gamma \in G_x, $$ for $\xi \in L_c(G_x, E)$, and $g\in A_\fld (\G xx, E)$. The products within the above sums are, of course, occurrences of the bundle multiplication, rather than multiplication of scalars, but, apart from that, these are exactly the same formulas as in the standard (untwisted) case.

In particular, if we are given a left module $V$ over $A_\fld (\G xx, E)$, then $$ L_c(G_x, E)\mathop {\medotimes }\limits _{A_\fld (\G xx, E)} V, $$ is the \emph {induced module} over $\Ak $.

The reader may wish to verify any details eventually left out above, but these will nevertheless come as a consequence of what we are about to do.

Turning now to the main goal of this section, the following result is intended to fit twisted Steinberg algebras within our theory.

\state Proposition Given a twisted, ample groupoid $(G,E)$, the following holds: \iStyle i \item $\Lc $ is an abelian algebra generated by its idempotent elements. \item $\Lc $ admits a canonical embedding as a subalgebra of $\Ak $. \item For every continuous local section $\xi $ of $E$, whose domain is a compact open bisection, one has that $\tilde \xi $ is a normalizer of $\Lc $ in $\Ak $. \item If $\xi $ is as above, then the local section $\xi \st $ defined on $\dom (\xi )^{-1}$ by $$ \xi \st (\gamma ) = \clauses { \cl {\xi (\gamma ^{-1})^{-1}}{\text {if } \xi (\gamma ^{-1})\neq 0_{\gamma ^{-1}}, } \cl {0_{\gamma }}{\text {otherwise},} } $$ is continuous, and the partial inverse of $\tilde \xi $ is given by $\widetilde {\xi \st }$. \item $\Lc $ is a regular subalgebra of $\Ak $.

\Proof We leave (i) as an easy exercise.

Regarding (ii), recall from (2.11.ii) that the restiction of $E$ to $G\ex 0$ is a trivial line bundle with a distinguished global section, namely $x\mapsto 1_x$. Therefore one may view any given $f\in \Lc $ as a section of $E$ over $G\ex 0$, namely $$ \tilde f:x\in G\ex 0 \mapsto f(x)1_x\in E. $$ This said it is easy to see that the correspondence $f\mapsto \tilde f$ gives an embedding of $\Lc $ into $\Ak $.

Given $\xi $ as in (iii), and using (2.12.iii), we have that $\xi \st $ is continuous on $\supp (\xi )^{-1}$, and it is clear that it is also continuous on $\dom (\xi )^{-1}\setminus \supp (\xi )^{-1}$. Recalling that $\supp (\xi )$ is both open and closed relative to $\dom (\xi )$ by (2.3.v), we see that $\xi \st $ is continuous everywhere on its domain. We then leave it for the reader to do the easy computations needed to verify that $\tilde \xi $ is indeed a normalizer, with partial inverse $\widetilde {\xi \st }$.

This said, it is clear that $\Ak $ is spanned by the set of all normalizers, so in order to prove (v), we just need to show that $\Lc $ is an s-unital subalgebra of $\Ak $. For this we fix a finite subset $F$ of $\Ak $, and we will provide a dedicated unit for $F$, belonging to $\Lc $. Assuming, as we may, that $$ F = \{\tilde \xi _1, \tilde \xi _2, \ldots , \tilde \xi _k\}, $$ where the $\xi _i$ are continuous local sections defined on compact open bisections $U_i$, let $$ K = \bigcup _{i=1}^k s(U_i)\cup r(U_i), $$ so that $K$ is a compact open subset of $G\ex 0$, and hence the characteristic function $\1K$ belongs to $\Lc $. An easy computation shows that $\1K\tilde \xi _i= \tilde \xi _i= \tilde \xi _i\1K$, so we see that $\1K$ is a dedicated unit for $F$. This concludes the proof. \endProof With the above result we see that the theory of induced representations developed in sections (5-10) may be applied in the context of Steinberg algebras. However so far it might not be clear that the induction process described above relates to other induction theories, such as e.g.~the one developed by Steinberg in \cite {SteinbergOne} for the case of untwisted groupoids.

In what follows we will therefore discuss the main players arising from our theory when applied to the twisted groupoid situation, with emphasis on describing the isotropy algebra and the imprimitivity bimodule. The goal is to show that these reduce to the twisted group algebra of the relevant isotropy group and to the standard imprimitivity bimodule, respectively.

Our next immediate step will be to study the partial homeomorphisms $\beta _n$, introduced in (5.10), for the normalizers mentioned in (13.2.iii).

\state Proposition Let $\xi $ be a continuous local section of $E$, whose domain is a compact open bisection, and let $n=\tilde \xi $. Then, viewing $n$ as a normalizer, one has that: \item $\src (n)= s(\supp (\xi ))$, and $\tgt (n) = r(\supp (\xi ))$, \item $\beta _n\big (s(\gamma )\big ) = r(\gamma ), \for \gamma \in \supp (\xi )$, \item For any $x, y\in G\ex 0$, one has that $n\in \N yx$ if and only if $\supp (\xi )\cap \Suber Gyx$ is nonempty.

\Proof Replacing $\xi $ by the restriction of $\xi $ to its support (which is again a compact open bisection by (2.3.v)), nothing else is affected in the statement, so we may assume that $\xi $ vanishes nowhere on its domain. In particular $$ \xi \st (\gamma ) = \xi (\gamma ^{-1})^{-1}, $$ for every $\gamma $ in $\dom (\xi )^{-1}$. An easy computation then shows that $n^*n$ coincides with the characteristic function on $s(\supp (\xi ))$, and that $nn^*$ is the characteristic function on $r(\supp (\xi ))$, from where (i) follows.

For any $a$ in $\Lc $ we then have that $$ (n^*an) \big (s(\gamma )\big ) = \xi \st (\gamma ^{-1}) a\big (r(\gamma )\big )\xi (\gamma ) = \xi (\gamma )^{-1} a\big (r(\gamma )\big )\xi (\gamma ) =a\big (r(\gamma )\big ), $$ proving (ii). As for the last point, we have by definition that $n\in \N yx$ if and only if $x\in \src (n)$, and $\beta _n(x)=y$, which is in turn equivalent to the condition in (iii). \endProof

Our next goal will be to identify the isotropy algebra $\B xx$ defined in (5.25), as well as the imprimitivity module $M_x$ introduced in (6.2).

For this, we will first present the following technical result that could be thought of as a replacement for \cite [Lemma 3.6]{LisaOne}, a result that requires the groupoid to be Hausdorff, and hence cannot be applied to our present situation.

\state Lemma Let $Z\subseteq G$ be any subset, and let $\{\xi _i\}_{1\leq i\leq k}$ be a finite set of continuous local sections whose domains are compact open bisections, and such that the function $$ b := \sum _{i=1}^k\tilde \xi _i $$ vanishes on $Z$. Suppose also that $Z\cap \dom (\xi _i)$ has at most one point for every $i$. Then there is a finite set of continuous local sections $\{\eta _j\}_{1\leq _j\leq l}$ whose domains are compact open bisections, such that $$ b := \sum _{j=1}^l\tilde \eta _j, $$ and each $\tilde \eta _j$ vanishes on $Z$.

\Proof We begin by splitting $$ \{1, 2, \ldots , k\} = I_0\sqcup I_1, $$ where a given $i$ lies in $I_0$ if and only if $\dom (\xi _i)\cap Z$ is empty. For each $i\in I_1$, we then let $\gamma _i$ be the single element of $\dom (\xi _i)\cap Z$, and we put $$ \Gamma = \{\gamma _i : i\in I_1\}. $$ Observing that there is no reason for the correspondence $i\mapsto \gamma _i$ to be injective, we let $$ J_\gamma = \{i\in I_1: \gamma _i=\gamma \}, \for \gamma \in \Gamma . $$

As a first case, let us assume that $I_0$ is empty and that $\Gamma $ is a singleton. In other words, this is to say that there is some $\gamma $ in $G$, such that $\dom (\xi _i)\cap Z=\{\gamma \}$, for all $i$. We then put $$ U = \medcap _{i=1}^k \dom (\xi _i), $$ noticing that $U$ is an open\fn {The reader should not be fooled into believing that $U$ is compact, as the intersection of compact sets might not be compact in a non-Hausdorff space.} bisection containing $\gamma $. Since open bisections in an \'etale groupoid are always locally compact, Hausdorff, totally disconnected topological spaces, we may find a compact open neighborhood $W$ of $\gamma $ contained in $U$, and hence also contained in each $\dom (\xi _i)$. For every $i\leq k$, we then define \Item {$\bullet $} $\eta _i=\xi _i|_{W}$, \Item {$\bullet $} $V_i=\dom (\xi _i)\setminus W$, \Item {$\bullet $} $\zeta _i=\xi _i|_{V_i}$.

\medskip \noindent It is then clear that $W$, as well as each $V_i$ are compact open bisections, that each $\eta _i$ and each $\zeta _i$ are continuous local sections, and finally that $$ \tilde \xi _i = \tilde \eta _i + \tilde \zeta _i, $$ for all $i$. We conclude that $$ b = \sum _{i=1}^k \tilde \eta _i + \sum _{i=1}^k \tilde \zeta _i = \tilde \eta + \sum _{i=1}^k \tilde \zeta _i, $$ where $\eta =\sum _{i=1}^k \eta _i$. Of course we are benefiting from the fact that the $\eta _i$ are local sections defined on the same domain, namely $W$, and hence their pointwise sum defines a local section on $W$.

Notice that the domains of the $\zeta _i$, namely the $V_i$, contain no elements from $Z$, so clearly $\tilde \zeta _i$ vanish on $Z$. On the other hand, the domain of $\eta $ contains a single element of $Z$, namely $\gamma $, so $\tilde \eta $ vanishes on $Z\setminus \{\gamma \}$. Moreover, $$ 0 = b(\gamma ) = \tilde \eta (\gamma ) + \sum _{i=1}^k \tilde \zeta _i(\gamma ) = \tilde \eta (\gamma ), $$ so we see that $\tilde \eta $ actually vanishes on all of $Z$, and we have therefore completed the proof of this case.

In the general case, write $$ b = \sum _{i\in I_0}\tilde \xi _i + \sum _{\gamma \in \Gamma }\sum _{j\in J_\gamma }\tilde \xi _j = \sum _{i\in I_0}\tilde \xi _i + \sum _{\gamma \in \Gamma }b_\gamma , $$ where each $b_\gamma =\sum _{j\in J_\gamma }\tilde \xi _j$.

For $i$ in $I_0$, notice that $\tilde \xi _i$ vanishes on $Z$, because $\dom (\xi _i)\cap Z$ is empty. On the other hand, for $\gamma $ in $\Gamma $, and for $j$ in $J_\gamma $, we have that $\tilde \xi _j$ vanishes on $Z\setminus \{\gamma \}$, because $\dom (\xi _j)$ has empty intersection with $Z\setminus \{\gamma \}$, so it follows that $b_\gamma $ also vanishes on $Z\setminus \{\gamma \}$. We claim that $b_\gamma $ in fact vanishes on all of $Z$, because $$ 0 = b(\gamma ) = \sum _{i\in I_0}\tilde \xi _i(\gamma ) + \sum _{\alpha \in \Gamma }b_\alpha (\gamma ) = b_\gamma (\gamma ). $$ The result then follows upon applying the case already proved to each $b_\gamma $. \endProof

The following depends on the above Lemma in a crucial way and is intended to identify important subspaces of $B$ related to the construction of $\B xx$ and $M_x$.

\state Proposition With $B$ as in (13.1), fix $x$ in $G\ex 0$, and $b$ in $B$. One then has that \item $b$ vanishes on $G_x$ if and only if $b \in BJ_x$. \item $b$ vanishes on $\G xx$ if and only if $b \in \L xx \explica {(5.25.1)}{=} J_xB + BJ_x$.

\Proof (i) By (5.26.i) we have that $J_x$ is s-unital. So, given any element $b$ in $BJ_x$, we may write $b = bv$, for some $v$ in $J_x$, thanks to (5.14.ii). Given $\gamma $ in $G_x$, we then have that $$ b(\gamma ) = (bv)(\gamma ) = b(\gamma )v(s(\gamma )) = b(\gamma )v(x) = 0, $$ so $b$ vanishes on $G_x$, as required. Conversely, if $b$ vanishes on $G_x$, write $$ b=\sum _{i=1}^k\tilde \xi _i, \reqno{(13.5.1)} $$ where the $\xi _i$ are continuous local sections whose domains are compact open bisections. Therefore, for each $i$, we have that $G_x\cap \dom (\xi _i)$ has at most one point, allowing us to invoke (13.4), and hence assume that each $\tilde \xi _i$ vanishes on $G_x$. Consequently, recalling that $\tilde \xi _i$ is a normalizer, we have that $$ \big ((\tilde \xi _i)^*\tilde \xi _i\big )(x) = \sum _{\gamma \in G_x}(\tilde \xi _i)^*(\gamma ^{-1})\tilde \xi _i(\gamma ) = 0. $$ This says that $x\notin \src (\tilde \xi _i)$, so $\tilde \xi _i$ belongs to $BJ_x$, by (5.27.i), and since $i$ is arbitrary, we also have that $b\in J_x$, concluding the proof of (i).

\medskip \noindent (ii) Assuming that $b \in J_xB + BJ_x$, we may use (5.18.i) to find $u,v\in J_x$, such that $b = ub +bv - ubv$, so if $\gamma \in \G xx$, we have $$ b(\gamma ) = (ub)(\gamma ) +(bv)(\gamma ) - (ubv)(\gamma ) = $$ $$ = u\big (r(\gamma )\big )b(\gamma ) +b(\gamma ) v\big (s(\gamma )\big ) - u\big (r(\gamma )\big )b(\gamma )v\big (s(\gamma )\big ) = $$ $$ = u(x)b(\gamma ) +b(\gamma ) v(x) - u(x)b(\gamma )v(x) = 0, $$ so $b$ vanishes on $\G xx$. Conversely, assuming that $b$ vanishes on $\G xx$, write $$ b=\sum _{i=1}^k\tilde \xi _i, \reqno{(13.5.2)} $$ where the $\xi _i$ are continuous local sections whose domains are compact open bisections. Therefore, for each $i$, we have that $\G xx\cap \dom (\xi _i)$ has at most one point, allowing us to invoke (13.4), and hence assume that each $\tilde \xi _i$ vanishes on $\G xx$. Consequently $\supp (\xi _i)\cap \G xx$ is empty and, recalling that $\tilde \xi _i$ is a normalizer, we have by (13.3.iii) that $\tilde \xi _i\notin \N xx$. Thus, either $x\notin \src (\tilde \xi _i)$, in which case $$ \tilde \xi _i \explica {(5.27.i)}{\in} BJ_x \subseteq J_xB + BJ_x, $$ or $x\in \src (\tilde \xi _i)$, but $\beta _{\tilde \xi _i}(x)\neq x$, in which case $$ \tilde \xi _i \explica {(5.27.iii)}{\in} J_xB + BJ_x. $$ This completes the proof. \endProof

With this we may prove another of our main results:

\state Theorem Let $(G,E)$ be a twisted, ample groupoid, and put $A=\Lc $ and $B = \Ak $. Fixing $x$ in $G\ex 0$, we have that: \item $M_x$ is naturally isomorphic to $L_c(G_x,E)$, as left $B$-modules. \item The isotropy algebra $\B xx$ is naturally isomorphic to the twisted group algebra of the isotropy group $\G xx$, namely $A_\fld (\G xx, E)$. \item Identifying $\B xx$ with $A_\fld (\G xx, E)$, according to (ii), the left $B$-module isomorphism of (i) is right $\B xx$-linear. Therefore $M_x$ is isomorphic to $L_c(G_x,E)$, as $B$-$\B xx$-bimodules.

\Proof (i) Observe that the restriction map $$ \psi : b\in B\mapsto b|_{G_x} \in L_c(G_x,E) $$ is left $B$-linear and we claim that it is also surjective. Indeed, since $G_x$ is discrete, we have that $L_c(G_x,E)$ is linearly spanned by sections whose support is a singleton, so it suffices to show that these are in the range of $\psi $. Given such a section $\eta $, write its support as $\{\gamma \}$, and let $u=\eta (\gamma )$. By (2.3.ii) we can find a local section $\xi $ of $E$ passing through $u$, and we may clearly suppose that $\dom (\xi )$ is a compact open bisection containing $\gamma $. Therefore $\psi (\tilde \xi )=\eta $, proving surjectivity.

By (13.5.i) the null space of $\psi $ is $BJ_x$, so $\psi $ factors through $M_x = B/BJ_x$, providing the desired isomorphism.

\itmProof (ii) By (13.5.ii) the null space of the surjective function $$ \varphi :b\in B\mapsto b|_{\G xx} \in A_\fld (\G xx, E) $$ is $\L xx$, so it factors through $\B xx \explica {(5.21)}{\simeq} B/\L xx$, providing a $\fld $-linear isomorphism $$ \tilde \varphi : \B xx \to A_\fld (\G xx, E). $$ It therefore suffices to prove that $$ \varphi (b_1b_2) = \varphi (b_1)\varphi (b_2), \reqno{(13.6.1)} $$ for all $b_1,b_2\in \B xx$.

Observe that $B$ is spanned by the set $N$ of all normalizers of the form $\tilde \xi $, where $\xi $ is a local section whose domain is a compact open bisection. Therefore, according to (5.28.iii), $\B xx$ is spanned by the set of all $\tilde \xi +\L xx$, with $\tilde \xi $, as above, for which $\tilde \xi $ lies in $\N xx$. In proving (13.6.1), we may then suppose that $$ b_1=\tilde \xi _1 +\L xx, \and b_2=\tilde \xi _2+\L xx, $$ where $\xi _1$ and $\xi _2$ are as above. Incidentally, notice that, by (13.3.iii), each $\tilde \xi _i$ lies in $\N xx$ if and only if there exists some element $\gamma _i\in \supp (\xi )\cap \G xx$, which is clearly unique, given that $\supp (\xi _i)$ is a bisection.

The restriction of each $\tilde \xi _i$ to $\G xx$ is thus the function supported on the singleton $\{\gamma _i\}$, taking on the value $\xi _i(\gamma _i)$ at $\gamma _i$. With this, the verification of (13.6.1) becomes elementary.

Proving the last point amounts to showing that $$ \psi (fg) = \psi (f)\varphi (g), \reqno{(13.6.2)} $$ where $f\in M_x$, and $g\in \B xx$, and we may clearly assume that $$ f=\tilde \xi + BJ_x, \and g=\tilde \eta + L_x, $$ where $\xi $ and $\eta $ are local sections whose domains are compact open bisections, and moreover $\supp (\eta )\cap \G xx$ is nonempty, hence necessarily a singleton. The left-hand-side of (13.6.2) then becomes $$ \psi (fg) = \psi ( \widetilde {\xi \eta } + BJ_x) =\widetilde {\xi \eta }|_{G_x}, $$ while the right-hand-side is $$ \psi (f)\varphi (g) =\tilde \xi |_{G_x} \tilde \eta |_{\G xx}, $$ from where (13.6.2) follows with no difficulty. \endProof

Now that we know that the inclusion ``$\Lc \subseteq \Ak $'' is regular, and now that we have identified the main ingredients of the induction process, we may apply all of the results proved for regular inclusions, above, to twisted Steinberg algebras. The following gives a sample of this procedure:

\state Theorem Let $(G,E)$ be a twisted, ample groupoid and let $V$ be an irreducible left module over $\Ak $. Assuming that $\Resx V$ is nonzero for some $x$ in $G\ex 0$, e.g.~when $V$ is finite dimensional as a $\fld $-vector space (c.f. (7.6)), one has that $\Resx V$ is an irreducible left module over the twisted group algebra of the isotropy group $\G xx$, namely $A_\fld (\G xx, E)$, and $V$ is naturally isomorphic to the module induced by $\Resx V$.

\Proof This is just a special case of (10.2). \endProof

\state Theorem Let $(G,E)$ be a twisted, ample groupoid, and let $x\in G\ex 0$. Given any unital left $A_\fld (\G xx,E)$-module $W$, let $V$ be the $\Ak $-module induced by $W$. Then: \item $W$ is naturally isomorphic to $\Resx V$, \item Every submodule $V_1\subseteq V$ is induced by a unique submodule $W_1\subseteq W$.

\Proof The first point is a special case of (8.4), while the second one follows from (9.3). \endProof

Regarding our theory of ideals, here is the application of (12.14) for twisted Steinberg algebras:

\state Theorem Let $(G,E)$ be a twisted, ample groupoid. Then \item Every ideal of $\Ak $ coincides with the intersection of a family of annihilators of induced modules. \item Every primitive ideal of $\Ak $ coincides with the annihilator of some induced module.

\relax

\overfullrule 0pt

\references

\Article AraBosa P. Ara, J. Bosa, R. Hazrat, A. Sims; Reconstruction of graded groupoids from graded Steinberg algebras; Forum Math., 29 (2017), no. 5, 1023--1037

\Bibitem BeckyOne B. Armstrong, L. O. Clark, K. Courtney, Y. Lin, K. Mccormick, J. Ramagge; Twisted Steinberg algebras; \sl J. Pure Appl. Algebra, \bf 226 \rm (2022), no. 3, Paper No. 106853, 33 pp

\Article BeckyTwo B. Armstrong, G. G. De Castro, L. O. Clark, K. Courtney, Y. Lin, K. Mccormick, J. Ramagge, A. Sims, B. Steinberg; Reconstruction of twisted Steinberg algebras; Int. Math. Res. Not., 2023 (2023), no. 3, 2474-2542

\Article LisaOne L. O. Clark, C. Farthing, A. Sims, M. Tomforde; A groupoid generalisation of Leavitt path algebras; Semigroup Forum, 89 (2014), no. 3, 501-517

\Article Demeneghi P. Demeneghi; The ideal structure of Steinberg algebras; Adv. Math., 352 (2019), 777-835

\Article DokuchaExel M. Dokuchaev and R. Exel; The ideal structure of algebraic partial crossed products; Proc. London Math. Soc., 115 (2017), 91-134

\Article EffrosHahnOne E. G. Effros and F. Hahn; Locally compact transformation groups and C*-algebras; Bull. Amer. Math. Soc., 73 (1967), 222--226

\Bibitem EffrosHahnTwo E. G. Effros and F. Hahn; Locally Compact Transformation Groups and C*-Algebras; Mem. Amer. Math. Soc. 75, American Mathematical Society, Providence, RI, 1967

\Article Tight R. Exel; Inverse semigroups and combinatorial C*-algebras; Bull. Braz. Math. Soc. (N.S.), 39 (2008), 191-313

\Article ExelPardo R. Exel and E. Pardo; The tight groupoid of an inverse semigroup; Semigroup Forum, 92 (2016), no. 1, 274-303.

\Bibitem ExelPitts R. Exel and D. R. Pitts; Characterizing groupoid C*-algebras of non-Hausdorff \'etale grou\-poids; Lecture Notes in Mathematics, vol. 2306, Springer, 2022

\Article FS T. Fack, G. Skandalis; Sur les repr\'esentations et id\'eaux de la C*-alg\`ebre d'un feuilletage; J. Oper. Theory, 8 (1983), 95-129

\Article FeldmanMoore J. Feldman and C. Moore; Ergodic equivalence relations, cohomology, and von Neumann algebras. II; Trans. Amer. Math. Soc., 234 (1977), no. 2, 325-359

\Article GR E. C. Gootman and J. Rosenberg; The structure of crossed product C*-algebras: a proof of the generalized Effros-Hahn conjecture; Invent. Math., 52 (1979), no. 3, 283-298

\Article IonWill M. Ionescu and D. Williams; The generalized Effros-Hahn conjecture for groupoids; Indiana Univ. Math. J., 58 (2009), no. 6, 2489-2508

\Bibitem Jacobson N. Jacobson; Structure of Rings; volume 37 of Colloquium Publications. American Mathematical Society, New Haven, CT, 1956

\Article Keimel K. Keimel; Alg\`ebres commutatives engendr\'ees par leurs \'El\'ements idempotents; Can. J. Math., 22 (1979), 1071-1078

\Article Kumjian A. Kumjian; On C*-diagonals; Canad. J. Math., 38 (1986), 969-1008

\Bibitem Lawson M. V. Lawson; Inverse semigroups, the theory of partial symmetries; World Scientific, 1998

\Article Nystedt P. Nystedt; A survey of s-unital and locally unital rings; Revista Integraci\'on - Universidad Industrial de Santander, 37 (2019), 251-260

\Bibitem Paterson A. L. T. Paterson; Groupoids, inverse semigroups, and their operator algebras; Birkhauser, 1999

\Bibitem Renault J. Renault; A groupoid approach to C*-algebras; Lecture Notes in Mathematics vol.~793, Springer, 1980

\Article RenaultStruId J. Renault; The ideal structure of groupoid crossed product C*-algebras. With an appendix by Georges Skandalis; J. Operator Theory, 25 (1991), no. 1, 3-36

\Article RenaultCartan J. Renault; Cartan subalgebras in C*-algebras; Irish Math. Soc. Bulletin, 61 (2008), 29-63

\Article Sauvageot J. L. Sauvageot; Ideaux primitifs de certains produits croises; Math. Ann., 231 (1977), 61-76

\Bibitem SimsWill A. Sims and D. Williams; The primitive ideals of some \'etale groupoid C*-algebras; preprint, arXiv:1501.02302 [math.OA], 2015

\Article SteinbergOne B. Steinberg; A groupoid approach to discrete inverse semigroup algebras; Adv. Math., 223 (2010), no. 2, 689-727

\Article SteinbergDisintegr B. Steinberg; Modules over \'etale groupoid algebras as sheaves; J. Aust. Math. Soc., 97 (2014), 418-429

\Article SteinbergEffrosHahn B. Steinberg; Ideals of \'etale groupoid algebras and Exel's Effros-Hahn conjecture; J. Noncommut. Geom., 15 (2021), n.o 3, 829-839

\endgroup

\vskip 1cm

\def \Address #1#2#3{{\bigskip {\tensc #2} \par \it E-mail address: \tt #3 \par }}

\Address {M. Dokuchaev} {Departamento de Matem\'atica, Universidade de S\~ao Paulo, Rua do Mat\~ao, 1010, S\~ao Paulo, SP, 05508-090, Brazil} {dokucha@ime.usp.br}

\Address {R. Exel} {Departamento de Matem\'atica, Centro de Ci\^encias F\'{\i }sicas e Mate\-m\'a\-ticas, Universidade Federal de Santa Catarina, Florian\'opolis, SC, 88040-900, Brazil} {ruyexel@gmail.com}

\Address {H. Pinedo} {Escuela de Matem\'aticas, Universidad Industrial de Santander, Bucaramanga, Santander, Colombia} {hpinedot@uis.edu.co}

\close

\bye